\documentclass[12pt%
,draft%
]{article}
\usepackage{amsmath,amsfonts,amsthm,amssymb}
\usepackage{verbatim}
\newif\ifnotation  
\newif\iflabelsinmargin 
\newcommand\iwlog{{\it wlog\/}}

\newcommand\ga{\alpha}
\newcommand\gb{\beta}
\newcommand\gd{\delta}

\renewcommand\gg{\gamma}

\newcommand\gi{\iota}
\newcommand\gk{\kappa}
\newcommand\gl{\lambda}
\newcommand\gs{\sigma}
\newcommand\gS{\Sigma}
\newcommand\gth{\theta}
\newcommand\gw{\omega}

\newcommand\cb{{\mathcal B}}

\newcommand\cd{{\mathcal D}}
\newcommand\ce{{\mathcal E}}
\newcommand\cf{{\mathcal F}}

\newcommand\cl{{\mathcal L}}
\newcommand\cm{{\mathcal M}}

\newcommand\catx[2]{\mskip-2mu\raise#1\hbox{$#2\smallfrown$}\mskip-2mu}
\newcommand\cat{%
    \mathchoice{\catx{5pt}{\displaystyle}}{\catx{4.5pt}{\textstyle}}%
		{\catx{2.5pt}{\scriptstyle}}{\catx{2pt}{\scriptscriptstyle}}}

\newcommand\nothing{\varnothing}
\renewcommand\setminus{\smallsetminus}
\newcommand\restrict{{\restriction}} 

\newcommand\tri{\vartriangleleft}
\newcommand\trige{\trianglerighteq}

\newcommand\sat{\models}
\newcommand\ps{{\mathcal P}}

\newcommand\sh{^\sharp}
\newcommand\card[1]{|#1|}

\newcommand\set[1]{\{\,#1\,\}}
\newcommand\seq[1]{(\,#1\,)}

\newcommand\image{\raise1.5pt\hbox{\lq\lq}\kern-.7pt}

\newcommand\union{\bigcup}
 
\newcommand\cof{\operatorname{cf}}
\newcommand\range{\operatorname{range}}
\newcommand\domain{\operatorname{domain}}
\renewcommand\sup{\operatorname{sup}}
\newcommand\lub{\operatorname{lub}}
\newcommand\ult{\operatorname{ult}} 
\newcommand\crit{\operatorname{crit}}
\newcommand\len{\operatorname{len}}
\newcommand\inx{\operatorname{index}}

\newcommand\lt{{<}}
\newcommand\comp{{\circ}}

\newtheorem{theorem}{Theorem}[section]
\newtheorem{cor}[theorem]{Corollary}
\newtheorem{lemma}[theorem]{Lemma}
\newtheorem{prop}[theorem]{Proposition}
\newtheorem*{claim}{Claim}
\theoremstyle{definition}
\newtheorem{defin}[theorem]{Definition}
\newtheorem*{unnumthm}{Theorem}
\theoremstyle{remark}

\newtheorem{case}{Case}
\newcommand\startcases{\setcounter{case}0}
\newtheorem{quest}{Question}

\def\lineto#1{%
	\def\xxx#1{height .4pt depth 0pt width #1}
	\expandafter\vrule\xxx#1
}

\newcommand{\paranumber}[1]{\marginpar{\thesection.#1}}

\newcommand\notate[1]{}
\newcommand\usage[1]{}

\newcommand\eqdef{\stackrel{\scriptscriptstyle\text{def}}=}
\newcommand\indiscbelong{\operatorname{\gb}}

\newcommand\OO{\operatorname{O}} 
\def\SO(#1){\sup(\OO(#1))}	
\newcommand\PO{\OO'}		
\def\SPO(#1){\sup(\PO(#1))}
\newcommand\ot{\operatorname{otp}}
\newcommand\tcf{\operatorname{\text{\rm tcf}}}
\newcommand\pcs{precovering set}
\newcommand\pri{principal indiscernible}
\newcommand\is{indiscernible sequence}
\newcommand\bis{basic \is}
\newcommand\pris{principal \is}

\newcommand\core{K}

\newcommand\mse{\mathfrak}
\newcommand\barh{\bar h}	
\newcommand\barmse[1]{\overline{\mse #1}}

\newcommand\bari{\bar\imath}
\newcommand\bark{{\overline\core}}
\newcommand\bare{\overline E}
\newcommand\barf{\bar F}
\newcommand\barm{\overline M}

\newcommand\lis{\ell}

\newcommand\lessb{<_b}
\newcommand\eqb{=_b}
\newcommand\geb{\ge_b}
\newcommand\leb{\le_b}
\newcommand\trieq{\trianglelefteq}
\newcommand\gsless{\prec}
\newcommand\gslesseq{\preccurlyeq}
\newcommand\tand{\text{ and }}

\renewcommand\vec[1]{\expandafter\boldsymbol #1}
\begin{document}
%

\date{\today}
\title{Indiscernible Sequences for Extenders, and\\ 
The Singular Cardinal Hypothesis}
\author{Moti Gitik
\thanks{ 
Some of the results were obtained while Gitik was visiting 
Los Angeles in Fall 1991.  He would like to thank A.~Kechris,
D.~Martin and J.~Steel for their hospitality.}
\and 
William J. Mitchell
\thanks{Mitchell was
partially supported by grant number DMS-9240606 from the National
Science Foundation.}}

\maketitle
\begin{abstract}
We prove several results giving lower bounds for the large cardinal
strength of a failure of the singular cardinal hypothesis.  The main
result is the following theorem:
\begin{unnumthm}
Suppose $\gk$ is a singular strong limit cardinal
and $2^\gk\ge\gl$ where $\gl$ is not the successor of a cardinal of
cofinality at most $\gk$.  If $\cof(\gk)>\gw$ then it follows that
$o(\gk)\ge\gl$, and if $\cof(\gk)=\gw$ then either $o(\gk)\ge\gl$ or
$\set{\ga:K\sat
o(\ga)\ge\ga^{+n}}$ is cofinal in $\gk$ for each $n\in\gw$.
\end{unnumthm}
We also prove several results which extend or are related to this
result, notably
\begin{unnumthm}
If 
$2^{\gw}<\aleph_{\gw}$  and  $2^{\aleph_\gw}>\aleph_{\gw_1}$
then there is a sharp for a model with a strong cardinal.
\end{unnumthm}

In order to prove these theorems we give a detailed analysis of
the sequences of indiscernibles which come from applying
the covering lemma to
nonoverlapping sequences of extenders.
\end{abstract}


\write16{Do for final version:}
\write16{******** Remove the note at the beginning.}
\write16{******** change definition of\string\paranumber\space to {\par}.}

%

%
\renewcommand\thetheorem{\arabic{theorem}}
\paranumber{1}
The covering lemma asserts, roughly, that for any uncountable set $x$
of ordinals there is a set $y\supset x$  such that $\card y=\card x$
and $y\in K[C]$ where $C$ is some sequence of indiscernibles. 
In many  applications of the covering lemma, such as in the proof that
$\gl^+=(\gl^+)^K$ whenever $\gl$ is a singular cardinal, 
the indiscernibles don't
pose a problem: the covering lemma is used in an interval
where there are no measurable cardinals in $K$, and
thus 
there are no indiscernibles.
For other applications, such as the singular cardinal hypothesis, this
is not possible.  If $\gk$ is singular then the covering lemma
asserts, in effect, that the number of subsets of $\gk$ is determined
by the number of cofinal sequences of indiscernibles in $\gk$. Thus
giving an upper limit to $2^{\gk}$ entails giving an upper limit on
the number of sequences of indiscernibles, which requires a detailed
understanding of these sequences.
It is this second class of applications which we will be considering
in this paper.

Work on this class of problems began with the work of Dodd and
Jensen on the model $L[\mu]$.  Their ideas were
extended to models for sequences of measurable cardinals by Mitchell
\cite{mitchell.applications,mitchell.defsch} and Gitik
\cite{gitik.mble-cards-not-ch}.  In this paper we extend this analysis to
models containing
sequences of nonoverlapping extenders, including models up to a strong
cardinal.
Our main application is the following theorem:

\paranumber{2}
\begin{theorem}[\ref{app-main-thm}]
Suppose that $\gk$ is a strong limit cardinal with 
$\cof(\gk)=\gd<\gk$, and that $2^\gk\ge\gl>\gk^+$, where $\gl$ is not
the successor of a cardinal of cofinality less than~$\gk$.
\begin{enumerate}
\item If $\gd>\gw_1$ then $o(\gk)\ge\gl+\gd$.
\item If $\gd=\gw_1$ then  $o(\gk)\ge\gl$.
\item If $\gd=\gw$ then either $o(\gk)\ge\gl$ 
or else $\set{\ga:K\sat o(\ga)\ge\ga^{+n}}$ is cofinal
in $\gk$ for each $n<\gw$.
\end{enumerate}
\end{theorem}

\paranumber{3}
Woodin \cite{cummings.GCH-successors} has constructed models
of $2^\gk=\gl$ and $\cof(\gk)=\gd>\gw$ from a model of $o(\gk)=\gl+\gd$,
so clause~(1) cannot be strengthened.  Another approach to the same
conclusion has been taken
by Miri Segal in \cite{segal.thesis}.
For $\gd=\gw$, Gitik and Magidor
\cite{gitik.negate-sch,gitik-magidor.ext-bas-force} 
show that the condition $o(\gk)\ge\gl$ cannot be improved in
clause~(3), and recent work of Gitik \cite{gitik.hidden} 
makes it unlikely that the second
condition in clause~(3) can be eliminated.  We will also show that
if there is an $n$ such that $\set{\ga:o(\ga)\ge\ga^{+n}}$ is bounded
in $\gk$ then the conclusion to clause~(2) can be strengthened to
match clause~(1), but it is not known whether this is true without the
added hypothesis.

If we assume that the GCH holds below $\gk$ then we can get slightly
more:
\begin{cor}[\ref{gch-main-cor}]
Suppose that $n>0$ and $\gk$ is a cardinal of cofinality $\gw$ such that
 $2^{\gk}\ge\gk^{+(n+2)}$ while $2^{\ga}=\ga^+$ for all $\ga<\gk$,  and
assume that there is an $m<\gw$ such that $\set{\ga:K\sat
o(\ga)\ge\ga^{+m}}$ is bounded in $\gk$.  Then
$o(\gk)\ge\gk^{{n+2}}+1$.
\end{cor}
Results in \cite{gitik-magidor.sing-cards-revisited-I} show that
$o(\gk)=\gk^{n+2}+1$ is sufficient.  
The restriction to $n>0$ is necessary here since  by results of
Woodin and Gitik \cite{gitik.negate-sch} $o(\gk)=\gk^{++}$
is enough to obtain a model of GCH with $\cof(\gk)=\gw$ and
$2^{\gk}=\gk^{++}$.

\paranumber{4}
For the case $\gk=\gw_\gw$ we have the following theorem:
\begin{theorem}
[\ref{aleph-gw1}]
If 
$2^{\gw}<\aleph_{\gw}$  and  $2^{\aleph_\gw}>\aleph_{\gw_1}$
then there is a sharp for a model with a strong cardinal.
\end{theorem}

\paranumber{5}
\medskip
The results concerning sequences of indiscernibles are much more
difficult to state.  The Dodd-Jensen covering lemma for $L[\mu]$
asserts that if $L[\mu]$ exists, but $0^{\dagger}$ does not exist,
then either every uncountable set $x$ of ordinals is contained in a
set in $L[\mu]$ of the same cardinality as $x$, or else there is a
sequence $C$ which is Prikry 
generic over $L[\mu]$ such that every uncountable set $x$ of ordinals
is contained in a set in $L[\mu,U]$ of the same cardinality as $x$.
Furthermore, the sequence $C$ is unique except for finite segments.
Uniqueness may fail if there are more measures in the core model:
starting from a model with inaccessibly many measurable cardinals it
is possible \cite{mitchell.skies}
 to construct a model in which each of the measurable
cardinals of $K$ has a Prikry sequence and hence is singular, but
there is no single system of indiscernibles for all of the cardinals.
A weaker uniqueness property is established in
\cite{mitchell.applications,mitchell.defsch}, however.
It is shown there that for each each uncountable
set $x$ of ordinals there is a function $h\in K$, 
a ``next indiscernible'' function $n$, and an ordinal $\rho$ of
cardinality at most $\card x^{\gw}$  such that $x$ is contained in
the smallest set $X_{\rho,h,n}$ containing 
$\rho$ and closed under the functions $h$ and $n$.
The function $n$ is somewhat complicated.  If $o(\gk)\le 1$ for all
$\gk$ then $n(\ga,\gg)$ is just the least indiscernible larger than
$\gg$ 
for the
measure $\ce_\ga$, but for larger cardinals it also must
generate certain limits of indiscernibles.  It is shown in
\cite{mitchell.applications,mitchell.defsch} that the
function $n$ is 
unique in the sense that for any other choice $\rho',h',n'$ there is
an ordinal $\eta<\sup x$ such that $n'(\ga,\gg)=n(\ga,\gg)$ whenever
$\ga,\gg\in X_{\rho,h,n}\cap X_{\rho',h',n'}$ and $\gg\ge\eta$.  In
this paper we extend these results up to a strong
cardinal, in the case
$\cof(\sup x)>\gw$, and to cardinals $\gk=\sup(x)$ such that
$\set{\ga<\gk:o(\ga)\ge\ga^{+n}}$ is bounded in $\gk$ for some $n<\gw$
in the case $\cof(\gk)=\gw$. 

\medskip
Section~1 is a brief introduction to the core model $K$ for sequences of
extenders and to its covering lemma.  It is intended to describe
notation used in the rest of the paper as well as to establish some
basic results concerning indiscernibles relative to a particular
covering set.  Most of the arguments which  require
a detailed reference to the proof of the covering lemma have been
gathered into this section,
so that with a few exceptions (mainly in subsection~3.2) 
the rest of the paper can be read in a
black box fashion, 
referring to results from section~1 rather than to basic core model
theory external to this paper.

\paranumber{6}
Section~2 covers results concerning sequences of indiscernibles.
The basic result is that such sequences are, except on a bounded
set, independent of the particular covering set used to obtain the
indiscernibles.  The applications to the singular cardinal
hypothesis are given in section~3, and some open problems are stated
in section~4.
\renewcommand\thetheorem{\thesection.\arabic{theorem}}

%


%


\errorcontextlines=100

\section{Introduction and Notation}
\paranumber{1}
We assume throughout this paper that there is no sharp for an inner
model with
a strong cardinal, so that a core model is guaranteed to exist. 
Expositions of these models include
 \cite{koepke.strong},
\cite{mitchell-steel.inner-model} and \cite{steel.core-model}.
The first of these uses somewhat different notation, and the latter two
are primarily concerned with 
larger cardinals and hence involve complications which are, from our
point of view, unnecessary.  Fortunately
our arguments will not make serious use of fine structure and hence are not
heavily dependent on the exact construction of the core model.

\paranumber{2}
The proof is heavily dependent on the covering lemma, and indeed on
the proof of the covering lemma.  We will begin this section with an
outline of this proof, partly to orient the reader and partly to
introduce the notation which will be used later in the paper.
Most of our 
references to the proof of covering lemma will be concentrated in this
section, so that a 
reader who is not fully comfortable with the details of the proof will
be able to get something out of the rest of the paper.

\paranumber{3}
\subsection{Extenders and the Core Model}
An $\gk,\gl$-extender~$E$ is a sequence of ultrafilters,
$E=\set{E_a:a\in[\gl]^{<\gw}}$, with
$E_{a}$ an ultrafilter on $^{a}\gk$.
An extender may be obtained from 
an embedding $\pi$ by
$$E_a=\set{x\subset{{}^a\gk}:\dot a\in \pi(x)},$$
where $\dot a=\pi^{-1}\restrict(\pi(a))$.
We will frequently identify a finite function $\gs\in{^a\gk}$ with
$\set{\gs(\xi):\xi\in a}\in[\gk]^{\card a}$, so that the equation
above could be written
$$E_a=\set{x\subset [\gk]^{\card a}:a\in \pi(x)},$$

\paranumber{4}
Going the other direction, an embedding $\pi$ can be generated from
the extender $E$
and a model $M$ which is to be the domain of $\pi$:
$$\pi\colon M\to\ult(M,E)=\set{[a,f]:a\in[\gl]^{<\gw}\text{ and }f\in
M\text{ and }f\colon {^a\gk}\to M}$$
where $[a,f]=[a',f']$ if and only if
$\set{\gs\in {^{a\cup a'}\gk}:f(\gs\restrict a)=f'(\gs\restrict
a')}\in E_{a\cup a'}$.
This will define an embedding on $M$ provided that $E_{a}$ is an
ultrafilter  on at least the subsets of $[\gk]^{\card a}$ which are in
$M$.  

\paranumber{5}
If $E$ is a $\gk,\gl$-extender then we call $\gk$ the {\it critical point} of
$E$, written $\crit(E)$.  If $\eta\le\gl$ then we write
$E\restrict\eta$ for the 
restriction of $E$ to the support $\eta$, that is,
$E\restrict\eta=\seq{E_a:a\in[\eta]^{<\gw}}$.
The {\it natural  length} of $E$, written
$\len(E)$, is defined to be the least ordinal $\eta\ge\gk^+$ such that 
$\ult(M,E)=\ult(M,E\restrict\eta)$ for any model $M$ such that $E$ is
an extender on $M$.

\medskip
\paranumber{6}
The core model, $\core$, is a model of the form $L[\ce]$, where $\ce$ is a
sequence of extenders and partial extenders on $L[\ce]$. 
Each member $\ce_\gg$ of the sequence $\ce$ is an extender on
$L(\ce\restrict\gg)$.  The set $\ce_\gg$ may or may not be a full
extender on all sets in  
$L[\ce]$: this depends on whether there are any 
subsets of $\crit(\ce_\gg)$ in $L[\ce]$ which are constructed after
$L_{\gg}[\ce]$ and hence are not measured by $\ce_\gg$. 
The ordinal $\gg$ is called the {\it index} of $E=\ce_\gg$,
written $\gg=\inx(E)$.
It is defined by $\inx(E)=\len(E)^{+}$ as evaluated in $L[\ce\restrict\gg]$.

\paranumber{6a}
The following theorem lists some of the properties of $K$ which we
shall need.  The proof can be found in the references.

\begin{theorem}
\label{K-prop}
The core model $K=L[\ce]$ is maximal among all iterable models $L[\cf]$ in
the following three senses:
\begin{enumerate}
\item
If $\mse m$ is a mouse which is coiterable with $K$ and agrees with
$K$ up  to the projectum of $\mse m$ then $\mse m\in K$.
\item
If $E$ is an extender such that $\ce\restrict\gg\cat E$ is good and
$\ult(K,E)$ is iterable then $E=\ce_\gg$.
\item
If $M=L_{\nu}[\cf]$ or $M=L[\cf]$ is iterable then there is an
iterated ultrapower of $K$ such that $M$ is a (possibly proper)
initial segment of the last model of the iteration.

\paranumber{6b}
Furthermore, if there is any elementary embedding $j\colon L[\ce]\to
M$ then 
this iterated ultrapower 
does not drop, and $j$ is the canonical embedding of this iterated
ultrapower. 
\end{enumerate}
\end{theorem}

For the  models in this paper the iterability properties referred to
above are all be guaranteed by 
countable completeness and hence are not a problem.  For core
models much larger than those considered here countable completeness is not
enough for iterability, so that iterability does become a serious problem.

Clause~3 is actually a combination of the global maximality property
that $K$ is not shorter than any model $L[\cf]$ which which it may be
compared, together with clauses~1 and~2.  The form of clause~3 which
we give here is not always true in core models larger than those
considered here.

\smallskip
\paranumber{6c}
The relation $\tri$ is defined for extenders in the same way  as for measures:
$E\tri E'$ in a model $M$ if and only if $E\in\ult(M,E')$.  
This
ordering is a well founded partial ordering.

We will write $\OO(\gk)$ to indicate
the set of $\gg$ such that $\ce_\gg$ is defined and is a full extender on
$\gk$ in $K$, and we will write 
$o(\gk)$ for the order type of $\OO(\gk)$.
We will will also write $\PO(\gk)$ for $\OO(\gk)\cup\{\gg\}$ where
$\gg$ is the strict sup of $\OO(\gk)$, that is,
$\gg=\sup\set{\nu+1:\nu\in\OO(\gk)}$. 

\paranumber{6d}
In some respects the models which we will use sit
uncomfortably between models with overlapping extenders and those with
no extenders other than measures. Our attention in later sections of
this paper will be largely devoted to the major difference, the
greater complexity of the indiscernibles, but there 
is one other difference which is more of an annoyance than a
problem and should be discussed here.
This problem is that if $\ce_\gg$ is an extender in $\ce$ with 
critical point $\gk$ then $i^{\ce_\gg}(\ce)$ may have extenders with
critical point $\gk$  which are not in $\ce$.
Suppose, for example, that $\ce_\gg$ is a measure on $\gk$ which
concentrates 
on cardinals $\ga<\gk$ such that $o(\ga)>\ga^{++}$.  Then
$\gg=\gk^{++}$ in 
$\ult(L[\ce],\ce_\gg)$ since $\ce_\gg$ is a measure, 
but $o(\gk)>\gk^{++}$ in $\ult(L[\ce],\ce_\gg)$.  
Thus, if we set $\ce'=i^{\ce_\gg}(\ce)$ then $\ce'_{\gg'}$ exists for
some ordinals $\gg'$ with
$\gg<\gg'\in\OO^{\ce'}(\gk)\setminus\OO^{\ce}(\gk)$.  If we had taken
the ultrapower by $\ce_\gg$ during the course of a comparison, because
$\ce_\gg$ was not a member of the other model in the comparison, then
it may be well be that some of then new extenders $\ce'_{\gg'}$ are
also not in the other model, requiring a second ultrapower by another
extender with the same critical point.

\paranumber{7}
The reader who is familiar with inner models for
overlapping extenders will recognize this situation
as a trivial example of an
iteration tree: one which is linear except that it
has side branches of length one in addition to
the main trunk.  For the less sophisticated reader we will sketch a
second solution.
This solution is simply to expand the sequence $\ce$, for the purpose
of the comparison lemma, so that if
$\ce_\gg$ has critical point $\gk$ then every extender on $\gk$ in the
sequence 
$\ce'=i^{\ce_\gg}(\ce)$ is also in the sequence $\ce$.  If we do this
then it is no longer true that $\ce_{\gg}\tri\ce_{\gg'}$ if and only
$\gg<\gg'$, but this failure is not such as to cause a serious problem.
  Notice that 
\[\sup(\OO^{\ce'}(\gk))<i^{\ce_{\gg}}(\gk)<(\gg^+)^{L[\ce\restrict\gg+1]},\]
which is smaller than the index of the next extender in $\ce$ on
$\gk$.  Thus the new extenders which appear in the ultrapower by
$\ce_\gg$ all lie between 
$\ce_\gg$ and the next extender on the original sequence.  The
expanded sequence will satisfy that $\ce_\gg\tri\ce_{\gg'}$ if and
only if $i^{\ce_\gg}(\gk)<i^{\ce_{\gg'}}(\gk)$.  We will write
$\gg\tri\gg'$ to mean that $\ce_\gg$ are extenders on the expanded
sequence with the same critical point, and $\ce_{\gg}\tri\ce_{\gg'}$.
In
addition we will write $\gg\tri\SPO(\ga)$ for all $\gg\in\OO(\ga)$.

\smallskip
\paranumber{8}
A cardinal $\gk$ is strong if  for all $\gl>\gk$ there is an
elementary embedding $i\colon V\to M$ such that $V_{\gl}\in M$.  Thus
$\gk$ is strong in a model $L[\cf]$ if and only if $O^{\cf}(\gk)$ is
unbounded in the ordinals.
The assumption that there is no sharp of a strong cardinal means
that there does not exist a pair $(\cf,I)$ such that $L(\cf)$
satisfies that there is a strong cardinal and that
$I$ is a proper class of
indiscernibles for $L[\cf]$.  The lack of such a sharp implies that
the extenders of $\ce$ never overlap, that is, there are no ordinals
$\gg$ and 
$\gg'$ in the domain of $\ce$ such that
$\crit(\ce_\gg)<\crit(\ce_{\gg'})<\gg$. 

\medskip
\notate{$\ult(M,\pi,\nu)$}
Although all of the extenders $E$ which we will be explicitly
considering are 
complete in the sense that each ultrafilter $E_a$ in $E$ is
$\gk$-complete, where $\gk=\crit(E)$,
we will use ultrapower constructions to define
extensions of elementary embeddings, and these constructions 
implicitly use extenders which are not complete.   If $\pi\colon N\to X$ 
then we will write $\ult(M,\pi,\nu)$ for the ultrapower of $M$ by the
extender of length $\nu$ generated by $\pi$, that is,
\[\ult(M,\pi,\nu)=\set{[a,f]:a\in[\nu]^{<\gw}\text{ and }f\in M}\]
where $[a,f]=[a',f']$ if and only if 
$$\dot a\cup\dot a'\in\pi\bigl(\set{\gs:f(\gs\restrict
a)=f'(\gs\restrict\ga')}\bigr).$$ 

\paranumber{9}
In order for $\ult(M,\pi,\nu)$ to exist, $N$ much contain all of the
sets which need to be measured in the ultrapower:
\begin{prop}
\label{ult-exists}
Let $\pi\colon N\to X$, with $\nu\in X$ and 
 let $\nu'$ be the least ordinal such that
$\pi(\nu')\ge\nu$.  Then $\ult(M,\pi,\nu)$ is defined whenever
\begin{align*}
&\ps(\nu')\cap M\subset N&\qquad&\text{if $\nu>\sup\pi\image\nu'$,}\\
&\forall\ga<\nu'\left(\ps(\ga)\cap M\subset N\right)&&\text{if
$\nu=\sup\pi\image\nu'$.}
\end{align*}
\qed
\end{prop}

\subsection{ The Covering Lemma}
It is assumed that the reader is familiar with the Dodd-Jensen
covering lemma \cite{dodd-jensen.cl-L[U],dodd-jensen.core} and with
the covering lemma for sequences 
of measures \cite{mitchell.coreI,mitchell.coreII}.
It will be recalled that Jensen's covering lemma for $L$ asserts that if
$0\sh$ does not exist and $x$ is any uncountable set then there is a set
$y\in L$ such that $x\subset y$ and $\card y=\card x$.  As 
larger cardinals become involved  the generalizations of the
covering lemma become more complex and less satisfactory, but
the proof remains essentially the same; indeed these generalizations
are still called ``covering lemmas'' not so much because their
statement looks like Jensen's covering lemma for $L$ as
because their proof looks like Jensen's original proof.  

\paranumber{10}
The first
step in the proof of the covering lemma is to replace the set $x$ with
a nicer set $X\supset x$ having the same cardinality as $x$:
\begin{defin}
A\notate{\pcs} {\it $\gd$-closed \pcs\ $X$\/} for $\gk$ is a set
$X\prec H_\tau$, for some $\tau\ge (2^{\gk})^+$,
such that $X^{\gd}\subset X$,  $\card X<\card\gk$
and $X$ is cofinal in $\gk$.
\end{defin}
Usually we will have $\gk=\sup(x)$ and  $\gd=\card x=\cof(\gk)$ and
$\tau=(2^{\gk})^{+}$, and in this case we will simply refer to $X$
as a \pcs.  On the few occasions when we use more or less closure, or
require 
$\tau$ to be larger than $(2^{\gk})^+$, we will so specify.

\begin{prop}
If $\gd<\gk$ are cardinals, $x\subset\gk$, and  $(\sup(\cof(\gk),\card
x)^{\gd}<\card\gk$ then there is a $\gd$-closed precovering set
$X\supset x$.\qed
\end{prop}
 
\medskip
\paranumber{10a}
In order to simplify notation, we will assume throughout the rest of
this section that we have a fixed 
\pcs\ $X$.  Later in the paper, when it may not be clear which \pcs\
is meant, we will modify the notation either by adding a subscript or
by specifying ``in $X$'' to indicate
which \pcs\ is intended.  Thus, in this section we will use
$\pi\colon N\cong X\prec H_\tau$ to denote the Mostowski
collapse of $X$, but if there were more that one \pcs\ involved we would
write
$\pi^X\colon N^{X}\cong X\prec H_{\tau}$.

\paranumber{11}
We will consistently use an over-bar to relate members of the
collapse $N$ of $X$ with the corresponding 
members of $X$.  If $x\in X$ then we write $\bar x$ for
$\pi^{-1}(x)$.  When $\bar x$ is used for some object which is not a
member of $N$ then the corresponding object $x$ will need to be
defined on a case by case basis, but it will always follow the rule
that $x$ is related to $\bar x$ via the embedding $\pi$.
\medskip

\paranumber{12}
By theorem~\ref{K-prop} there is an iterated ultrapower of $K$ with
final model $M_{\gth}$ having $\bark_{\bar\gk}$ as an initial
segment. 
Let $\seq{M_{\xi}:\xi\le\nu}$ be the sequence of models of the
iteration and $j_{\xi,\xi'}\colon M_{\xi}\to M_{\xi'}$ the
corresponding embeddings.

For most ordinals $\xi<\gth$ we will have
$M_{\xi+1}=\ult(M_{\xi},E)$ where $E$ is the least extender which is
in $M_\xi$ but not in $\bark$, but
for finitely many ordinals $\xi<\gth$
the iteration may {\it drop to a mouse.}
This means that
$M_{\xi+1}=\ult(M_{\xi}^{*},E)$ where $M_{\xi}^*$ is a mouse such that
$M_\xi^*\in M_{\xi}$.
This happens whenever there is a subset $x\subset\rho$, with 
$x\in M_{\xi}\setminus \bark$, for some ordinal $\rho$ which is
less than or equal to the
critical point of the first extender on which $M_\xi$ and $\bark$
disagree. 

The next step depends on whether the iteration ever does drop to a
mouse before reaching a model $M_{\gth}$ which agrees with $\bark$
up to $\bar\gk$.  Jensen, in his proof of the covering lemma for $L$,
was able to prove outright that this must happen
 by observing that otherwise the embedding $\pi$ could be
extended to a nontrivial embedding from $L$ into $L$, which implies
that $0{\sh}$ exists.  The argument works for sequences of measures,
but can fail for extenders.  We will sketch a proof that shows that if
there are no overlapping extenders and the iteration does not drop
then full covering holds over $K$ for cofinal subsets of $K$, so that
$2^{\gk}=\gk^+$ by the same proof as for $L$.  This proof is given in
detail (for overlapping extenders) in \cite{mit-sch-ste.CovLemWoodin}.

Suppose that the iteration does not drop.  Then $M_{\gth}$ is a proper
class and $j_{0,\gth}\colon K\to M_{\gth}$ exists.  Let $k\colon
M_{\gth}\to \barm=\ult(M_{\gth},\pi,\gk)$ be the canonical
embedding. By theorem~\ref{K-prop} there is an iterated ultrapower
$i\colon K\to \barm$ such that $i=k\comp j_{0,\gth}$.  Now
$\crit(\tilde\pi\comp j_{0,\gth})\le\crit(\pi)=\eta$, so
$\crit(i)\le\eta$.  The first ultrapower in $i$ uses an extender $E$
in $K$ which is not in $\barm$, and since $\barm$ agrees with $K$ at
least up to $\gk$ it follows that $\len(E)\ge\gk$.  Since there is no
model with overlapping extenders it follows that there are no
measurable cardinals $\mu$ in the interval $\eta<\mu\le\gk$.  Then the
iteration $i_{0,\gth}$ involves only finitely many ultrapowers before
reaching $\bar\gk$ (\textit{cf} the proof of lemma~\ref{cl}), so there
is an ordinal $\bar\rho<\bar\gk$ such 
that any member of $\bar\gk$ can be expressed in the form
$i_{0,\gth}(f)(\gg)$ for some $f\in K$ and $\gg<\bar\rho$.  It follows
that $X\subset y=\set{i^{E}(f)(\gg):f\in{^{\eta}\eta}\cap K\land
\gg<\rho}$.  Now $y\in K$ since $E\in K$, and $\card y\le\gg
\,2^{\card{\crit(E)}}\le \gg\, 2^{\eta}<\gk$.  Thus full covering holds for
subsets of $\gk$, whenever the iteration does not drop, which is what
we were trying to show.  For the rest of this paper we will assume
that the iteration does drop.

\paranumber{14}
In order to simplify notation it is convenient to use the critical
points of the extenders to index the 
models in the iterated
ultrapower, rather than indexing them 
sequentially as in the last paragraph.  Let
$j_{\xi,\xi'}$ be the canonical embedding from $M_{\xi}$ to
$M_{\xi'}$, which is defined provided that the iteration does not drop
to a mouse in 
the half-open interval $[\xi,\xi')$.  

\paranumber{15}
\begin{defin}
\notate{$\barmse
m_{\bar\nu}$}\notate{$\bari_{\bar\nu,\bar\nu'}$}\notate{$\bare_{\bar\nu}$}  
If $\nu$ is an ordinal in $N$ then $\barmse m_{\nu}\eqdef M_{\xi_{\nu}}$,
where $\xi_\nu$ is the least ordinal $\xi$ such that
$\crit(j_{\xi,\xi+1})\ge\nu$.  If there is no such ordinal  $\xi$ then
 $\barmse m_{\nu}=M_{\gth}$.  
We write ${\bari_{\nu,\nu'}}$ for the 
embedding $j_{\xi_\nu,\xi_{\nu'}}\colon \barmse m_{\nu}\to\barmse
m_{\nu'}$,  and $\bare_\nu$ for the extender used at stage $\xi_\nu$
of the iteration,
so  $\barmse m_{\nu+1}=\ult(\barmse m_{\nu},\bare_{\nu})$.  We write
$\barh_{\nu}$ for the canonical 
Skolem function of the premouse $\barmse m_{\nu}=M_{\xi_\nu}$.
\end{defin}

Note that if $\nu=\crit(j_{\xi,\xi+1})$,  then $\barmse
m_\nu=M_\xi$ and $\barmse m_{\nu+1}=M_{\xi+1}$.

\paranumber{16}
The embedding $\bar i_{\nu,\nu'}$ does not exist if the iteration
drops to a mouse somewhere in the half open interval
$\xi_\nu\le\xi<\xi_\gk$.   
Such drops only occur finitely often.  Those familiar with fine
structure in these models will recall that the iteration may also drop
in degree, but this also occurs only finitely often.
Thus there is an ordinal $\bar\nu_0<\bar\gk$ such
that the iteration never drops in the interval
$\xi_{\bar\nu_0}\le\xi<\xi_{\gk}$, so that $\bar i_{\nu,\nu'}$ is always
defined when $\bar\nu_0\le\nu<\nu'\le\bar\gk$.  Since we are only
interested in subsets of $\gk$, and are not concerned with what
happens on bounded subsets of $\gk$, it will be sufficient to restrict
ourselves to $\nu$ in this interval.

If we are doing fine structure in terms of $\gS_n$-codes, then we can
think of the models $\cm_{\nu}$ for $\bar\nu_0\le\nu<\gk$ as
$\gS_{n-1}$-codes, for some fixed $n<\gw$, for premice
$J_{\ga_\nu}[\cf_{\nu}]$.  All of the models have the same $\gS_1$-projectum
$\bar\rho<\bar\nu_0$, so that $j_{\nu,\nu'}(\bar\rho)=\bar\rho$ for
$\nu_0\le\nu<\nu'\le\bar\gk$.   The Skolem function $\bar f_{\nu}$ of
$\barmse m_{\nu}$ is just the canonical $\gS_1$-Skolem function and is
also preserved by the maps $j_{\nu,\nu'}$.  The reader who does not
full understand  this construction will not be mislead if he
thinks only of the case $n=1$, so that $\barmse
m_{\nu}=J_{\ga_\nu}(\cf_{\nu})$ and the embeddings $j_{\nu,\nu'}$ are
ordinary ultrapowers by functions in $\barmse m_{\nu}$.

\newcommand\bmsegk{\barmse m_{\bar\gk}}

\paranumber{18}
So far we have concentrated on the collapsed model $N$, but we
are really interested  in the uncollapsed model $X$.  The connection
between the two is made by using the collapse map $\pi$ as an
extender.  In particular we can define $\tilde\pi$ to be the canonical
embedding from  $\bmsegk$ into $\mse m=\ult(\bmsegk,\pi,\gk)$, which
exists by proposition~\ref{ult-exists} since $\gk=\sup\pi\image\bar\gk$.
This
embedding preserves the fine structure of $\bmsegk$, so  
$\tilde\pi\comp \bar h=h\comp \pi$ where $\bar h$ and $h$ are the
Skolem functions of $\bmsegk$ and $\mse m$ respectively.  

The usual proof of the next lemma consists mainly of the proof that
$\mse m$ is iterable.  The proof with extenders involves one
additional difficulty, and we include just enough of the proof to
indicate a solution to this problem.
It should be noted that $\inx(\bare_{\bar\nu})<\bar\gk$ for all
$\bar\nu<\bar\gk$.  If, to the contrary,
$\inx(\bare_{\bar\nu})\ge\bar\gk$ then $\barmse m_{\bar\nu}$ agrees with 
$\bark$ up to $\bar\gk$, so that the iteration was already complete at
$\barmse m_{\bar\nu}$ before $\bare_{\bar\nu}$ was chosen.

\paranumber{19}
\begin{lemma}
The structure $\mse m=\ult(\bmsegk,\pi,\gk)$ is a member of $K$. 
\end{lemma}

\begin{proof}[Sketch of proof]
The new difficulty 
is that there may be an extender $\bare$ on the extender sequence of
$\bmsegk$ such that $\crit(\bare)<\bar\gk\le\len(\bare)$.
This is no problem if $\bare$ is an actual member of $\bmsegk$, since
in this case $\tilde\pi(\bare)$ is defined and is in $K$.
Thus we need only worry about the
case when $\bare$ is the last extender in the sequence of $\bmsegk$.
In this case standard arguments show that the structure $\mse m'$
obtained by omitting the final
extender of $\mse m$ is a member of $K$. 

\paranumber{19a}
Set $\bar\mu=\crit(\bare)$ and $\mu=\pi(\bar\mu)$.
If $z\subset\ps(\bar\mu)$ is an arbitrary member of $\bark$ which
has 
cardinality $\bar\mu$ in $\bark$, then by amenability the set $E\cap
z=\set{\bare_a\cap z:a\in[\len(\bare)]^{<\gw}}$ is a member of
$\bmsegk$.  Thus we can define $E=\union_{z}\tilde\pi(\bare\cap z)$.
If $X$ is cofinal in $\mu^+{^{(K)}}$ then
$E$ is a full extender on $K$.  In that case standard arguments show
that $E$ is in $K$, so that $\mse m$ is
a member of $K$.

\paranumber{20}
The referee has pointed out that we can ensure that $X$ is cofinal in
$\mu^{+}$ of $K$ by choosing the
precovering set $X$ so that  $\mu^+\cap X$ is cofinal in $\mu^+$ whenever
$\mu$ is $<\gk$-strong in $K$.  This is possible since our assumption
that there are no overlapping extenders implies that there can be at
most one such cardinal $\mu$, and our assumption that 
$\gk$ is a strong limit cardinal implies that any subset of $\gk$ of
cardinality less than $\gk$ is contained in a precovering set.  

For the sake of the interested reader we will sketch a proof, without
this extra assumption on $X$, that $E\in K$ even when $X$ is not
cofinal in $\mu^+{^{(K)}}$.  In this case let $\mse n$ be the least
mouse which has projectum less then or equal to $\mu$ and such that
$\mse n$ is larger than every mouse in $X$ with projectum $\mu$.  Then
$\mse n$ is the least mouse in $K$ such that there is a subset of
$\ga=\pi(\bar\ga)$ definable in $\mse n$ which is is not decided by
$E$, so that $E$ is an extender on $\mse n$ and we can let $i^E\colon
\mse n\to\mse n'=\ult(\mse n,E)$ 
be the canonical embedding.  Then
$\mse n'$ is an iterable premouse which agrees with $K$ up to the
length of $E$, so $\mse n'$ is a member of $K$.  Now
$\range(i^E)=h^{\mse n'}\image \ga$, and hence the range of $i^E$ is
definable in $\mse n'$.  Then for any $a\in[\len(E)]^{n}$ the ultrafilter
$E_a$ is equal to the set of $x\subset[\mu]^n$ such that there is a set $y\in\range(i^E)$
such that $y\cap[\mu]^n=x$ and $a\in y$, and it follows that $E\in K$.
\end{proof}

Notice that $\mse m$ as defined in the last section is not a mouse,
since $E$ is not a complete extender on $\mse m$.  It is close enough
for our purposes, however, since its Skolem function $h$ still
satisfies the crucial identity $\tilde\pi\comp \bar
h=h\comp \pi$.

\paranumber{21a}
Now let $\bar h$ and $h$ be the Skolem functions of $\bmsegk$ and $\mse m$
respectively.  Then $h\in K$, and since $\tilde\pi\comp \bar
h=h\comp \pi$  we can use $h$ to cover the set $X$ as follows:
Set
$\rho=\rho^X=\pi(\bar\rho)$ where $\bar\rho$
is the projectum of $\bmsegk$, and let $I$ be the set of
ordinals $\pi(\xi)$ such that
$\crit(\bare_{\nu})\le\xi<\len(\bare_{\nu})$ for some extender
$\bare_{\nu}$ used in the iteration which gave $\bmsegk$.
Then $X\cap\gk=\pi\image\bar\gk\subset h\image(\rho\cup I)$.  

We will call the members of $I$ {\it indiscernibles} by analogy with
the simpler case of sequences of measures.
We will eventually need to make a detailed analysis of these
indiscernibles, but 
first we look at the  covering lemma to see what can be
obtained by looking at intervals in which there are no
measures and hence no indiscernibles:

\paranumber{21b}
\begin{lemma} {\rm (Covering Lemma without indiscernibles)}
\label{cl}
Assume that there does not exist a sharp for a model with a strong
cardinal,  that  $\gl^{\gw}<\gk$, and that there are no
measurable cardinals $\nu$ of $K$ in the half-open interval
$\gl<\nu\le\gk$.  Then for every subset  $y$ of $\gk$ such that
$\card y\le\gl$  there is a $z\in K$ such that
$y\subset z$ and $\card z\le\gl$.

\paranumber{21c}
In particular if $\gk>\gw_2$ is a regular cardinal in $K$ then
$(\cof(\gk))^{\gw}\ge\card\gk$, and if $\gl$ is a singular cardinal in
$V$ then either $\nu^{\gw}\ge\gl$ for some $\nu<\gl$ or else
$\gl^+=\gl^+{^{(K)}}$.  
\end{lemma}

\proof
The proof is by induction on $\gk$.  From the discussion above we know
that any subset $x$ of $\gk$ is contained in a set of the form
$h\image(\rho\cup I)$ where $h\in K$, $\rho<\gk$, and $I$ is the set of
indiscernibles.
We will show that there is an ordinal $\eta$, with $\rho\le\eta<\gk$, such
that $I\setminus\eta$ is finite.  It follows that $x$ is contained in
a set $y\in K$ such that $\card{y}^K\le\eta$.  By the  induction
hypothesis it follows that $x$  is contained in a set $y'\in K$ such
that $\card{y'}^K\le\gl$.

If there is some $\nu$ such that
$\crit(\bare_\nu)\le\rho<\pi(\len(\bare_\nu))$ then set
$\eta=\pi(\len(\bare_\nu))$,  and otherwise set $\eta=\rho$.  Thus
every member of $I\setminus\eta$ comes from an extender
$E=\bare_{\nu}$ such that $\eta\le\pi(\crit(E))<\gk$.
Notice that $E$ must be a measure, since
otherwise $\xi=\crit(E)$ is measurable in
$\barmse m_{\nu+1}=\ult(\barmse m_{\nu},\bare_{\nu})$, which implies
that $\xi$ is measurable in 
$\bark$ and hence $\pi(\xi)$ is measurable in $K$, contrary to assumption.
Furthermore $\bari_{\nu,\bar\gk}(\bar\xi)>\bar\gk$, since otherwise 
$\pi(\bari_{\nu,\bar\gk}(\bar\xi))$ is measurable in $K$.  If there
are infinitely many measures $\bare_\nu$ with
$\eta\le\xi_\nu=\crit(\bare_\nu)<\bar\gk$ then there must be an
infinite set $D$  of  $\nu$ so that
$\bari_{\nu,\nu'}(\xi_\nu)=\xi_{\nu'}<\bar\gk$ for $\nu<\nu'$ in $D$, but then
any limit point  $\gb$ of $D$ of cofinality $\gw$ is
measurable in $\bark$.
To see this, let $\vec d\subset D$ be a countable set with
$\gb=\sup\vec d$, and let $U$ be
the set of $x\in\ps(\gb)\cap \bark$ such that for all sufficiently
large $d\in\vec d$ 
\begin{align*}
d&\in x&&\text{if $o^{\bark}(d)=0$}\\
x\cap d&\in U_{d}&&\text{if $o^\bark(d)>0$}
\end{align*}
where $U_d$ is the order~$0$ measure on $d$ in
$\bark$.  Then $U\in N$ since $\vec d\in{^{\gw}}N\subset N$, and $U$
is a measure on $\bark$.
 Thus
$\eta\le\pi(\gb)\le\gk$ and $\gb$ is measurable in $K$, contrary to
assumption, and so $I\setminus \eta$ must be finite as required.
\endproof

\paranumber{21e}
Now we must prepare for the hard work of analyzing the indiscernibles.
The preparation will take up the rest of the section, and the actual
analysis will be carried out in the next section. 
So far we have concentrated on the collapsed model $N$, and on the
collapse $\bark$ of $K$, but  we are really interested in the 
uncollapsed models  $V$ and $K$.
In the rest of this section we will describe the relationship, induced
by the map $\pi$, between
objects in $X$ and objects of~$N$. 

\paranumber{21f}
We write $\rho^X$ for $\pi(\bar\rho)$, where $\bar\rho$ is the
$\gS_1$-projectum of $\barmse m^*_{\bar\nu_0}$, that is, $\bar\rho$ is
the least 
ordinal such that there  is a subset $x$ of $\bar\rho$ which is $\gS_1$
definable in the $\gS_{n-1}$-code $\barmse 
m_{\bar\nu_0}$ such that $x\notin \barmse m_{\bar\nu_0}$. This same
ordinal $\bar\rho$ will 
be the projectum of all of the models $\barmse m_{\nu}$ for
$\bar\nu_0<\nu\le\bar\gk$.  

Let $C$ be the set of ordinals $\nu\in X$ such that
$\bar\nu=\pi^{-1}(\nu)=\crit(\bari_{\bar\nu,\bar\nu+1})$.  If we were
dealing with 
measures then $C$ would be the set of indiscernibles, but we will 
call any member of
$\union\set{\pi\image(\len(\bare_{\bar\nu}\setminus\bar\nu)):\nu\in C}$
an {\it indiscernible.}  The members of $C$ are called {\it principal
indiscernibles}.

\paranumber{22}
\begin{defin}\label{ih-def}
Suppose that $\nu_0<\nu<\nu'\le\gk$ and $\nu,\nu'\in C$.
\begin{enumerate}
\item
$i_{\nu,\nu'}\eqdef\pi\comp \bari_{\bar\nu,\bar\nu'}\comp \pi^{-1}\restrict(X\cap\eta)$, 
where 
$\eta$ is the least inaccessible cardinal of $K$ above $\SO(\nu)$.
\item
$h_{\nu}\eqdef\pi\comp \barh_{\bar\nu}\comp  \pi^{-1}\restrict
\set{\xi\in X:\pi\comp  \barh_{\bar\nu}\comp \pi^{-1}(\xi)<\eta}$.
\end{enumerate}
\end{defin}

\paranumber{23}
Notice  that $\eta\in X$. The following proposition implies that
$\eta>\pi(\inx(\bare_{\bar\nu}))$.  It   works for $\nu=\gk$ as well
if we take $\bare_{\bar\gk}$ to be the $\tri$-least extender which is in
$\barmse m_{\bar\gk}$ but is not in $\bark$.
\begin{prop}\label{inxE-bnd}
Suppose that $\cf$ is an extender sequence and
$\tau=\crit(\cf_{\gg})$. Then
$\ult(L[\cf],\cf_{\gg})\sat\SO(\tau)^{+}\ge\gg$. 
\end{prop}
\begin{proof}
Recall that $\gg=\len(\cf_\gg)^+$ as evaluated in $L[\cf\restrict\gg]$
or, 
equivalently, as evaluated in the ultrapower $\ult(L[\cf],\cf_{\gg})$.  Thus it is enough to show
that $\SO(\tau)\ge\card{\len(\cf)}$ in $\ult(L[\cf],\cf_{\gg})$.
Consider the extenders $\cf\restrict\eta$ for $\eta<\len(\cf_{\gg})$.
All of these extenders are in $\ult(L[\cf],\cf_{\gg})$, so if
$\SO(\tau)<\card{\len(\cf)}$ then there is $\eta_0$ such that
$\ult(L[\cf],\cf_{\gg}\restrict\eta)=\ult(L[\cf],\cf_{\gg}\restrict\eta_0)$
for all $\eta$ in the interval $\eta_0<\eta<\len(\cf_{\gg})$.   It
follows that every ordinal in that interval can be written in
$L[\cf\restrict\gg]$ in the form $i^{\cf_\gg\restrict\eta_0}(f)(a)$ for
some $f\colon\gk^{n}\to\gk$ and some $a\in[\eta_0]^{\lt\gw}$.  Thus
$\card{\len(\cf_{\gg})}\le\eta_0$ in $L[\cf\restrict\gg]$.  
\end{proof}

\paranumber{23a}
Next we need to consider the image under $\pi$ of the extenders
used in the iterated ultrapower.
\begin{defin}
  Suppose that $\nu'\in C$  and $\nu'$ is a \pri\ for $\nu$ in~$X$.
\begin{enumerate}
\item\usage{$\barf_{\bar\nu,\bar\eta}$}
$\barf_{\bar\nu',\bar\nu}\eqdef\bari_{\bar\nu',\bar\nu}(\bar E_{\bar\nu'})$.
\item
\usage{$F_{\nu,eta}$}
$F_{\nu,\eta}\eqdef \pi(\barf_{\bar\nu,\bar\eta})$ if
$\barf_{\bar\nu, \bar\eta}\in\bark$, and it is undefined otherwise.
\end{enumerate}
\end{defin}

Eventually we will show, using proposition~\ref{inxE-bnd}, that for
many $\bar\nu$, $\barf_{\bar\nu,\bar\gk}$ is in $N$ and hence
$F_{\nu,\gk}$ exists and is in $K$.

\paranumber{23b}
Notice that $\barf_{\bar\nu',\bar\nu}$ is a member (or the last
extender) of $\barmse
m_{\bar\nu}$. 
If $\bar\nu<\bar\gk$ then $\barf_{\bar\nu',\bar\nu}$ is a member of
$\bark$ if and only if either 
$\bar\nu=
\crit(\bare_{\bar\nu})>\bari_{\bar\nu',\bar\nu}(\crit(\bare_{\bar\nu'}))$
(in which case $\bar\nu>\crit(\barf_{\bar\nu',\bar\nu})$ and  every
extender in $\barmse m_{\bar\nu}$ with critical 
point less than $\bar\nu$ is in $\bark$) or 
$\bar\nu=\bari_{\bar\nu',\bar\nu}(\crit(\bare_{\bar\nu'}))$ and
$\barf_{\bar\nu',\bar\nu}\tri\bare_{\bar\nu}$ (in which case
$\bare_{\bar\nu}$ is the $\tri$-least extender in $\barmse
m_{\bar\nu}$ which is not in $\bark$).
Essentially the same analysis works at $\bar\gk$:
$\barf_{\bar\nu,\bar\gk}$ is in $\bark$ if either every extender on
$\bar\gk$ in $\barmse m_{\bar\gk}$ is in $\bark$, or else
$\barf_{\bar\nu,\bar\gk}\tri\bare_{\bar\gk}$ where $\bare_{\bar\gk}$
is the 
$\tri$-least extender  which is in $\barmse
m_{\bar\gk}$ but not in $\bark$.

\paranumber{24}
\begin{lemma} \label{ftoi-lem}
Suppose that $\nu_0<\nu<\nu'\le\gk$ and that $\nu$ is a \pri\ for
$\nu'$ in $X$, that is, that $\nu\in C$ and
$\bar\nu'=\bari_{\bar\nu,\bar\nu'}(\bar\nu)$. Then
\begin{enumerate}
\item
$h_{\nu'}\restrict \nu=i_{\nu,\nu'}\comp  h_{\nu}$. 
\item
If $z$ is in $X\cap K_{\eta}$,  where
$\eta$ is the least inaccessible cardinal of $K$ above $o(\nu)$,
 then $z$ is in
$h_\nu\image\nu$.  Indeed $z=h_\nu(\vec d)$ where $d$ is a finite
sequence of ordinals, each of which is either in  $\rho$ or an
indiscernible smaller than~$\nu$.
\item
If $z\in h_{\nu'}\image\nu$, $\vec b\in\left[\pi\image\len(\bar
E_{\bar\nu})\right]^{\lt\gw}$, and
$\vec b'=i_{\nu,\nu'}(\vec b)$ then 
$z\in 
(\pi\image\barf_{\bar\nu,\bar\nu'})_{\vec b'}$ if and only if $\vec
b\in z$. 
\item
If $f\in X\cap K$, the ordinal $\nu$ is a limit of $C\cap\nu$, and
$\gg$ is the least member of $C$ above $\nu$ 
then $\gg\cap f\image \xi =h_\nu\image(\xi\cap f\image \nu)$
for every sufficiently large ordinal $\xi$ in $C\cap\nu$.
\item
If $y\subset\nu$ and $\card y\le\gd$ then there are functions
$i'$, $h'$ and $h''$ in $X\cap K$ such that $i'\restrict
y=i_{\nu,\nu'}\restrict y$, $h'\restrict y=h_{\nu}\restrict y$, and
$h''\restrict y = h_{\nu'}\restrict y$.
\end{enumerate}
\end{lemma}

\paranumber{25}
\begin{proof}
Clause~(1) follows  from the definition of $i_{\nu,\nu'}$ and
$h_{\nu}$, and clauses~(2) and~(3) 
follow from the corresponding facts about the iterated ultrapowers
$\barmse m_{\nu}$ and $\barmse m_{\nu'}$ and the fact that
$\len(\bare_{\bar\nu})$ is smaller than $\eta$.
Clause~(4) follows from the fact that $f$ is in the range of
$i_{\nu',\nu}$ for  some $\nu'\in C\cap\nu$.

\paranumber{26}

This leaves only clause~(5) to be proved.  By clause~(1) we can set
$i'=h''{\circ} (h')^{-1}$, so it will be enough to show that the
functions $h'$ and $h''$ exist.   The proof is
identical for
$h'$ and $h''$, so we will only give the proof for  $h'$. 
Now $h_{\nu}\image y\subset
X$ and $X$ is $\gd$-closed so 
$h_\nu\image y\in X$.  Thus we can apply lemma~\ref{cl},
the  covering lemma without indiscernibles,
inside $X$.  Since there are no measurable cardinals in $K$ between $\nu$ and
$\sup(\range(h_\nu))$, lemma~\ref{cl} asserts that there is a function
$f\in X\cap K$  such 
that $h_{\nu}\image y\subset f\image\nu$.  Then $\bar f=\pi^{-1}(f)\in
N$, and since the next member of $C$ above $\nu$ is larger than
$\sup(\range(h_\nu))$ it follows that $\bar f\in \barmse m_{\bar\nu}$.  

\paranumber{27}
The model $\barmse m_{\bar\nu}$ must have cofinality greater than
$\gd$.  To see why this is true, recall that $\barmse m_{\bar\nu}$ is
the $\gS_{n-1}$ code $(J_{\ga}[\ce],A)$ of some premouse, and has
$\gS_1$-projectum $\rho\le\bar\nu$.  Let $x\subset \rho$ be $\gS_1$
definable in $\barmse m_{\bar\nu}$, but not a member of $N$. Then
$x=\union\set{x_{\ga'}:\ga'<\ga}$, where $x_{\ga'}$ is the set of
$\xi<\rho$ such that there is a witness $z\in J_{\ga'}[\ce]$ of the
$\gS_1$ fact ``$\xi\in 
x$''.  Then each set $x_{\ga'}$ is in $N$, and
if $\cof(\ga)\le\gd$ then it would follow that $x\in N$. 

Since $\cof(\ga)>\gd$  and $\bar h_{\bar\nu}$ is $\gS_1$ definable in
$\barmse m_{\bar\nu}$, it
follows that there is a set $\bar y_0\in \barmse m_{\bar\nu}$ such
that $\bar y=\pi^{-1}(y)\subset\bar y_0$ and
$\bar h_{\bar\nu}\restrict\bar y_0\in\barmse m_{\bar\nu}$. 
Define a partial function $\bar g\colon \bar y_0\to \bar\nu$
by  letting $\bar g(\xi)$ be the least ordinal $\eta$ such that
$\bar h_{\bar\nu}(\xi)=\bar f(\eta)$. Then $\bar f\comp \bar g=\bar
h_{\bar\nu}\restrict\bar y_0$. 

\paranumber{28}
  Now we have to consider two
cases.  If $\nu_0<\nu<\gk$ then by the choice of $\nu_0$ the iteration
did not drop to a mouse at $\barmse m_{\bar\nu}$, that is,
$\ps(\bar\nu)\cap\barmse m_{\bar\nu}\subset N$.  In particular $\bar
g\in N$ and we can set $h'_{\nu}=f\comp\pi(\bar g)\in K\cap X$, so
that 
$h'(\xi)=h_{\nu}(\xi)$ for all $\xi\in X\cap\domain(h')\supset y$. 

The only other possibility is $\nu=\gk$, in which case $X$ is cofinal
in $\gk$ so that we can define $\tilde\pi\colon \barmse
m_{\bar\gk}\to\mse m=\ult(\barmse 
m_{\bar\gk},\pi,\gk)$.  Then $\mse m\in K$ as in the proof of
lemma~\ref{cl}, so $h^*=\tilde\pi(\bar g)\comp f\in K$.
Now $h^*(\xi)=h_{\gk}(\xi)$ for all $\xi\in X\cap\pi(y_0)$, so
$h^*\restrict y=h_\gk\restrict y$.  But
$h_{\gk}\restrict y\in X$, so by the elementarity of $X$
there is a function $h'\in X\cap K$ such that $h'\restrict
y=h_{\gk}\restrict y$. 
\end{proof}

\subsection{Indiscernibles}
\paranumber{29}
It only remains to  briefly discuss our notation for indiscernibles
before we can begin the analysis of sequences of indiscernibles.  As
stated before, we call $\nu$ a {\it principal indiscernible} if
$\nu\in C$, that is, if
$\nu=\pi(\bar\nu)$ where $\bar\nu$ is the critical point of
$\bari_{\bar\nu,\bar\gk}$ or equivalently if $\nu$ is the critical
point of $i_{\nu,\gk}$.

\paranumber{30}
We will say that $a$ is a {\it \pri\ for $\ga$} if $a\in C$ and
$i_{a,\ga}(a)=\ga$.  We say that $a$ is an \pri\ for the extender $E$
on $\ga$ if $E=F_{a,\ga}$.  Notice that if $a$ is a \pri\ for $\ga$
and $F_{a,\ga}\notin K$ then $a$ is not an indiscernible for any extender on
$\ga$. This differs from the way the term is usually used, but it is
useful here because we will spend a large part of the next section
showing that the relation ``$\nu$ is an indiscernible for $\ga$'' is
definable before we begin to look at  the definability of 
the relation ``$\nu$ is an indiscernible for the extender $F$ on
$\ga$.''

\paranumber{31}
As stated earlier, we will say that an ordinal $b$ is an indiscernible 
whenever there is a principal indiscernible $a$ 
such that $b\in \pi\image(\len(\bare_{\bar a}\setminus\bar a))$, where
$a=\pi(\bar a)$.  Since these indiscernibles will be used to
reconstruct the image of the extender used at stage $a$ 
it will be convenient to generalize this notation: 
\begin{defin}\label{indforb-def}
An ordinal $b$ is  an {\it indiscernible for $\gb$ belonging to
$(a,\ga)$} if (i)~$a$~is a \pri\ for $\ga$, (ii)~$\gb=i_{a,\ga}(b)$, 
(iii)~$\gb<\inf(C\setminus \ga+1)$,   and (iv)~$b$~is smaller than 
$\SO(a)^{++}$.
\end{defin}

Notice that if $\bare_{\bar a}$ is the extender used at the $\bar a$th
stage of the iteration then proposition~\ref{inxE-bnd} implies that
$\pi(\inx(\bare_{\bar a}))$ cannot be 
larger than the
upper bound on $b$ given in clause~(iv), and that this upper
bound is smaller than the next inaccessible cardinal
above $\SO(a)$  in $K$.

\paranumber{32}
We will consistently use Roman letters for indiscernibles and the
corresponding Greek letters for the ordinals for which they are
indiscernibles.  Thus $a$ will denote a \pri\ for $\ga$, and $b$ and $c$
will denote indiscernibles for $\gb$ and $\gg$ respectively.

\medskip
It was pointed out earlier that all of the definitions in this section
are relative to a fixed \pcs\ $X$.  Whenever it is not clear which
\pcs\ is being used we will specify the relevant \pcs, either by
adding a superscript or by using the words ``in $X$''.

\paranumber{33}
Unless otherwise specified, successors are always calculated in the
core model $K$.  Thus $\gk^{+n}$ means the $n$th successor as
calculated in $K$.  Other functions will still be calculated in $V$
unless otherwise stated, so that $\card x$ and $\cof(\gk)$ are the
cardinality of $x$ and cofinality of $\gk$ in the real world.

The letter $h$ will always be used to denote a Skolem function, and if
$x$ is a set then we
will write $h\image x$ to mean $\set{h(\vec\nu):\vec\nu\in
[x\cup\gw]^{<\gw}}$.


%

%

\newcommand\param{\vec e}
\newcommand\uparam{\Gamma}

\section{Definability and Uniqueness\\of Indiscernible Sequences}

\paranumber{1}
The covering lemma for one measure
\cite{dodd-jensen.core,dodd-jensen.cl-L[U]} asserts that if 
$0^{\dagger}$ does not exist then any uncountable set $x$ of ordinals 
is contained
in a set $y$ such that $\card y=\card x$ and either $y\in\core$
(where $K=L[\mu]$ if it exists and $K$ is the Dodd-Jensen core
model otherwise) or else
$y\in L[\mu,C]$ where 
$C$ is a Prikry sequence for the measure $\mu$.  Furthermore 
the Prikry sequence $C$ is unique up to initial segments: any other
Prikry sequence over $L[\mu]$ is contained in $C$ except for a finite
set.
If there are sequences of measures in $\core$ then it is still true
that each individual measure has a unique maximal Prikry sequence, 
but there need not be a uniform system of 
indiscernibles for the whole sequence of measures
\cite{mitchell.skies}.
It is true that any small set of measures has  a system of
indiscernibles, but 
the particular system
of indiscernibles to be used to cover a given set $x$ depends on 
the set $x$.
A modified version of the uniqueness of the sequence $C$ does extend
to sequences of measures, however. It is shown in 
\cite{mitchell.sing-card-hyp,mitchell.applications,mitchell.defsch},
that,
roughly speaking, if we
specify a small set of measures for which we want indiscernibles, then
the system of indiscernibles for that set of measures is
unique up to an initial 
segment.  In this section we will generalize these results to models
containing nonoverlapping extenders.

\paranumber{2}
\medskip
We have already specified what it means for $a$ to be an indiscernible
for $\ga$ in a particular \pcs\ $X$.  In this section we will be
interested in sequences of indiscernibles.  Like the individual
indiscernibles these sequences will be defined for a particular \pcs\
$X$, but unlike the case for individual indiscernibles we will show that under
fairly general hypotheses the sequences of indiscernibles are
independent of the choice of $X$.

\paranumber{3}
In the last section we fixed $\gd$ to be 
cofinality of $\gk$, and each \pcs\ $X$ was assumed to be closed under
sequences of length~$\gd$.  Unless otherwise specified we will use
boldface letters to designate sequences of length~$\gd$, so that for
example we write $\vec a=\seq{a_\gi:\gi<\gd}$.

The ordering on sequences is by eventual dominance.  We will 
indicate\usage{$\geb$}\usage{$\lessb$}\usage{$\eqb$}
this by a subscript $b$ on the ordering relation: $\vec\gb'\geb\vec\gb$,
and $\vec\gb'\eqb\vec\gb$ mean respectively that for every
sufficiently large $\gi<\gd$ we have $\gb'_\gi\ge\gb_\gi$, or that for
every sufficiently large $\gi<\gd$ we have $\gb'_\gi=\gb_\gi$.  The
relation $\vec\gb'\not\geb\vec\gb$ means that it is not true that
$\vec\gb'\geb\vec\gb$, that is, that $\gb'_\gi<\gb_\gi$ for unboundedly
many $\gi<\gd$.

\paranumber{4}
\begin{defin}
\begin{enumerate}
\item
The sequence $\vec a$ is an {\it {\pris}\/} for $\vec\ga$  in $X$
if $\vec a$ and $\vec\ga$ are nondecreasing  sequences of length $\gd$
such that $\sup\vec a=\sup\vec\ga$, and
$a_\gi$ is an principle indiscernible 
for $\ga_\gi$ in $X$ for every sufficiently large $\gi<\gd$.
\notate{\pris}
\item\notate{$\vec a$ is a \bis\ for $\vec\ga$}
The sequence $\vec a$ is a {\it\bis\ for $\vec\ga$ in $X$\/} if $\vec a$ is a
\pris\ for $\vec\ga$ and
$\ga_\gi=i_{a_{\gi},\ga}$ (or, equivalently,
$\ga_{\gi}=\inf(i_{a_\gi,\gk},\ga)\,$) 
for all sufficiently large $\gi<\gd$,
where $\ga=\sup\vec\ga$.
\item
The sequence $\vec b$ is an {\it \is\ for $\vec\gb$
belonging to $(\vec a,\vec\ga)$\/}  in $X$ if $b_\gi$ is an indiscernible for
\notate{\is\ for $\vec\gb$ belonging to $\vec a$}
$\gb_\gi$ belonging to $(a_\gi,\ga_\gi)$  in $X$ for every sufficiently large
$\gi<\gd$.
\end{enumerate}
\end{defin}

\paranumber{5}
We will say that a sequence $\vec b$ is an \is\ for
$\vec\gb$ (without the qualifier ``in $X$'') if $\vec b$ is an \is\
for $\vec\gb$ in  every \pcs\ $X$ such that $\vec b\subset X$.
We similarly drop the qualifier ``in $X$'' from the definitions of
a \pris\ and of a \bis\ if the statement of definition is satisfied
for every \pcs\ $X$.
Most of the rest of this section will be concerned with proving
(using, in the case $\gd=\gw$,
an additional assumption  on the size of the members of the sequences 
$\vec b$ and $\vec \gb$)
that
we always can drop the qualifier ``in $X$'':
$\vec a$ is a \bis\ for $\vec\ga$, or $\vec b$ is a \is\ of $\vec\gb$
belonging to $(\vec a,\vec\ga)$, in a particular
\pcs\ $X$ if and only if
the same thing is true in any \pcs\ $Y$ containing the relevant
sequences.
For each  property $P$ of interest we will find 
a first order formula $\phi$ such that if $a$ is any member of a
\pcs\ $X$ then $P(a,X)$ holds if and  only if 
$X\sat\phi(a)$.  It follows that $P$ is independent of
the choice of \pcs, since if $Y$ is any other \pcs\ then $X$ and $Y$
are elementary substructures of $H_{\tau}$ and hence satisfy the
same formulas.

\paranumber{6}
In order to to avoid superscripts we will continue to 
work with a fixed \pcs\ $X$, but the formulas we obtain will not
depend on $X$.  

\medskip
\begin{lemma}\label{pris-def-lem}
There is a formula $\phi(\vec a, \vec\ga)$ which holds in $X$ if and
only if $\vec a$ is a \pris\ for $\vec\ga$ in $X$.
\end{lemma}

\begin{proof}
Let $\phi'(\vec a,\vec\ga)$ be the conjunction of the formulas
$$\aligned
&\exists h\in\core\exists \gi_0\forall \gi>\gi_0\;\;\ga_\gi\in h\image
a_\gi\\
&\forall h\in\core\exists\gi_0\forall\gi>\gi_0\;\;h\image
a_\gi\cap\ga_\gi\subset a_\gi 
\endaligned$$
 By lemma~\ref{ftoi-lem} the formula
$\phi'(\vec a,\vec\ga)$ holds if and only if $\vec a$ is a
\bis\ for $\vec\ga$.
Let $\ga=\sup\vec\ga=\sup\vec a$. 
If $\vec a$ is a \pris\ for $\vec\ga$ but not a \bis\ for
$\vec\ga$ then  both $\vec a$ and $\vec\ga$ are \bis s for the
sequence $\vec\ga'$ defined by
$\ga'_{\gi}=i_{a_\gi,\ga}(a_\gi)=i_{\ga_\gi,\ga}(\ga_{\gi})$.
Thus the following
formula 
$\phi(\vec a,\vec\ga)$  will satisfy the requirements of the lemma:

\begin{multline*}
\exists\vec\ga'\bigl(\phi'(\vec a,\vec\ga')\text{ and}\\
\text{if } I =
\set{\gi:\ga_\gi\not=\ga'_\gi}
\text{ is unbounded in $\gd$ then } \phi'(\vec\ga\restrict
I,\vec\ga'\restrict I)\bigr).
\end{multline*}

\end{proof}

\medskip
\begin{lemma}\label{main-lem}
{\rm(}Main Lemma\/{\rm)}
 There is a formula $\phi(\vec
a,\vec\ga,\vec b,\vec\gb)$ which holds if and only if 
\begin{enumerate}
\item 
$\vec a$ is a
\pris\ for $\vec\ga$
\item
$\vec b$ is an \is\ for $\vec\gb$ belonging
to $(\vec a,\vec\ga)$
\item
If $\gd=\gw$ then there is an integer $n$ such that 
$\vec\gb\lessb \left(\ga_i^{+n}:i\in\gw\right)$ and $\vec
b\lessb\seq{a_i^{+n}:i\in\gw}$.
\end{enumerate}
\end{lemma}

\paranumber{7}
Before proving this, we look briefly at the problem of determining,
given $\vec\gb$ and $(\vec a,\vec\ga)$, whether there exists an
indiscernible sequence $\vec b$ for $\vec\gb$ 
belonging to $(\vec a,\vec\ga)$.  The harder problem of determining
whether a particular sequence $\vec b$ is this indiscernible sequence
will be deferred until this problem
has been settled.

The easier problem breaks down into
two problems.  
The first, deciding whether $\vec \gb$ has an indiscernible sequence belonging
to $(\vec a,\vec\ga)$ at least in the weak sense that
$\gb_\gi=i_{a_\gi,\ga_\ga}(b_i)$, is answered rather easily by the next lemma,
assuming that $\vec \gb$ is not too large.  
The second question
is to determine whether the sequences satisfy the boundedness
conditions of definition~\ref{indforb-def}, that is, whether
\begin{gather}
\gb_{\gi}<\inf(C\setminus(\ga_\gi+1))\tag1\\
b_\gi<\max(\SO(a_\gi)^+,a_{\gi}^{++})\tag 2
\end{gather}
\noindent
hold
for all sufficiently large
$\gi<\gd$.   The bound~(2) for $b_\gi$ is 
quite straightforward, but the bound~(1) for $\gb_\gi$ is not possible to 
determine directly in $K$.  Most of the work involved in proving
lemma~\ref{main-lem} will come in the proof of
lemma~\ref{aux-main-lem} below, which 
uses the additional assumptions that $\vec\gb$ satisfies~(1)
and that $\vec a$ is a \bis.
Lemma~\ref{aux-main-lem} implies corollary~\ref{faga-exists}, which
implies among
other things that for \bis s the bound~(2) implies the bound~(1).
This proves the
main lemma for the case of \bis s,
and the general case follows easily from this special case.

\begin{lemma} 
There is a formula $\phi$ such that  if
$\ga_{\gi}\le\gb_\gi<\inf(C\setminus(\ga_\gi+1))$ for all $\gi<\gd$ then
$\phi(\vec a,\vec\ga,\vec\gb)$
holds in $X$ if and only if $\vec a$ is a \bis\ for $\vec\ga$ and there is
a sequence $\vec b$ such that $\gb_\gi=i_{a_\gi,\ga_\gi}(b_\gi)$ for
all sufficiently large $\gi<\gd$.
\end{lemma}

\begin{proof}
First note that  if $\gb_\gi=i_{a_\gi,\ga_\gi}(b_\gi)$ then
since $\vec a$ is basic
there is a
function $f\in X\cap\core$ such that 
$\gb_\gi\in f\image\ga_\gi$ for all sufficiently large $\gi<\gd$.  In
the case $\ga_\gi=i_{a_\gi,\gk}(\ga_\gi)<\ga=\sup\vec\ga$ this is true by
elementarity, since 
$h_\gk^X\cap (\gk\times\gk)$ is a member of $\core$
satisfying the stated property.  In the case
$\ga_\gi=\ga$ we have $\gb_{\gi}\in h_{\ga}^X\image a_{\gi}$ for all
sufficiently large $\gi<\gd$, and by lemma~\ref{ftoi-lem} there is a
function $f\in K\cap X$ such that $f(\xi)=h_{\ga}^X(\xi)$ for all
$\xi\in h_{\ga}^{-1}\image\vec \gb$.
On the other hand the existence of such a function $f$ implies that
$\vec\gb$ has an indiscernible sequence $\vec b$: for sufficiently
large $\gi\in\gd$ we have $f\in\range(i_{a_{\gi},\ga_\gi})$, and we
set $b_{\gi}=i_{a_\gi,\ga_\gi}^{-1}(f)(f^{-1}(\gb_\gi))$.
Thus if we let $\phi$ be 
the conjunction of the formula
$$\exists f\in
K\exists\gi_0<\gd\forall\gi\,(\gi_0<\gi<\gd\implies\gb_\gi\in f\image
a_\gi)$$ 
with the formula asserting that $\vec a$ a \bis\ for
$\vec\ga$ then $\phi$ satisfies the requirements of the lemma.
\end{proof}

\paranumber{8}
\begin{lemma}\label{aux-main-lem}
There is a formula $\phi$ such that if
$\gb_\gi<\inf(C\setminus(\ga_\gi+\nobreak 1))$ for all sufficiently large
$\gi<\gd$ then $\phi(\vec a,\vec\ga,\vec b,\vec\gb)$ holds if and only
if clauses~(1--3) of lemma~\ref{main-lem} are satisfied. 
\end{lemma}

\begin{cor}\label{faga-exists}
If $\vec a$ is an \is\ for $\vec\ga$ and $\ga_\gi=\sup\vec a$ for all
sufficiently large $\gi<\gd$
then $F_{a_\gi,\ga_\gi}$ is in $X$ for all sufficiently large $\gi<\gd$.
\end{cor}

Before starting on the proof of lemma~\ref{aux-main-lem} and
corollary~\ref{faga-exists} we will verify that together they imply
the main lemma, lemma~\ref{main-lem}.

\proof[Proof of main lemma from \ref{aux-main-lem} and \ref{faga-exists}]
Suppose first that $\vec a$ is a \bis\ for $\vec\ga$.  Since
$F_{a_\gi,\ga_\gi}$  is always in $X$ when
$\ga_{\gi}=i_{a_\gi,\ga}(a_\gi)<\ga=\sup\vec\ga$, the corollary
implies that 
$F_{a_\gi,\ga_{\gi}}$ is in $X$ for all sufficiently large $\gi<\gd$.
Since $\len(F_{a_\gi,\ga_\gi})<\inx(F_{a_\gi,\ga_\gi})\in\OO(\ga_\gi)$ 
it follows that if $\vec b$ is an \is\ for $\vec\gb$
belonging to $(\vec a,\vec\ga)$
then $\gb_\gi<\SO(\ga_\gi)$ for almost all $\gi<\gd$.
Thus the formula of lemma~\ref{aux-main-lem} satisfies
the requirements of the main lemma whenever $\vec a$ is an \bis\ for
$\vec\ga$. 
\paranumber{9}
Now we can treat the general case as in the proof of
lemma~\ref{pris-def-lem}. Let $\ga'_{\gi}=i_{\ga_\gi,\ga}(\ga_\gi)$
and $\gb'_{\gi}=i_{\ga_\gi,\ga}(\gb_\gi)$.   Then $\vec a$ is an
\pris\ for
$\vec\ga'$ and $\vec b$ is a \is\ for $\vec\gb'$ belonging to $(\vec
a,\vec\ga')$, and if
$I\eqdef\set{\gi<\gd:\ga_\gi\not=\ga'_\gi}$ is unbounded in $\vec\gd$
then the restriction $\ga\restrict I$  of $\ga$ is a \pris\ for
$\vec\ga'\restrict I$, and the restriction
$\vec\gb'\restrict I$ of $\vec\gb'$ is an \is\ for $\vec\gb\restrict
I$ belonging to 
$(\vec\ga\restrict I,\vec\ga'\restrict I)$.  This completes the proof
of lemma~\ref{main-lem}, assuming lemma~\ref{aux-main-lem} and
corollary~\ref{faga-exists}.
\endproof
The proof of lemma~\ref{aux-main-lem} will be broken into two cases,
the first for $\gd>\gw$ and the second for $\gd=\gw$.

\paranumber{10}
\proof[Proof of lemma~\ref{aux-main-lem} for $\gd>\gw$]
In this case
we use the game introduced by Gitik in \cite{gitik.mble-cards-not-ch}
to obtain a rather straightforward extension of the results which were
proved for measure sequences in
\cite{mitchell.sing-card-hyp,mitchell.applications,mitchell.defsch}.
The major difference is that the restriction to $\gd>\gw$ means that
where  the  results given in those papers for sequences of measures
have finite sets
of exceptional points, the results given in this paper for sequences
of extenders  may have a
countable set of exceptional points.

\paranumber{11}
\newcommand\game{\operatorname{\mathcal G}}
\newcommand\gameb[1]{\game(#1,\vec\gb)}
\newcommand\gamebb{\gameb{\vec b}}
Our presentation of Gitik's  game will differ somewhat from that of
\cite{gitik.mble-cards-not-ch}.
We will define a game $\gamebb$ between two
players, who, following Mathias, we call Adam and Eve. 
The first player, Adam, will be trying
to show 
that the ordinals in
$\vec b$ are  too small for $\vec b$ be an \is\ belonging to
$\vec\gb$.  He will 
do so by proposing sets of ordinals $B_{n,\gi}\subset \gb_{\gi}$.  Eve
will be required to defend the proposition that $\vec b$ is an \pris\
for $\vec \gb$ by choosing indiscernibles for ordinals in the sets
$\seq{B_{n,\gi}:\gi<\gd}$ which are consistent with $\vec b$ and with her
earlier choices.

\paranumber{12}
In the next two propositions
 we will show that if $\vec b$ is an \is\ for $\vec\gb$ then Eve has a
winning strategy for the game $\gamebb$, while Adam has a winning
strategy for the game $\gameb{\vec b'}$ whenever
$\vec b'\not\geb\vec b$.
With these propositions we can complete the proof  of
lemma~\ref{aux-main-lem}, since the \pris\ $\vec b$ for $\vec\gb$ is
definable by a formula~$\phi$  asserting that $\gb$ is 
the least sequence $\vec b'$ such that Eve has a winning
strategy for $\gameb{\vec b'}$.
The definition of the game $\gamebb$, and hence the formula~$\phi$,
will not depend in any way on the particular
\pcs\ $X$.

\paranumber{13}
\begin{defin}\label{gamedef}
If $\vec a$ is a \bis\ for $\vec\ga$ then the game
$\game(\vec\gb,\vec b)$ is defined as follows:

The first player, Adam, plays on his $n$th move a sequence
$\seq{B_{n,\gi}:\gi<\gd}$ such that
$B_{n,\gi}\in\left[\gb_\gi\setminus\ga_\gi\right]^{\le\gd}$ for each
$\gi<\gd$.
The second player, Eve,  responds with a sequence
$\seq{\tau_{n,\gi}:\gi<\gd}$ such that 
\begin{enumerate}
\item  For each
$\gi<\gd$ the function $\tau_{n,\gi}$ is an order preserving function
mapping a subset of 
$\gb_\gi\setminus\ga_\gi$ into $b_\gi\setminus a_\gi$. 
\item If $h$ is any function in $K$ then  $B_{n,\gi}\cap
h\image a_\gi\subset\domain\tau_{n,\gi}$ for all but boundedly many
$\gi<\gd$.
\item $\tau_{n,\gi}\supset\tau_{n-1,\gi}$ for all $n>0$.
\end{enumerate}

\paranumber{14}
Adam wins if Eve is ever unable to play; otherwise Eve wins.
\end{defin}

The idea is that if $\gb\in B_{n,\gi}$ then $\tau_{n,\gi}(\gb)$ should
be an indiscernible for $(a_i,\ga_i)$ belonging to $\gb_i$.
For convenience, we will write $j_\gi$ for $i_{a_\gi,\ga_\gi}$, and if
$\vec\xi$ is a sequence then we will write $j(\vec\xi)$ for
$\seq{j_\gi(\xi_\gi):\gi<\gd}$.

\paranumber{15}
\begin{prop}\label{EveWins}
If $j(\vec b)\geb\vec\gb$ then Eve has a winning strategy for the
game $\gamebb$.
\end{prop}

\begin{proof}
Suppose that the proposition is false.  Since the game   $\gamebb$ is
closed it is determined, and hence
Adam must have a winning
strategy.     By the elementarity of $X$ there is a winning strategy
$\gs\in X$ for Adam.  Now suppose Eve plays, in $V$,  against the
strategy $\gs$ by playing 
$\tau_{n,\gi}=j_{\gi}^{-1}\restrict\left(B_{n,\gi}\cap j\image
b_\gi\right)$.
It is easy to see that these plays by Eve satisfy the 
first and last clauses of definition~\ref{gamedef}, 
and the second clause follows from clause~(4) of
lemma~\ref{ftoi-lem}.  Thus Adam loses this game, contradicting the
assumption that $\gs$ is a winning strategy for Adam.
\end{proof}

\begin{prop}\label{AdamWins}
If $\cof(\gd)>\gw$ and $j(\vec b)\not\geb\vec\gb$ then Adam has
a winning strategy for the game $\gamebb$.
\end{prop}

\begin{proof}
As before it will be sufficient to show that if $\gs$ is any strategy
for  Eve which is a member of $X$ then Adam can refute the strategy
$\gs$ by playing in~$V$.  Adam will let $B_{n,\gi}=\nothing$ whenever
$j_\gi(b_\gi)\ge\gb_\gi$, so we can assume that
$j_\gi(b_\gi)<\gb_\gi$ for all $\gi<\gd$.
The first move in Adam's refutation will be
the singleton sets $B_{1,\gi}\eqdef\{j_\gi(b_\gi)\}$.  
His $n$th move, for $n>1$, will be the sets $B_{n,\gi}\eqdef(\gb_\gi\setminus
\ga_\gi)\cap j_{\gi}\image(\range(\tau_{n-1,\gi}))$, where the
functions $\tau_{n-1,\gi}$ are taken from Eve's previous
move.
We need to show that if Eve plays by the strategy~$\gs$ then Adam wins
this play of the game.

\paranumber{16}
First we observe that all of the plays in this game are members of $X$.
For Adam's moves it is sufficient to show that the sets $B_{n,\gi}$
are subsets of $X$, since his moves have cardinality $\gd$ and
$X$ is $\gd$-closed.  The sets $B_{1,\gi}$ are contained in $X$ by
construction,  and for $n>1$ the sets 
$B_{n,\gi}$ will be contained in $X$ provided that Eve's $n-1$st move is in
$X$.  But  Eve's strategy $\gs$ and the game $\gamebb$ are
both members of $X$, so Eve's moves will be in $X$ because Adam's
preceding moves were in $X$.

\paranumber{17}
Now let $\ga=\sup\vec\ga$.
We claim that $B_{n,\gi}\subset h_\ga\image a_{\gi}$ for each
$\gi<\gd$, where $h_\ga$ is the function defined in
definition~\ref{ih-def}(2).
  First, we have $B_{n,\gi}\subset h_{\ga_\gi}\image a_{\gi}$
for each $\gi<\gd$ since $a_0>\rho^X$ and $B_{n,\gi}\subset
\image a_{\gi}=i_{a_\gi,\ga_\gi}\image a_{\gi}$.  If
$\ga_\gi=\ga$ then the claim is established, and if $\ga_\gi<\ga$ then
the assumption that $\vec a$ is a \bis\ for $\vec\ga$ implies that
$i_{\ga_\gi,\ga}(\ga_\gi)=\ga_{\gi}$ and it follows by clause~(1) of
lemma~\ref{ftoi-lem} that 
$h_{\ga}\circ h_{\ga_{\gi}}^{-1}$ is the identity.
It follows by 
lemma~\ref{ftoi-lem}(5) that there is a function $h\in K$ such that
$\union_nB_{n,\gi}\subset h\image a_{\gi}$ for all $\gi<\gd$.

\paranumber{19}
By clause~(2) of the definition  of $\gamebb$ it follows that for each
$n\in\gw$ there is $\gi_n<\gd$ such that
$B_{n,\gi}\subset\domain(\tau_{n,\gi})$ for all $\gi>\gi_n$.  Since
$\gd=\cof(\gd)>\gw$ and $j(\vec b)\not\geb\vec\gb$ there is an ordinal
$\gi<\gd$ such that $\gi>\gi_n$ for all $n\in\gw$.
We are now ready to reach the
contradiction and hence complete the proof of lemma~\ref{AdamWins}.
Define an infinite descending $\gw$-sequence $\vec\eta$ of 
ordinals by setting
$\eta_0=j_\gi(b_\gi)$ and 
$\eta_{n}=j_{\gi}\circ\tau_{n-1,\gi}(\eta_{n-1})$ for $n>0$.
Since $\eta_n\in B_{n,\gi}\subset\domain\tau_{n,\gi}$ the ordinal
$\eta_{n+1}$ is defined for all $n<\gw$.
We have $\eta_0<\gb_\gi$ by assumption, and an easy proof by induction,
using  clauses~(2) and ~(3) of the definition of the game $\gamebb$,
shows that $\eta_{n+1}<\eta_n$ for all $n<\gw$.
This contradiction completes the proof of
the proposition.
\end{proof}
Now let $\phi(\vec a,\vec\ga,\vec b,\vec\gb)$ be the conjunction of
three formulas, asserting 
\begin{enumerate}
\item
$\vec\gb$ has an indiscernible sequence belonging to $(\vec
a,\vec\ga)$.
\item
Eve wins the game $\gamebb$.
\item
Adam wins the game $\gameb{\vec b'}$ for all $\vec b'\not\geb\vec b$. 
\end{enumerate}
Lemmas~\ref{EveWins} and~\ref{AdamWins} imply that $\phi(\vec
a,\vec\ga,\vec b,\vec\gb)$ is true whenever $\vec b$ is an
indiscernible sequence for $\vec\gb$.  On the other hand, if $\vec b'$
is any sequence such that  $\phi(\vec
a,\vec\ga,\vec b',\vec\gb)$ is true then $\vec \gb$ has an
indiscernible sequence $\vec b$ by clause~(1), so  $\phi(\vec
a,\vec\ga,\vec b,\vec\gb)$ is  true as well. From clauses~(2) and~(3)
it 
follows that $\vec b\leb\vec b'$ and $\vec b'\leb\vec b$, so $\vec
b'\eqb\vec b$.  This completes the proof of the case $\gd>\gw$ of
lemma~\ref{aux-main-lem}. 
\qed

\paranumber{20}
\begin{proof}[Proof of lemma~\ref{aux-main-lem} for $\gd=\gw$]
When $\gd=\gw$ the situation becomes much more difficult, and in this
case we only know how to reconstruct the embeddings under the
assumption that there is no inner model of $\exists\ga\, o(\ga)=\ga^{+\gw}$.

Define $\ga_{n,k}$,  for
integers $n$ and $k$, to be the smaller of  $\ga^{+n}_k$ and
$i_{a_k,\ga_k}(a_k^{+n})$, and let $a_{n,k}$ be
$i_{a_k,\ga_k}^{-1}(\ga_{n,k})$, provided that it exists. 
For each $n\in\gw$  we claim that $a_{n,k}$ exists for all but
finitely many $k\in\gw$. This is immediate if
$\ga_k=\ga=\sup\vec\ga$ for all sufficiently large $k<\gw$, since any
member of $X\cap\inf(C\setminus\nobreak(\ga+\nobreak 1))$ is in 
the range of $i_{\nu,\ga}$
for all sufficiently large $\nu<\ga$.  Now suppose that
$\ga_k=i_{a_k,\ga_k}(a_k)=i_{a_k,\ga}(a_k)$, and note that if $a_{n,k}$ does
not exist then $\ga_k^{+n}<i_{a_k,\ga_k}(a_k^{+n})$.  Now
$\pi^{-1}(\ga^{+n}_k)$ is the $n$th successor of 
$\pi^{-1}(\ga_{n})$ in $\barmse m_{\bar\ga_k}$, since
$\bare_{\bar\ga_k}$ does not exist, 
and $\pi^{-1}(i_{a_k,\ga_k}(a_k^{+n}))$  is the $n$th  successor of
$\pi^{-1}(\ga_{k}^{+n})$ in $\ult(\barmse m_{\bar\ga_k},\barf_{\bar
a_k,\bar\ga_k})$.  Since  $\ult(\barmse m_{\bar\ga_k},\barf_{\bar
a_k,\bar\ga_k})$ is smaller than $\barmse m_{\bar\ga_k}$ it follows
that $\pi^{-1}(i_{a_k,\ga_k}(a_k^{+n}))\le\pi^{-1}(\ga_{n})$ so that 
$i_{a_k,\ga_k}(a_k^{+n})\le\ga_{k}^{+n}$, and hence $a_{n,k}$ does
exist.

\paranumber{21}
Note that if $n\in\gw$ then the sequence $\vec a_n\eqdef\seq{a_{n,k}:k\in\gw}$ 
 is an \is\  for
$\vec\ga_n\eqdef\seq{\ga_{n,k}:k\in\gw}$ belonging to $(\vec
a,\vec\ga)$. 

We prove lemma~\ref{aux-main-lem} by induction on $n$, with the
 induction step relying on the following lemma.  
Since $\ga_{1,k}$ is always equal to $\ga^+$ the case $n=1$ could be
handled by standard methods, but for convenience we treat it as
part of the general induction.

\paranumber{22}

\begin{lemma}\label{ctble-ind-step}
There is a first order formula $\psi$ with the following property:
Assume that $X$ is a \pcs,
$n$ is an integer and the sequence $\vec a_n$ and $\vec\ga_n$ are as
defined above.  Let $T\in X$ be a set such that  if $n=0$ then
$T=\set{(\vec\gd,\vec\gd):\vec\gd<\vec\ga}$, and if $n>0$ then
$$
T\cap X=\set{(\vec d,\vec \gd):\vec\gd\lessb\vec\ga_n\text{ and }
\vec d\text{ is an \is\
for }\vec\gd\text{ in }X}.$$
 Then for all sequences $\vec\gb$ and $\vec
b$ in $X$ the formula $\psi(T,\vec a_n,\vec\ga_n,\vec b,\vec\gb)$
is true
if and only if 
$\vec\gb\lessb\vec\ga_{n+1}$
and $\vec b$ is an \is\ for $\vec\gb$.

\end{lemma}

\paranumber{23}
\def\iso(#1,#2,#3){f_{#1,#2}(#3)}
\begin{proof}
If $\ga$ and $\gb$ are any two ordinals such that
$\card\gb^{\core}\le\ga<\gb$ 
then let $f_{\ga,\gb}$ be 
the least map $f$ in the natural ordering of
$\core$ such that $f\colon\ga\cong\gb$. 
 We will define a  function $S(\gg,\gb, \xi)$ by recursion on $\gb$.
The domain of $S$ is the set of triples of ordinals
$\xi,\gg$ and $\gb$
such that $\xi<\gg<\gb$ and 
$\gg\ge\card{\gb}^{\core}$,  and $S$ is defined recursively as
follows:
$$S(\gg,\gb,\xi)=\cases
0&\text{if $\iso(\gg,\gb,\xi)<\gg$}\\
S\bigl(\gg,\iso(\gg,\gb,\xi),\xi\bigr)+1&\text{if $\gg\le
f_{\gg,\gb}(\xi)<\gb$.} 
\endcases$$
Now let  $\psi(T,\vec a_n,\vec\ga_n,\vec b,\vec\gb)$ be the conjunction
of the following three formulas:
\begin{align}
&\exists k_0\forall
k>k_0\quad(\card{b_k}<a_{n,k}\land\card{\gb_k}<\ga_{n,k})\tag {i}\\
&\exists g\in\core\exists k_0\forall k>k_0\quad\gb_k\in g\image a_k \tag{ii}\\
&\forall(\vec d,\vec\gd)\in T\,\exists k_0\forall k>k_0\quad
S(\ga_{n,k},\gb_k,\gd_k)=S(a_{n,k},b_k,d_k).\tag{iii}
\end{align}

\paranumber{24}
It is clear that $\psi(T,\vec a_n,\vec\ga_n,\vec b,\vec\gb)$ holds in
$X$ whenever $\vec b$ is an \is\ for $\vec\gb$ in $X$.  Now suppose
that 
\begin{equation}
\psi(T,\vec a_n,\vec\ga_n,\vec b',\vec\gb)
\tag{1}
\end{equation}
 is true in $X$ for
some sequence $\vec b'\not\eqb\vec b$.  If we set $\vec\gb'=j(\vec b')$
then $\vec b'$ is an \is\ for $\vec\gb'$ belonging to $(\vec
a,\vec\ga)$ and hence
\begin{equation}
\psi(T,\vec a_n,\vec\ga_n,\vec b',\vec\gb')
\tag{2}
\end{equation}
 is also true in $X$.  We will show that (1) and~(2) lead to a
contradiction.  We can assume \iwlog\ that 
$\gb'_k<\gb_k$ for unboundedly many $k<\gw$.  
For each such $k$ set 
$\gd_k=f^{-1}_{\ga_{n,k},\,\gb_n}(\gb'_k)$, so that $\gd_k<\ga_{n,k}$.
The sequence $\vec\gd$ has an
\is\ $\vec d$ belonging to $(\vec a,\vec\ga)$, since $\gd_k$
is defined in $\core$ from the parameters $\ga_{n,k}$,  $\gb_k$ and
$\gb'_k$ and each of the sequences $\vec\ga_n$, $\vec\gb$ and
$\vec\gb'$ has an \is\ in $X$.   Then $S(\ga_{n,k},\gb_i,\gd_k)\not=0$
since $f_{\ga_{n,k}}(\gd_k)=\gb'_{k}>\ga_{n,k}$, so for sufficiently
large $k$ such that $\gb'_k<\gb_k$ 
\begin{align*}
S(a_{n,k},b'_{n,k},d_k)&=
S(\ga_{n,k},\gb_{k},\gd_k)&&\text{by (1)}\\
&=S(\ga_{n,k},\gb'_k,\gd_k)+1&&\text{by the choice of
$\operatorname{S}$ and $\gd_\gk$}\\ 
&= S(a_{n,k},b'_{k},d_k)+1&&\text{by (2).}
\end{align*}
This
contradiction
completes the proof of lemma~\ref{ctble-ind-step}.
\end{proof}

\paranumber{25}
Lemma~\ref{aux-main-lem} will follow easily from
lemma~\ref{ctble-ind-step} once we verify
that it is possible to define the
sequence $\vec\ga_{n+1}$ and its \is\ $\vec a_{n+1}$. 
This is straightforward: $\vec\ga_{n+1}$ is the minimal sequence which has
an indiscernible sequence but does not have an indiscernible sequence
satisfying $\psi$, that is, $\vec\ga_{n+1}$ is the only sequence
$\vec\ga'=\seq{\ga'_k:k<\gw}$, up to bounded segments, which satisfies
the conjunction of the following three formulas:
\begin{align*}
&\exists g\exists k_0\forall k>k_0\;\ga'_{k}\in g\image (a_k)\\
&\lnot\exists\vec a'\;\psi(T,\vec a_n,\vec\ga_n,\vec a',\vec\ga')\\
&\forall\vec\gb\forall g\in\core\;(\text{if }
  I=\set{k:\gb_k<\ga'_{k}\land\gb_k\in g\image a_k}
  \text{ is infinite }\\
&
\qquad\qquad\qquad\qquad\qquad
\text{then }\exists\vec b\;
  \psi(T,\vec a_n\restrict I,\vec\ga_n\restrict I,\vec b\restrict
  I,\vec\gb\restrict I). 
\end{align*}

Similarly $\vec a_n$ is the minimal sequence $\vec a'$
which is not an \is\ for any sequence $\vec\ga'$
satisfying the formula $\psi$ and is hence is definable up to an initial
segment.

\paranumber{26}
This completes the proof of lemma~\ref{aux-main-lem}.
\end{proof}

\newcommand\gve{\varepsilon}
\proof[Proof of corollary~\ref{faga-exists}]
The hypothesis of corollary~\ref{faga-exists} asserts that 
$\vec a$ is an \is\ for $\vec\ga$ in $X$ 
such that $\ga_\gi=\ga=\sup\vec\ga$ for sufficiently large $\gi<\gd$,
and the conclusion asserts that  $F_{a_\gi,\ga}\in X$ for sufficiently large
$\gi<\gd$.  If the hypothesis is true and the conclusion is false
then we can assume without loss of generality that
$F_{a_\gi,\ga}\notin X$ for all $\gi<\gd$. This means that
$\barf_{\bar a_\ga,\bar\ga}\notin \bark$, so that
either $\ga=\gk$ or $i_{\ga,\gk}(\ga)>\ga$,  and in either case
$\bare_{\bar\ga}\trieq\barf_{\bar a,\bar\ga}$. 

We will define, in $X$, a set $G$ such that
$\pi\image\bare_{\bar\ga}=G\cap X$.  It will follow that
$\bare_{\bar\ga}=\pi^{-1}(G)\in N$, so $\bare_{\bar\ga}\in \bark$,
contradicting the choice of $\bare_{\bar\ga}$ as the least extender in
$\barmse m_{\bar\ga}$ which is not in $\bark$.
\paranumber{27}
In order to define $G$ we need to decide inside $X$ whether a pair
$(\gve,z)$ is in $\pi\image\bare_{\bar\ga}$.  Now notice that if we
set $\gg=\pi\left(\inx\left(\bare_{\bar\ga}\right)\right)$ then
lemma~\ref{aux-main-lem} 
implies that  for
each ordinal $\gb$ with $\ga\le\gb\le\gg$ 
 there is an indiscernible sequence $\vec b$ for the constant sequence
$\gb$
belonging to $(\vec a,\ga)$, and there is a formula $\phi$ picking out
these pairs $(\vec b,\gb)$.  In order to use the lemma we have to
check that $\gb$ is less than the least member of $C$ above $\ga$, but
this is immediate:  Since $\bare_{\bar\ga}$ is the extender used on $\barmse
m_{\bar\ga}$ and the models do not have overlapping extenders, the
critical points of 
extenders used later will be greater than the index $\pi^{-1}(\gg)$ of
$\bare_{\bar\ga}$. 

Thus we can choose indiscernible sequences $\vec
e$ and $\vec c$ belonging to $(\vec a,\ga)$ for the constant
sequences $\gve$ and $\gg$ respectively.  For sufficiently large
$\gi\in D$ we will have $i_{a_\gi,\ga}(e_\gi,c_\gi)=(\gve,\gg)$ and
$z\in\range(i_{e_{\gi},\ga})$, and for all such ordinals~$\gi$ we will
have $(\gve,z)\in \pi\image\bare_{\bar\ga}$ if and only if 
\begin{align*}
e_\gi\in z&\qquad\text{if $\bare_{\bar\ga}=\barf_{\bar
a_\gi,\bar\ga}$}\\ (e_\gi,z\cap a_{\gi})\in\ce_{c_\gi}&\qquad\text{if
$\bare_{\bar\ga}\tri\barf_{\bar a_\gi,\bar\ga}$.}
\end{align*}
Since the  indiscernible sequences can be defined inside $X$ and
$^{\gd}X\subset X$, this definition of $G$ can be carried out in $X$.
This completes the proof of corollary~\ref{faga-exists} and of the
main lemma.
\endproof
\paranumber{28}
\bigskip
For the rest of the section we will assume that if $\gd=\gw$ then 
$\set{\nu<\gk:o(\ga)\ge\ga^{+n}}$ is bounded in $\gk$ for some $n<\gw$.
The next task
is to extend our notion of \is\ to sequences of indiscernibles for
particular extenders.
\begin{defin}
\begin{enumerate} 
\item
The ordinal $a$ is a {\it\pri\ in $X$ for the extender $F$} on $\ga$
if $a$ is a  \pri\ for $\ga$ in $X$ and $F=F^X_{a,\ga}$.
\item
The sequence $\vec a$ is a {\it\pris\ in $X$ for the sequence $\vec F$} of
extenders if $a_\gi$ is a \pri\ in $X$ for $F_\gi$ for every
sufficiently large $\gi<\gd$.
\item
The sequence $\vec a$ is a {\it \pris\ for the sequence $\vec F$ of
extenders}  if for  every \pcs\
$X$ containing $\vec a$, the sequence $\vec a$ is a \pris\ in $X$ for
$\vec F$.
\end{enumerate}

\end{defin}

\paranumber{28a}
\begin{lemma} \label{isForExtCor}
There is a formula $\psi$ such that  $\psi(\vec
a, \vec\ga,\vec F)$ is true in a \pcs\ $X$
 if and only if $\vec a$ is a \pris\ in $X$ for the sequence $\vec F$ of
extenders on $\vec\ga$.  Thus  if $\vec a$ is a \pris\ for $\vec F$ in
some \pcs\ $X$
then it is a \pris\ for $\vec F$.
\end{lemma}

\paranumber{28b}
\begin{proof} 
By lemma~\ref{main-lem} there is a first order formula $\psi(\vec
a,\vec\ga,\vec F)$
over $X$ asserting that the following statements are true.
We write $\gg_\gi$ for $\inx(F_{\gi})$ and $\vec\gg$ for
$\seq{\gg_{\gi}:\gi<\gd}$. 

\begin{enumerate}
\item
$\vec\ga=\seq{\crit(F_\gi):\gi<\gd}$, and $\vec a$ is a \pris\ for
$\vec\ga$. 
\item 
There is an indiscernible sequence $\vec c$ for  $\vec\gg$ belonging to
$(\vec a,\vec\ga)$, and $c_\gi\notin \OO(a_{\gi})$ for sufficiently
large $\gi<\gd$.
\item
If $\vec c'$ and $\vec \gg'$ are any other sequences
such that $\vec c'$ is
an \is\ for $\vec\gg'$, then, with at most boundedly many exceptions, 
$c_{\gi}\in\OO(a_{\gi})$ for all $\gi$ such that 
$\gg'_{\gi}\not\tri\gg_{\gi}$.
\item
If $f$ is any function in $K$, $\vec\gve\lessb\vec\gg$, and 
 and $\vec e$ is an \is\ for $\vec\gve$
belonging to $(\vec a,\vec\ga)$ then there is an ordinal $\gi_0<\gd$
such that for all $\gi_0<\gi<\gd$ and all $z\in f\image a_{\gi}$ we
have $z\in (F_{\gi})_{\gve_\gi}$ if and only if $e_\gi\in z$.
\end{enumerate}

If $\vec a$ is an \is\ for $\vec F$ in $X$ then $\psi(\vec a,\vec\ga,\vec F)$
will be true in $X$, and hence in $V$.  
Now we will show that if 
$\vec F'$ is any 
sequence of extenders in $K$ such that $\psi(\vec a,\vec\ga,\vec F')$ then
$\vec F'\eqb\vec F$.  By clauses~(2) and~(3) it is enough to show that,
with at most boundedly many exceptions, one of $F'_\gi$ and $F_{\gi}$ is
an initial segment of the other.

If neither of $F_\gi$ and $F'_\gi$ is an initial segment of the other
then let $(z_\gi,\gve_\gi)$ be the least pair such that 
$z_\gi\in (F_{\gi})_{\gve_\gi}\iff  z_\gi\notin(F'_{\gi})_{\gve_\gi}$.
Since there are indiscernible sequences for $\vec\gg$ and $\vec\gg'$,
there is an indiscernible sequence $\vec e$ for $\vec\gve$, but then
clause~(4) cannot be true for both $\vec F$ and $\vec F'$.
\end{proof}

\paranumber{28c}
We are now able to define the version of the ``next indiscernible''
function which is appropriate to the sequences which we are
considering.   We will define three separate functions:
The function $s^X$
gives the next principle indiscernible, $a^X_{\xi}$ 
gives the $\xi$th-next accumulation point, and
$\indiscbelong^X$  gives the indiscernible  for an ordinal $b$
belonging to a pair
$(a,\ga)$.  Definition~\ref{NextIndiscDef} below
has the formal definitions for these functions, together with
another function $\lis^X$ which is a useful variant of $s^X$.  
 The definition is relative to a
particular \pcs\ $X$, but we will finish up this section by showing
that any two \pcs s $X$ and $X'$ agree on the values of these functions
for all sufficiently large ordinals $\nu\in X\cap X'$ below $\gk$.

Recall that we use $\gg\tri\gg'$ to mean that either
$\ce_{\gg}\tri\ce_{\gg'}$ or $\gg\in\OO(\ga)$ and $\gg'=\SPO(\ga)$.

\paranumber{29}
\begin{defin}\label{NextIndiscDef}
If $X$ is a \pcs\ then
\begin{enumerate}
\item \usage{$s^X(\gg,\nu)$}
$s^X(\gg,\nu)$ is the least ordinal $a>\nu$ such that $a$ is a
\pri\ in $X$ for an ordinal $\ga$ such that $\ce_\gg=F^X_{a,\ga}$.
\item \usage{$\lis^X(\gg,\nu)$}
$\lis^X(\gg,\nu)$ is the least ordinal $a$, with $\gg>a\ge\nu$, such that
$a=s^X(\gg',\nu)$ for some $\gg'\trige\gg$.
\item \usage{accumulation point in $X$}
$a$ is an {\it accumulation point in $X$\/} for $\gg$ if
$\ga<\gg\in\PO(\ga)$, where either $\ga=a$ or $a$ is a \pri\ for $\ga$
in $X$, and $\lis^X(\gg',\nu)<a$ for every ordinal $\nu\in a\cap X$ and
every $\gg'\tri\gg$ in $h^X_\ga\image a\cap \OO(\ga)$.
\item \usage{$a^X_\xi(\gg,\nu)$}
If $\xi<\gw_1$ then $a^X_\xi(\gg,\nu)$ is the $\xi$th accumulation point for
$\gg$ above $\nu$.  We write $a^X(\gg,\nu)$ for $a^X_1(\gg,\nu)$.
\item\usage{$\indiscbelong^X(\gb,a,\ga)$}
$\indiscbelong^X(\gb,a,\ga)$ is equal to the ordinal $b$, if there
is one, such that $b$ is an indiscernible in $X$ for $\gb$ belonging to
$(a,\ga)$.
\end{enumerate}

\end{defin}

\paranumber{30}
Notice that if $a$ is a \pri\ for $\ga$ then there is a $\tri$-largest
ordinal $\eta\in\PO(\ga)$ such that $a$ is an 
accumulation point for $\eta$, and that the set of accumulation points for
an ordinal $\eta$ is closed in $X$.

\paranumber{31}
\begin{lemma}
\label{semi-cont}
Suppose that $Y$ is a \pcs\ and that $\ga\in Y$ has cofinality $\gd$.
Then for all but boundedly many $\nu<\ga$, if
 $\inf(\gd,\gw_1)\le\cof(\nu)\le\gd$ and
$\nu$ is an accumulation point in $Y$ for some 
$\eta\in\PO(\ga)$ then there is $\eta'\ge\eta$ and
$\gg<\nu$ such that $\nu=\lis^Y(\eta',\gg)$.
\end{lemma}
\begin{proof}
By lemma~\ref{faga-exists}, for all sufficiently large ordinals
$\xi<\ga$ which are 
principle indiscernibles for $\ga$ in $Y$, there is
$\eta_\xi\in\OO(\ga)$ such that $\xi$ is an indiscernible for 
$\ce_{\eta_\xi}$ in $Y$.  Let $\nu$ be as in the hypothesis so that
$\nu$ is a principle indiscernible for $\ce_{\eta_\nu}$.   Using
lemma~\ref{faga-exists} again, all but boundedly many of the principle
indiscernibles $a$ for $\ga$ below $\nu$ are indiscernibles for some extender
on $a$. For all such $a$ we have $\eta_a\tri\eta_{\nu}$, so that
$\nu=s^{Y}(\eta_{\nu},\gg)$ for some $\gg<\nu$. 
\end{proof}

\paranumber{32}
\begin{lemma}\label{next-ind-fctn}
For all \pcs s $X$ and $X'$ there is an ordinal $\eta<\gk$ such that
if $\eta<\nu<\ga\le\gk$, $\xi\le\gd$,  and $\ga<\gg\in\PO(\ga)$ then
\begin{align*}
s^X(\gg,\nu)&=s^{X'}(\gg,\nu)\tag{i}\\
\ell^X(\gg,\nu)&=\ell^{X'}(\gg,\nu)\tag {ii}\\
a_\xi^{X}(\gg,\nu)&=a_\xi^{X'}(\gg,\nu)\tag {iii}\\
\indiscbelong^X(\gb,a,\ga)&=\indiscbelong^{X'}(\gb,a,\ga)\tag{iv}
\end{align*}
\noindent
whenever the arguments are members of $X\cap X'$. The equality sign
here means that if either side is defined then both sides are
defined and they are equal.
\end{lemma}

\begin{proof}
If the lemma fails then one of the equations~(i--iv) must fail
cofinally often.  Suppose first that equation~(i) fails cofinally
often, say for $\vec
\gg=\seq{\gg_\gi:\gi<\gd}$ and $\vec\nu=\seq{\nu_\gi:\gi<\gd}$.  This
means that  $\sup\vec\nu=\gk$ and 
$\nu_\gi<\crit(\ce_{\gg_\gi})\le\gk$  and
$s^X(\gg,\nu_\gi)\not=s^{X'}(\gg,\nu_\gi)$
for each
$\gi<\gd$.  We may suppose
without loss of generality that $s^X(F_{\gi},\nu_\gi)$ exists for all
$\gi<\gd$.  If we set $a_\gi=s^X(F_{\gi},\nu_\gi)$ then $\vec a$ is a
sequence in $X$ such
that $\vec\nu\lessb\vec a$ and $\vec a$ is a \pris\ for $\vec F$ where
$F_{\gi}=\ce_{\gg_\gi}$. 
Since this is a first order assertion about $\vec a$ in
$X$, 
there must be some sequence
$\vec a'$ in $X'$ which satisfies the same assertion in $X'$,
and hence $s^{X'}(\gg_\gi,\nu_\gi)$ also exists for all  sufficiently
large $\gi<\gd$.
Let $a'_\gi=s^{X'}(\gg_\gi,\nu_\gi)$ for each $\gi<\gd$.  Then $\vec
a'\geb\vec a$, since otherwise it is true in $V$ that there is a
\pris\ for $\vec F$ which is smaller than $\vec a$ cofinally often.
Then the same statement is true in $X$, contradicting the choice of~
$\vec a$.  Similarly $\vec a\geb\vec a'$ so $\vec a\eqb\vec a'$, which
means that equation~(i) holds for all but boundedly many $\gi<\gd$,
contrary to the choice of $\vec\gg$ and $\vec \nu$.

\paranumber{33}
The proof of equation~(ii) and (using lemma~\ref{main-lem})
equation~(iv) is similar.  The proof of equation~(iii) is also
similar, but slightly more complicated because of the extra quantifiers
in the definition of the function $a^X$ and the possibility of
different subscripts $\xi_\gi$.
\end{proof}

\begin{lemma}\label{reachable-lem}
If $X$ is any \pcs\ then there is $\nu<\gk$ 
such that  for every ordinal $\gb\in X$ with
$\nu<\gb\le\SO(\gk)$, at least of the following holds:
\begin{enumerate} 
\item
$\gb\in h_{\gk}^X\image\gb$.
\item
$\gb=\indiscbelong^X(\gb',a,\ga)$ for some ordinals $\gb'$, $a$ and
$\ga$ such that 
 $a<\gb$ and
$\ga<\gb'<\SO(\ga)$ with $\ga$ and $\gb'$
in $h_{\gk}^X\image a$.
\item
$\gb=s^X(\gg,\nu)$ for some ordinals $\nu<\gb$ and $\gg\in h_{\gk}^X\image\gb$.
\item
$\gb=a^{X}_{\xi}(\gg,\nu)$ for some $\nu<\gb$ and $\gg\in
h_{\gk}^X\image\nu$ and some countable ordinal~$\xi$.
Furthermore, if $\gd=\gw$ then $\xi$ may be taken to be $1$.
\end{enumerate}

\end{lemma}

\paranumber{33a}
\begin{proof}
Set $h=h^{X}_{\gk}$.
If $\gb$ cannot be written in the first form then $\gb$ is an
indiscernible in $X$, and if it also cannot be written in the second form
then it must be a principal indiscernible in $X$ for some extender
$F_\gb\in h\image\gb$.

Now let 
$\eta\in h\image\gb$ be the $\tri$-largest ordinal such that $\gb$ is an
accumulation point for $\eta$.
Then for some $\nu_0<\gb$, 
$\lis^X(\eta,\nu_0)$ either does not exist or is greater than or
equal to $\gb$.  Define $\nu_{\gi}=a_{\gi}^X(\eta,\nu_0)$ for each
$\gi\le\gd$.  If $\nu_\gi=\gb$ for some $\gi<\gd$ then $\gb$ falls into case~(4).  Otherwise $\nu_{\gi}$ exists and
$\nu_{\gi}\le\gb$ for all $\gi\le\gd$.  In this case
$\nu_{\gd}=\lis^Y(\eta,\nu_0)$ by lemma~\ref{semi-cont}, so $\gb=\nu_\gd$
and so $\gb$ 
falls into case~(3).  

If $\gd=\gw$ and $\gb$ falls into case~(4) then $\xi$ is a
successor since $\xi<\gw$, so $\gb=a^X(\eta,\nu_{\xi-1})$ as required by the second sentence of
clause~(4). 
\end{proof}

\paranumber{34}
\begin{cor}\label{next-acc}
Let $\ce_\eta$ be $\tri$-largest such that $\gg$ is an accumulation point for
$\eta\in\PO(\ga)$ in $X$, and suppose  
$\cof(\gg)>\inf(\gw_1,\gd)$.
 Then there is $\nu<\gg$ such that $\gg=a^X(\eta,\nu)$.
\end{cor}
\begin{proof}
Define $\vec\nu$ as in the last proof, and
let $\eta=\inf(\gw_1,\gd)$.  
If $\nu_{\gi+1}=\gg$ for some $\gi<\eta$ then
$\gg=a^X(\eta,\nu_{\gi})$.  Otherwise $\nu_{\eta}<\gg$ since
$\cof(\gg)>\eta$, but this is
impossible
because lemma~\ref{semi-cont} implies that
$\nu_{\eta}\ge\lis^X(\eta,\nu_0)$. 
\end{proof}

%

%

%
\section{Applications}
The main result to be proved in this section is the following theorem:

\paranumber{1}
\begin{theorem}\label{app-main-thm}
Suppose that $\gk$ is a strong limit cardinal with 
$\cof(\gk)=\gd<\gk$, and that $2^\gk\ge\gl>\gk^+$, where if $\gl$ is a
successor cardinal then the predecessor of $\gl$ has cofinality greater
than $\gk$.
\begin{enumerate}
\item If $\gd>\gw_1$ then $o(\gk)\ge\gl+\gd$.
\item If $\gd=\gw_1$ then  $o(\gk)\ge\gl$.
\item If $\gd=\gw$ then either $o(\gk)\ge\gl$ 
or else $\set{\ga:K\sat o(\ga)\ge\ga^{+n}}$ is cofinal
in $\gk$ for each $n<\gw$.
\end{enumerate}
\end{theorem}

The proof of theorem~\ref{app-main-thm} is like that in Gitik
\cite{gitik.negate-sch}.
It has two ingredients:
the first is the analysis of indiscernibles which was given in
section~2, and   the second is a result of
Shelah which is given below, following
some preliminary definitions, as theorem~\ref{shelah.thm}.

\paranumber{2}
As in the last section, if $\vec c$ and $\vec c'$ are in $\prod\vec b$
then we will write $\vec c\lessb\vec c'$ to mean that $\set{b:c_b\ge
c'_b}$ is bounded in $\sup(\vec b)$, and $\vec c\eqb\vec c'$ to mean
that $\set{b:c_b\not=c'_b}$ is bounded in $\sup(\vec b)$.  If $\vec b$
is a sequence of cardinals then a subset $\cd$ of $\prod\vec b$ is
said to be {\it cofinal} in $\prod\vec b$ if for each sequence $\vec
c\in\prod\vec b$ there is a sequence $\vec d\in\cd$ such that $\vec
c\lessb\vec d$.  The set $\prod\vec b$ is said to have {\it true
cofinality} $\gl$, written $\tcf\left(\prod\vec b\right)=\gl$, if
there is a sequence $\seq{\vec c_\nu:\nu<\gl}$ of members of
$\prod\vec b$ which is cofinal in $\prod\vec b$ and linearly ordered
by $\lessb$.

\begin{theorem}\label{shelah.thm}
{\rm(}Shelah{\rm)}
Suppose that $\cof(\gk)=\gd<\gk$ and $2^\gk\ge\gl$, where $\gl$ is a
regular cardinal.  If $\gd=\gw$ then also assume that $\gl<\gk^{+\gw}$.
\begin{enumerate}
\item 
\cite[chap. IX, 5.12 and 5.10(1)]{shelah.card-arith}
There is a sequence $\vec a\subset\gk$ of regular cardinals such
that 
$\tcf(\prod\vec a)=\gl.$
\item 
\cite[chap. II, 1.2]{shelah.card-arith}
Any strictly increasing sequence from $\prod\vec a$ of
length less than $\gl$ and cofinality greater than 
$\gk$ has a least upper bound.
\end{enumerate}
\end{theorem}

\paranumber{3}
In subsection~3.1 we apply the techniques of section~2 to the
sequence given 
by Shelah's theorem.  For $\gd>\gw$ this analysis leads directly to
the proof of 
lemma~\ref{uncnt-seq} below, which is clauses~1
and~2 of 
theorem~\ref{app-main-thm} except that 
clause~1 is weakened by replacing $\gl+\gd$ with $\gl$.
For clause~3, the case $\gd=\gw$, the analysis yields
lemma~\ref{cnt-seq}, which is used in subsection~3.2 to prove
clause~3 of theorem~\ref{app-main-thm}.
In subsection~3.3 we prove various further results,
including the full strength of theorem~\ref{app-main-thm}(1).

\paranumber{4}
As usual, all successors are computed in $K$.

\begin{lemma}
\label{uncnt-seq}
If $\gk$ is a strong limit cardinal with $\gw<\gd=\cof(\gk)<\gk$, and 
$\gk^+<\gl\le 2^{\gk}$
where $\gl$ is not the successor of a cardinal of
cofinality less than $\gk$, then $o(\gk)\ge\gl$. 
\end{lemma}
Notice that if lemma~3.2 is true for successor cardinals $\gl$ then it is
true for all limit cardinals.  Thus it will be sufficient to prove
lemma~\ref{uncnt-seq} for regular cardinals $\gl$.

\begin{lemma}\label{cnt-seq}
Suppose that $\gw=\cof(\gk)<\gk$ and $\gk^+<o(\gk)<\gl\le 2^{\gk}$,
and assume that $\gl$ is regular and $\set{\ga:o(\ga)>\ga^{+n}}$ is
bounded in $\gk$ for some $n<\gw$.
Then there is a countable sequence $\vec b$, cofinal in $\gk$, along
with  continuous, nondecreasing functions
$f_b$, and ordinals $\gg_b$, $\ga_b$ and $\gs_b$ for $b\in\vec b$
such that $\gs_b<b$ and $\gs_b\in\vec b$
for all but 
boundedly many $b\in\vec b$ 
and the set  $\cl\in\prod\vec b$, defined below,  is cofinal in
$\prod\vec b$
and has true
cofinality $\gl$.

A sequence  $\vec c\in\prod\vec b$ is in $\cl$ if and only if
for some \pcs\ $Y$, and all
sufficiently large $b\in \vec b$,
\begin{equation}
\label{e.bdef}
c_{b}=
\begin{cases}
\lis^Y(f_b(c_{\gs_b}),\gg_{b})&\quad\text{if $b$ is a limit of principle
indiscernibles}\\ 
\indiscbelong^Y(f_b(c_{\gs_b}),\gg_b,\ga_b)&\quad\text{otherwise.}
\end{cases}
\end{equation}

\paranumber{5}
Furthermore any strictly increasing, non-cofinal subsequence of $\cl$
 of cofinality greater than $\gk^+$ has a least upper bound in
 $\prod\vec b$.
\end{lemma}

\subsection{Proof of lemmas~\ref{uncnt-seq} and~\ref{cnt-seq}}
The main goal of this subsection is to prove 
lemma~\ref{cnt-seq}. This is true for the case of uncountable
cofinality, $\gd>\gw$, as well as for countable cofinality---the
difference is that in the 
case $\gd>\gw$ we immediately reach a easy contradiction and
hence do not  need to explicitly state an intermediate result corresponding
to lemma~\ref{cnt-seq}.

\paranumber{6a}
We will write $S^{f,Y}_b(\nu)$ for the function given in
equation~\eqref{e.bdef} of lemma~\ref{cnt-seq}.  Thus a sequence
$\vec c=\seq{c_b:b\in\vec b}$ in $\prod\vec b$ is in $\cl$ if and 
only if there is a \pcs\ $Y$ so that 
$c_b=S^Y_b(c_{\gs_b})$  for all sufficiently large
$b\in\vec b$. 

\paranumber{6b}
If we had required $c_b=S^Y_b(c_{\gs_b})$ for \emph{all} $b\in\vec b$ such
that $\gs_b\in \vec b$, 
and if $S^Y_b$ did not depend on $Y$, then it would follow that
a sequence $\vec c\in\cl$ is determined by 
$\seq{c_b:\gs_b\notin\vec b}$.  Since $\set{b\in\vec
b:\gs_b\notin\vec b}$ is bounded in $\gk$ there are fewer than $\gk$
choices for $\set{c_b:\gs_b\notin\vec b}$ and it would then 
follow that
$\tcf(\cl)<\gk$, contradicting the assertion that $\tcf(\cl)=\gl$ and
completing the proof of the theorem.

\paranumber{6c}
For the case $\cof(\gk)>\gw$ this is nearly what happens. We  show
that $\vec b$ has order type $\gd$,  and then Fodor's theorem implies that
$\gs_b$ is constant on an unbounded subset $y$ of $\gk$.  Since
$S^{f,Y}_b=S^{f,Y'}_b$ for sufficiently large $b\in\vec b$,
this implies that $\set{\vec c\restrict y:\vec c\in\cl}$ had  fewer
than $\gk$ members,
 modulo the relation $\eqb$, and this  contradicts the
assertion that $\tcf(\cl)=\gl$.

\paranumber{6d}
The case $\cof(\gk)=\gw$  is more difficult. The strategy is to try to
show that $\cof(c_b)=\cof(c_{\gs_b})$ for $\vec c\in\cl$ and
sufficiently large $b\in\vec b$, which would lead to essentially 
the same contradiction as in the case $\gd>\gw$.
In fact, however, it becomes necessary to
look at sequences $\vec d$ which are the least upper bound  for
certain subsets of $\cl$, 
instead of working with the sequence $\vec c$ in $\cl$ directly.
This argument is in subsection~3.2.

\smallskip
\paranumber{6e}
To understand the proof of  lemma~\ref{cnt-seq}, it will be helpful to
consider four levels of data:
\begin{description}
\item[Level 1.]
The functions $\lis^Y$ and $\indiscbelong^Y$.
\item[Level 2.]
The sequence $\vec b\supset\vec a$ and the parameters $\gs_b$,
$\gg_b$ and 
$\ga_b$.  Also included in this level is a procedure for defining the
functions $f_b$ from a single function $f\in K$---this is the function
$f$ which appears as a subscript in the notation $S^{f,Y}_b$.  The
procedure involves additional parameters $p_b$ and $\eta_b$.
\item[Level 3.]
The function $f$ used to define the functions $f_b$.
\item[Level 4.]
The set $\cl$, and the sequences $\vec c\in\cl$.
\end{description}

\paranumber{6f}
The items in level~4 is already defined in the lemma, using data from 
from   levels~1--3.

\paranumber{6g}
The functions $\lis^Y$ and
$\indiscbelong^Y$ of level~1 were defined, and their properties
proved, in the last section.  In particular we use 
lemma~\ref{next-ind-fctn}, which asserts that these functions are essentially
independent of $Y$, and lemma~\ref{reachable-lem} which provides
the inspiration for 
lemma~\ref{cnt-seq}.  Note, for example that the first case,
\begin{equation}
c_{b}=\lis^Y(f_b(c_{\gs_b}),\gg_{b})
\tag{$*$}\label{e.cb-def}
\end{equation}
of lemma~\ref{cnt-seq} comes from   case~(3) of lemma~\ref{reachable-lem}:
\begin{equation}
\nu=s^{Y}(h^{Y}_{\gk}(\nu'),\gg)\qquad\text{for some $\nu',\gg<\nu$}.
\tag{$**$}\label{e.nu=s}
\end{equation}
The equation~\eqref{e.cb-def} has 
 $\lis^Y$, which is more convenient to work with, instead of $s^Y$.
The parameter  
$\gg=\gg_b$ is made to  depend only on $b$.
Equation~\eqref{e.cb-def} asserts that 
whenever $c_b$ is a member of $\vec c$ then 
the ordinal $\nu'$ of equation~\eqref{e.nu=s} is also a 
member $c_{\gs_b}$ of $\vec c$ (unless $\gs_b\notin \vec b$).
Furthermore
the coordinate $\gs_b$ at which
$c_{\gs_b}$ appears in $\vec c$ has been fixed and does not depend on 
the sequence $\vec c$.  
Finally, the function $h^Y_\gk$, which depends on $Y$ and which
need not be in $K$, is replaced with a function $f_b$ in $K$
which again does not depend on $\vec c$ or $Y$.

\paranumber{6h}
 The data in level~2 is defined by working in a fixed  \pcs\
$X$, with the aim of finding parameters so that
 the set of restrictions $\vec c\restrict
\vec a=\seq{c_b:b\in\vec a}$ of sequences $\vec
c\in\union\set{\cl^{f,X}:f\in X\cap K}$ is cofinal in $\prod\vec a\cap X$.
For most of this construction we let $h^X_\gk$ play the role of $f$,
but  at the end we use the covering lemma to 
to show that there are suitable approximations to $h^X_\gk$ 
in $X\cap K$.

\paranumber{6i}
 The argument for level~3 begins with the  observation that  by
elementarity the set of restrictions $\vec c\restrict \vec a$ of
sequences 
\[\vec c\in\union\set{\cl^{f,Y}:f\in K\tand Y\text{ is a \pcs}}\] 
is cofinal in $\prod\vec a$. In order to prove the crucial fact that
there is a 
single function 
$f$ such that $\cl^f=\union\set{\cl^{f,Y}:Y\text{ is a \pcs}}$ is
similarly cofinal we use the assumption that
$(o(\gk)^{\gk})^K<\gl=\cof(\gl)$, and this is the only place where
this assumption is used.
In section~3.3 we prove slightly stronger versions of
theorem~\ref{app-main-thm} by
using a modification of this assumption which also implies the
existence of a single function $f$ so that $\cl^f$ is cofinal.  Thus
this modified assumption
leads to the same contradiction.
\smallskip

\medskip
\paranumber{7}
We are now ready to begin the proof of lemmas~\ref{uncnt-seq}
and~\ref{cnt-seq}. 
First we need a couple of preliminary results.  This first lemma will be
applied to sequences $\vec b$ which may not be increasing.
\begin{prop}
\label{cof-unbnded-prop}
If $\vec b$ is a sequence of cardinality at most $\vec\gd$ and
$\eta<\gk\le\tcf(\prod\vec b)$ then $\set{b\in \vec
b:\cof(b)\le\eta}$ is bounded in $\vec b$.
\end{prop}
\begin{proof}
Suppose to the contrary that $\eta<\gk$ but  $\vec b'=\seq{b\in\vec
b:\cof(b)<\eta}$  is cofinal in $b$.
Then $\tcf(\prod\vec b')=\tcf(\prod\vec b)\ge\gk$.  Now let
$y_b$ be a cofinal subset of $b$ of cardinality
at most $\eta$ for each $b\in\vec b'$.   Then
$\prod_{b\in\vec b'}y_b$ is cofinal in $\prod\vec b'$, but this is
impossible since
$\gk$ is a strong limit
cardinal and hence $\card{\prod_{b\in\vec b'}y_b}\le\eta^{\gd}<\gk$.
This contradiction proves the proposition. 
\end{proof}

Now let $X$ be a \pcs.  This \pcs\ will remain fixed through the rest
of this subsection.
\begin{defin}
\label{welladj-def}
We say that an ordinal $b\in X$ is {\it well adjusted} in $X$ 
\notate{well adjusted}
if 
\begin{enumerate}
\item $b>\rho^X$, and if $\gd=\gw$ then there is an $n<\gw$ such that
$o(\ga)<\ga^{+n}$ whenever $b<\ga<\gk$.
\item
$b$ is regular in $\core$.
\item
$b\cap X$ is not cofinal in $b$, 
\item
The indiscernibles (including nonprincipal indiscernibles) of $X$ are
cofinal in $X\cap b$.
\end{enumerate}
\end{defin}

\paranumber{8}
\begin{prop} \label{well-adj-prop}
 If $\vec b$ is a sequence of regular cardinals of $\core$ such
that $\tcf(\prod\vec b)>\gk^+$ then every sufficiently large member of
$\vec b$ is well adjusted in $X$.
\end{prop}
\begin{proof}
First, $\vec b$ is unbounded in $\gk$ by
proposition~\ref{cof-unbnded-prop},
so clause~(1) of definition~\ref{welladj-def}
 is satisfied for all sufficiently large $b<\gk$.
Clause~(2) is satisfied by hypothesis, and 
proposition~\ref{cof-unbnded-prop} implies that 
$\cof(b)>\card{X}$ for all but boundedly many members of $b$,  so that 
$X\cap b$ is not cofinal in $b$ for $\card X<b<\gk$.  Thus we only
need to verify clause~(4).

\paranumber{9}
Let $\vec b'$ be the set of ordinals $b$ in
$\vec b$ such that the 
indiscernibles of $X$ are not cofinal in $X\cap b$, and suppose that,
contrary to clause~(3), $\vec b'$  is cofinal in $\vec b$.
For each member $b$ of $\vec b'$
pick an ordinal $\xi_b<b$ in $X$ which is larger than all of the
indiscernibles of $X$ below $b$, so that $h_\gk^X\image\xi_b$ is
cofinal in $b\cap X$.  If $\vec\nu$ is any member of $\prod\vec b'$ in
$X$ then there is a sequence $\vec\nu'\in\prod_{b\in
\vec b'}\xi_b$ such that $h^X\circ\vec\nu'\geb\vec\nu$.  Now $h_\gk^X$
need not be in $K$, but this construction only uses the restriction
$h_\gk^X\cap(\gk\times\gk)$ of 
$h_\gk^X$ to ordinals below $\gk$, which is a member of $K$. 
Thus it is true in
$V$, and hence by elementarity it is true in $X$, that there is a 
 function $f\in{^{\gk}\gk}\cap K$ and a sequence $\vec\nu'\in
\prod_{b\in\vec b'}\xi_b$ such that
$f\circ\vec\nu'\geb\vec\nu$.  Since $\vec \nu$ was arbitrary 
 it is true in $X$, and by
elementarity again it is true in $V$, that the set of sequences of the
form 
$f\circ\vec\nu'$ for some $f\in K$ and some $\vec\nu'\in\prod_{b\in\vec
b'}\xi_b$ is cofinal in $\prod\vec b'$.  Since $\tcf(\prod\vec
b')=\tcf(\prod\vec b)>\gk^+$ and $\card{{^{\gk}\gk}\cap K}=\gk^+$
 there must be a single function $f$ such that the
set of sequences $f\circ\vec\nu'$ for $\nu'\in\prod_{b\in\vec b'}\xi_b$
is cofinal in $\prod\vec b'$, but  this is impossible because the members
$b$ of $\vec b$ are regular in $K$ and hence $f\image\xi_b$ is bounded
in $b$ for all $b\in\vec b'$.
\end{proof}

\paranumber{10}
It follows that every sufficiently large member
of $\vec a$ is well
adjusted, and we can assume without loss of generality that every
member of $\vec a$ is well adjusted.

We are now ready the define the sequence $\vec b$ and
the associated parameters.  For each well adjusted ordinal $b\in
X\cap\gk$ we will define $\gs_b<b$ along with a function $S^{f,Y}_b$,
depending on an arbitrary precovering set $Y$ and  function $f\in K$
as well as the ordinal $b$. 
The function $S^{f,Y}_b$ also depends on several parameters which will
be 
fixed in the 
course of this definition.  The function $S^{f,Y}_b$ is the function
appearing in equation~\eqref{e.bdef} of lemma~\ref{cnt-seq}.  We will
show that if we take 
$\cl^{f,Y}$ to be the set of 
sequences $\vec c\in\prod \vec b$ such that
$c_{b}=S^{f,Y}_b(c_{\gs_b})$ for all sufficiently large $b\in\vec b$,
then the union over functions $f\in K$ and \pcs{s} $Y$ is cofinal
in $\prod\vec b$ and hence has true cofinality $\gl$.

\begin{defin}\label{def-stuff}
We define an ordinal $\gs_b$ for all well adjusted $b\in X\cap\gk$,
and a function $S_b^{f,Y}$ for all well
adjusted $b\in X\cap\gk$, all \pcs{s} $Y$, and $f\in K$.  We also
define several auxiliary parameters.  The definition depends on the
fixed \pcs\ $X$ and is broken into two cases,
depending on whether or not the principle indiscernibles of $X$ are cofinal
in $b\cap X$.

\paranumber{11}
\startcases
\begin{case}
(The principle indiscernibles of $X$ are cofinal in $X\cap b$.)
In this case we define
\begin{enumerate}
\item
$\ga_b=i^X_{b,\gk}(b)$.  Thus either $b=\ga_b$ or $b$ is a 
\pri\ for $\ga_b$.  In either case
$b$ is a  limit in $X$ of principle
indiscernibles for $\ga_b$.
\item
Since $\cof(b)>\gd$, corollary~\ref{next-acc} implies that there is
$\eta\le\SPO(\ga_b)$ and $\gg<b$ so that $b=a^X(\eta,\gg)$.  We let
$\eta_b$ be this ordinal $\eta$ and let $\gg_b$ be the least 
ordinal $\gg$ such that $b=a^X(\eta_b,\gg)$ and $\lis^X(\eta_b,\gg)\not< b$.
\item
$\gs_b$ is the least ordinal $\gs$ in $X$ such that
$\set{h_\gk^X(\nu,p):\nu\in\gs\cap X}$ is cofinal in $X\cap\eta_{b}$ for
some finite sequence $p$ of ordinals.
\item
$p_b$ is the least finite sequence $p$ of ordinals in
$X$ such that  
$\set{h_\gk^X(\nu,p):\nu\in\gs\cap X}$ is cofinal in $X\cap\eta_b$.
\item
If $f$ is any function in $K$ then $f_b$ is the function defined by
$f_b(\nu)=\sup(\eta_b\cap f\image(\nu\times\{p_b\}))$.
\item
$S^{f,Y}_b(\nu)=
\lis^Y(f_b(\nu),\gg_b)$, if it is defined and less than $b$.
Otherwise
$S^{f,Y}_b(\nu)$ is undefined.
\end{enumerate}
\end{case}

\paranumber{12}
\begin{case}
($b$ is not a limit of principle indiscernibles)

Since $b$ is a limit of indiscernibles, but not a
limit of principal indiscernibles, there is a largest principle
indiscernible below $b$.  Let $\gg_b$ be this principal indiscernible,
and set $\ga_b=i_{\gg_b,\gk}(\gg_b)$.  Then
$\gg_b$ is a principal indiscernible for $\ga_b$ and every
ordinal in $X\cap(\gg_b,b]$ is an indiscernible belonging to
$(\gg_b,\ga_b)$.

Now let $\eta_b=i_{\gg_b,\ga_b}(b)$, so that
$b$ is an indiscernible in $X$ for $\eta_b$ belonging to
$(\gg_b,\ga_b$). The ordinals  $\gs_b$ and $p_b$, and the function
$f_b$,
are defined exactly as in
case~1. 

\paranumber{13}
Finally, set 
$S^{f,Y}_b(\nu)=\indiscbelong^Y(f_b(\nu),\gg_b,\ga_b)$ if it exists
and is less than $b$. Otherwise
$S^{f,Y}_b(\nu)$ is undefined.
\end{case}
\end{defin}

\begin{prop}
If $b$ is well adjusted then $\gs_b<b$.
\end{prop}
\begin{proof}
If $b$ is not a limit of principle indiscernibles then
$\gs_b\le\gg_b<b$, so  suppose that $b$ is a limit of principle
indiscernibles and that, contrary to the proposition, 
$\gs_b=b$. Define a sequence $(c_{\gi}:\gi<\gd)$ by recursion on $\gi$:
\begin{align*}
c_0&=\gg_b\\
c_{\gi+1}&=\ell^X(\xi_\gi,\gg_b)&&
  \text{where }\xi_\gi=\sup\left(\eta_b\cap h_\gk^X\image(X\cap c_i)\right)\\
c_{\gi}&=\sup \set{c_{\gi'}:\gi'<\gi}&&
  \text{if $\gi$ is a limit ordinal}.
\end{align*}

\paranumber{14}
If $\gi$ is a limit ordinal then $c_\gi$ is in $X$, since $X$ is
$\gd$-closed, and $c_\gi<b$ since $X\cap b$ is not cofinal in $b$.  
If $\gs_b=b$ then it follows that $c_\gi<b$ for each
$\gi\le\gd$.  Set $\zeta=\inf(\gd,\gw_1)$.  
Then
$c_{\zeta}=a^X(\xi_{\zeta},\gg_b)$ and by proposition~\ref{semi-cont} it
follows that $c_{\zeta}=\ell^{X}(\xi_{\zeta},\gg_b)=s^X(\xi,\gg_b)$ for some
$\xi\ge\xi_{\zeta}$.   But $\xi\in h_\gk^X\image(c_{\zeta})$, so $\xi\in h_\gk^X\image(c_{\gi})$ for some $\gi<{\zeta}$ and hence
$c_{\zeta}>c_{\gi+1}>c_{\zeta}$.  This 
contradiction completes the proof that $\gs_b<b$.
\end{proof}

\paranumber{15}
\begin{defin}\begin{enumerate}
\item
 $\vec b$ is the smallest set
such that  $\vec a\subset\vec b$ and 
$\gs_b\in \vec b$ for all $b\in\vec b$ such that $\gs_b$ is well adjusted
in $X$.  
\item
If $f$ is as above and $Y$ is a \pcs\ containing everything relevant
then $\cl^{f,Y}$ is the set 
of sequences $\vec c\in\prod \vec b$ such that
$c_b=S^{f,Y}_b(c_{\gs_b})$ for all sufficiently large $b\in\vec b$
such that $\gs_b\in\vec b$.
\item
If $f\in K$ then 
$\cl^f=\union\set{\cl^{Y, f}:Y\text{ is a \pcs}}$.
\end{enumerate}
\end{defin}
Notice that  by lemma~\ref{next-ind-fctn}
$\cl^{f}$ is first order definable in any $Y$ containing all
of the data, and that $\cl^{Y, f}=\cl^{f}\cap Y$.

\begin{lemma}\label{Ufclf-cofinal}
The set $\union_{f}\set{\vec c\restrict a:\vec c\in\cl^f}$ is
cofinal in $\prod\vec a$. 
\end{lemma}
\paranumber{16}
\begin{proof}
Since $\cl^{f}$ is first order definable, it is sufficient to
show that 
the lemma is true in $X$; that is,
to produce, given any sequence $\vec d$ in $\prod\vec a\cap X$, 
a function $f\in\core\cap X$ and a 
sequence 
$\vec c\in\cl^f\cap X$ such that $\vec c\restrict\vec a\geb\vec d$. 
For the function $f$ we will use $h_\gk^X$, or rather a
function in $X\cap K$
which is nearly equal to $h_\gk^X$.
We begin by  defining a  sequence $\vec c_n=\seq{c_{n,b}:b\in\vec b}$
for each $n<\gw$ so that
\begin{align*}
c_{0,b}&=d_b&&\text{if $b\in\vec a$}\\
c_{0,b}&=0&&\text{if $b\notin\vec a$}\\
S_b^{X,h_\gk^X}(c_{n+1,\gs_b})&\ge c_{n,b}
  &&\text{if $\gs_b\in\vec b$}\\
c_{n+1,b}&\ge c_{n,b}&&\text{for all $n$ and $b$.}
\end{align*}
We define $c_{n,b}$ by recursion.  Suppose that $c_{n,b}$ has been
defined for all $b\in\vec b$, and $c_{n+1,\gs_b}$ has been defined if
$\gs_b$ is in $\vec b$.
In order to define $c_{n+1,b}$, define 
$\xi_{b'}$ for each $b'\in\vec b$ to 
be the least ordinal $\xi$ such that $S_{b'}^{X,h_\gk^X}(\xi)\ge
c_{n,b'}$ if  
$b=\gs_{b'}$, and let $\xi_{b'}=0$ otherwise.  Then
$\set{\xi_{b'}:b'\in\vec b}\in 
X$ since $^{\gd}X\subset X$ and we can set
$c_{n+1,b}=\sup(\set{\xi_{b'}:b'\in\vec b}\cup\{c_{n,b}\}$). 

\paranumber{17a}
Set $y=\set{(c_{n,b},p_b):n\in\gw\text{ and }b\in\vec b}$.
By lemma~\ref{ftoi-lem}(5) there is a function $f\in X\cap K$ 
 such that
$f\restrict y=h_\gk^X\restrict y$.  Define 
the sequence 
$\vec c\in\cl^f$ by recursion on $b\in\vec b$: 
\[c_b=
\begin{cases}
\union_nc_{n,b}&\text{if $\gs_b\notin \vec b$}\\
S_b^{X,f}(c_{\gs_b})&\text{if $\gs_b\in \vec b$}.
\end{cases}
\]
We claim that  $ c_{n,b}\le c_{b}<b$ for each
$n\in\gw$.  The proof is a simple  recursion on $b\in\vec b$.  It is
true immediately if $\gs_b\notin\vec b$, while if $\gs_b\in\vec b$
then $c_{\gs_b}\ge c_{n+1,\gs_b}$ so 
$c_b=S^{X,f}_b(c_{\gs_b})\ge S^{X,f}_b(c_{n+1,\gs_b})\ge c_{n,b}$.

\paranumber{17b}
In particular, $c_b\ge c_{0,b}=d_b$ for $b\in a$, so $\vec
c\restrict a \ge\vec d$ as required.
\end{proof}

\paranumber{17c}
This completes the construction at level~2 as described in the
introduction to this subsection.  The
following corollary gives us level three, the choice of the function
$f$, and is thus much more important than its length suggests.  Notice
that this is the only place where we use the assumption 
that $(o(\gk)^{\gk})^{\core}<\cof(\gl)$.
\begin{cor}\label{tcfb=gl}
There is a function $f\in K$ such that $\set{\vec c\restrict\vec
a:\vec c\in\cl^f}$ is cofinal in $\prod\vec a$.
\end{cor}
\begin{proof}
The last lemma implies that $\union\set{\vec c\restrict \vec a:\vec
c\in\union_f\cl^f}$ is cofinal in $\prod\vec a$.  Now the relevant
functions $f\in K$ have domain contained in $\gk\times\gk^{<\gw}$ and
range contained in $\PO(\gk)$, so there are only
$(o(\gk)^{\gk})^K<\gl$ of them.
Since $\tcf(\prod\vec a)=\gl=\cof(\gl)$
 it follows that there is a single
function $f$ such that 
$\set{\vec c\restrict \vec a:\vec c\in\cl^f}$ is cofinal in
$\prod\vec a$. 
\end{proof}

\paranumber{18}
\begin{cor}\label{gs-unbnded}
The set $\set{b\in\vec b:\gs_b\notin \vec b}$ is bounded in $\vec b$, and
if $\nu<\gk$ then $\set{b\in\vec b:\gs_b<\nu}$ is bounded in $\vec b$.
\end{cor}
\begin{proof}
Recall that the functions $S_b^{f,Y}\colon\gs_b\to b$ are nondecreasing
and are cofinal in $b\cap Y$, whether or not $\gs_b\in\vec b$.
If we set  $S^{f,Y}(\vec
d)=\seq{S^Y_b(d_b):b\in\vec b}$ then it follows that
\[
\set{s^{f,Y}(\vec d):\vec d\in\prod_{b\in\vec b}\gs_b\tand Y \text{
is a \pcs}}
\]
is cofinal in $\prod\vec b$, 
and since $\vec d\lessb\vec d'$ implies $S^{f,Y}(\vec d)\leb
S^{f,Y'}(\vec d')$ for 
any \pcs{s} $Y$ and 
$Y'$ it follows that 
$\tcf(\prod_b\gs_b)=\tcf(\prod\vec b)=\gl$, and the
corollary follows by propositions~\ref{cof-unbnded-prop}
and~\ref{well-adj-prop}. 
\end{proof}

\medskip
\paranumber{19}
At this point we will treat the cases $\gd=\gw$ and $\gd>\gw$
separately. We 
begin with $\gd>\gw$,  finishing
the 
proof of lemma~\ref{uncnt-seq} by assuming that $\gd>\gw$ and
 showing that the properties
which we have established for the sequence $\vec b$ lead to a
contradiction.

\begin{proof}
[Proof of corollary~\ref{uncnt-seq}]
First we show that $\ot(\vec b)=\gd$.  
Set $\vec a_0=\vec a$ and for
$n\ge 0$ set $\vec a_{n+1}=\vec a_{n}\cup\set{\gs_{b}:b\in\vec
a_n\text{ and }\gs_b\in\vec b}$.  Since $\ot(\vec a)=\gd$, 
corollary~\ref{gs-unbnded} implies
that each $\vec a_n$ has order type $\gd$.  But $\vec b=\union_n\vec
a_n$, and since $\cof(\gd)>\gw$ it follows that $\vec b$ has order
type $\gd$.

Now since $\gs_b<b$ and $\gs_b\in\vec b$ for every sufficiently large
$b\in\vec b$, Fodor's theorem implies that there is an unbounded subset
$B$ of $\vec b$ such that $\gs_b$ is constant for $b\in B$.  But this
contradicts corollary~\ref{gs-unbnded}, and this contradiction shows
that it is not possible that $o(\gk)<\gl$.
\end{proof}

\medskip
\paranumber{20}
We now finish this subsection by completing the proof of
lemma~\ref{cnt-seq}.  We assume that $\gd=\gw$, and that the
hypothesis of lemma~\ref{cnt-seq} holds.  

\begin{proof}
[End of proof of lemma~\ref{cnt-seq}]
We have proved all of this lemma except for the last paragraph, which
asserts that that every non-cofinal subsequence $\cb$ of  $\cl^f$ of cofinality
at least $\gk^+$ has a least upper bound in $\prod\vec b$. 
Given such a
subset $\cb$, let $\vec d$ be the least upper bound of $\set{\vec c\restrict
\vec a:\vec c\in\cb}$, which exists by
clause~(2) of theorem~\ref{shelah.thm}.  

\paranumber{21a}
Define $b\gsless b'$ if for some $m> 0$ 
there is a chain $b=b_0<b_1<\dots<b_m=b'$
such that $b_{k}=\gs_{b_{k+1}}$ for $k<m$.
Let $Y$ be a \pcs\ with $\vec d$ and $\cb$ in $Y$, and 
for $b\gsless b'$ define
\[
S^{f,Y}_{b,b'}=
S^{f,Y}_{b_m}\circ S^{f,Y}_{b_{m-1}}\circ\dots\circ
S^{f,Y}_{b_1}\colon b\to b'
\]
We will extend this to $b\gslesseq b'$ by setting
$S^{f,Y}_{b,b}(\nu)=\nu$. 

\paranumber{21b}
Define $\vec d'\in\prod b$ by
letting $d'_{b}$ be the least ordinal $\nu$ such that
$\nu\ge d_b$ if $b\in\vec a$ and
$S^{f,Y}_{b,a}(\nu)\ge d_{a}$ for all $a$ such that $b\gsless a\in\vec
a$. This is possible since $S^{f,Y}_{b,a}$ is cofinal in $a\cap Y$,
which has cofinality greater than $\gd=\gw=\card {\vec a}$.
We claim that $\vec d'=\lub(\cb)$.

\paranumber{21c}
Any member of $\cb$ must be less than $\vec d'$ except on a bounded
set, so it will be sufficient to prove that $\vec d'$ is minimal.
We need to show that if   $\vec c$ is any sequence such that $\vec
c\lessb\vec d'$,  then 
$\vec c\lessb\vec d''$ for some  sequence $\vec d''\in\cb$.

\paranumber{21d}
To find $\vec d''$, define $\vec c'\in\prod\vec a$ by setting
$c'_a=\sup\set{S^{f,Y}_{b,a}(c_b):b\gslesseq a}$.
Each of the ordinals $S^{f,Y}_{b,a}(c_b)$, for $b\gslesseq a$, is
smaller than $d_a$ by the choice of $\vec d'$.  But $\cof(d_a)>\gw$
for all but boundedly many $a\in\vec a$ by
proposition~\ref{cof-unbnded-prop}, since $\tcf(\prod\vec d)\ge\gk^+$,
and hence $\vec c'\lessb\vec d$.

Since $\vec d$ is the least upper bound of $\set{\vec c\restrict\vec
a:\vec c\in\cb}$ it follows that there is $\vec d''\in\cb$
such that $\vec c'\lessb\vec d''\restrict\vec a\lessb\vec d'$.  Since
$S^{f,Y}_{b}$ is increasing, it follows that $\vec c\lessb\vec d''$,
as required.
\end{proof}
\paranumber{20}
\subsection{Countable cofinality: the proof of theorem \ref{app-main-thm}(3)}
Except for the need to consider nonprincipal extenders the proof of
theorem~\ref{app-main-thm}(3) is essentially the same as in
\cite{gitik.strength-not-sch}.  
We assume
that theorem~\ref{app-main-thm} is false with $\gd=\cof(\gk)=\gw$, and
let $\vec b$, $\cl$, and the associated ordinals be as given by
lemma~\ref{cnt-seq}.  We will assume that $\vec b$ and $\cl$ are
members of every \pcs\ mentioned.

\paranumber{22a}
\begin{defin}
Let $\cd$ be the class of sequences $\vec d\in\prod\vec b$ such that
$\vec d$ is the least upper bound of an increasing subsequence of
$\cl$ of order type $\gk^+$.
\end{defin}

\paranumber{22b}
Note that $\ot(\cd,\lessb)=\gl$ by lemma~\ref{cnt-seq}.
For sequences  $\vec d$ and $\vec d'$ in $\cd$  let $g(\vec d,\vec d')$ be
the least ordinal $b_0\in\vec b$, if there is one, such that
\begin{enumerate}
\item
either
$d_{b}<d'_b$ 
for all $b>b_0$, or $d_b>d'_b$ for all $b>b_0$, and
\item
If there are only boundedly many $b\in\vec b$ such that
\begin{equation}
\cof^K(d_b)=\cof^K(d'_{b})>\gg_b\tand
\cof^{K}(d_{\gs_b})\not=\cof^{K}(d'_{\gs_b})\tag{$*$}
\end{equation}
then $(*)$ is false for all $b>b_0$.
\end{enumerate} 
 Since
$\tcf(\cd)=\gl$ there is a subset $\cd'$ of $\cd$ of cardinality $\gl$ which
is linearly ordered under $\lessb$, so that $g(\vec d,\vec
d')$ is defined for all $\vec d,\vec d'\in\cd'$.  Since
$\gl>(2^{\gw})^+$ and $\card{\vec b}=\gw$, 
the Erd\H os-Rado theorem implies that there is a
sequence $D=\seq{\vec d_\gi:\gi<\gw_1}$ such that $g$ is constant on
$[D]^{2}$.  We can assume \iwlog\ that $g(\vec d_\gi,\vec d_{\gi'})=0$
for all $\gi<\gi'<\gw_1$, so that $d_{\gi,b}<d_{\gi',b}$ whenever
$\gi<\gi'$ and $b\in\vec b$.  

\paranumber{22c}
\newcommand\opA{\operatorname{\textbf Q}}
We will say that some property  $\opA(\gi,b)$
holds \emph{for almost all $(\gi,b)$}
if for all but countably many $\gi<\gw_1$
there is $\nu_{\gi}<\gk$ such that 
$\opA(\gi,b)$ holds for all $b\in\vec b\setminus \nu_\gi$.

\paranumber{22d}
\begin{lemma}
\label{t.dcases}
For almost all $(\gi,b)$ the relations in the following table hold.
Here $I$ is the  set of
$b\in\vec b$ such that $b$ is a limit of principle indiscernibles.  In
case~2c, ``almost all'' means that there is $\gi_0<\gw_1$ such that
for all $\gi,\gi'>\gi_0$ there is $\nu_{\gi,\gi'}<\gk$ such that the
conclusion holds whenever the hypothesis is true and
$b>\nu_{\gi,\gi'}$.
\begin{center}
\newcommand\xstrut{\vphantom{${\Bigl)}$}}
\begin{tabular}
{|ll|l|l|}
\hline
\multicolumn{3}{|c|}{Hypothesis}&\multicolumn{1}{c|}{Conclusion}\\
\hline
\raisebox{-1.8ex}[0pt][0pt]{(1)}&
\multicolumn{2}{l|}
{\raisebox{-1.8ex}[0pt][0pt]{$b\in I$}}&
    \xstrut$\cof(d_{\gi,\gs_b})=\cof(d_{\gi,b})$\\
  &\multicolumn{2}{l|}{} &\xstrut $d_{\gi,b}$  is regular in $K$\\
\hline
\xstrut(2a)&&
 $\cof^K(d_{\gi,b})<\gg_b$&$\cof^K(d_{\gi,\gs_b})=\cof^K(d_{\gi,b})$\\
\cline{3-4}
\xstrut(2b)&$b\notin I$&$\cof^K(d_{\gi,b})=\gg_b$& impossible\\
\cline{3-4}
\xstrut(2c)&& $\cof^K(d_{\gi,b})=\cof^K(d_{\gi',b})>\gg_b$&
$\cof^K(d_{\gi,\gs_b})=\cof^K(d_{\gi',\gs_b})$\\
\hline
\end{tabular}
\end{center}
\end{lemma}

\paranumber{22e}
Before proving lemma~\ref{t.dcases} we will show that it implies  the
theorem.   As before, we write $b\gsless b'$ if there is a 
chain $b=b_0<\dots<b_m=b'$ with $b_{i}=\gs_{b_{i+1}}$ for $i<m$, and
we write $S^Y_{b,b'}$ for the composition $S^Y_{b_1}\circ\cdot\circ
S^Y_{b_m}$.  We write $b\gslesseq b'$ if $b\gsless b'$ or $b=b'$, and we
set $S^Y_{b,b}(\nu)=\nu$.

\paranumber{22f}
\begin{proof}
[Proof of theorem~\ref{app-main-thm}(3), assuming lemma~\ref{t.dcases}]
By throwing out countably many sequences from $\seq{\vec
d_{\gi}:\gi<\gw_1}$ we can assume  without loss of generality that
lemma~\ref{t.dcases} is valid for all $\gi<\gw_1$, 
for all sufficiently large $b\in\vec b$.
By the
definition of $g(\vec d,\vec d')$ it follows that case~(2c) is valid for
all $\gi,\gi'\in\gw_1$ and for all  $b\in\vec b$, and by dropping to an
uncountable subset of 
$\seq{\vec d_\gi:\gi<\gw_1}$ and throwing out a bounded part of $\vec
b$ we can assume that the other cases are
also valid for all $\gi<\gw_1$ and all $b\in\vec b$.
			
\begin{claim}
For each $\gi<\gw_1$, the set of $b\in \vec b$ such that $d_{\gi,b}$
falls into case~(2c), 
that is, such that 
$b'\notin I$ and $\cof^{K}(d_{\gi,b})>\gg_b$, is unbounded in
$\vec b$.
\end{claim}
\begin{proof}
Suppose the contrary, that there is an $\gi<\gw_1$ and a $b_0\in\vec b$
such that $d_{\gi,b}$ falls into either case~(1) or case~(2a) for all
$b>b_0$. 
Then from the conclusions to these cases given in the table
we can conclude that $\cof(d_{\gi,b'})=\cof(d_{\gi,b})$ whenever
$b_0<b'\gslesseq b$.
By corollary~\ref{gs-unbnded}, $\set{b\in\vec b:\gs_b<b_0}$ is bounded
in $\vec b$, say by $b_1>b_0$.  Then for all $b>b_1$ in $\vec b$ there
is $b'\gsless b$ such that $b_0<b'<b_1$: namely the 
least member $b'$ of $\vec b$ such that $b_0<b'\gslesseq b$. Then
$\gs_{b'}<b_0$, so $b'<b_1$ by the choice of $b_1$.

\paranumber{22g}
It follows that $\cof(b)=\cof(b')<b_1$ for all but boundedly many
$b\in \vec b$, 
contradicting proposition~\ref{cof-unbnded-prop}.
\end{proof}
\begin{claim}
There is an unbounded subset $y$ of $\vec b\setminus I$ such that
$\gs_b\in y$ whenever $b\in y$ and $\gs_b\in\vec b$.
\end{claim}
\begin{proof}
If this claim is false then by the last claim there are, for every
$\gi<\gw_1$, ordinals $b_{0,\gi}\gsless b_{1,\gi}$ such that
$b_{0,\gi}\in I$ and $d_{\gi,b_{1,\gi}}$ falls into case~(2c).  Since
$\vec b$ is countable, there are ordinals $b_0\gsless b_1$ in $\vec b$ and an
uncountable set $x\subset \gw_1$ so that
$b_{0,\gi}=b_0$ and $b_{1,\gi}=b_1$ for all $\gi\in x$.
By the hypothesis of the theorem we have 
$\gg_{b_1}<b_1\le\gg_{b_1}^{+n}$ for some $n\in\gw$.  Since 
$\gg_{b_1}<\cof^K(d_{\gi,b_1})<b_1$
it follows that there are only 
finitely many possible  values for $\cof^K(d_{\gi,b_1})$, so there
must be ordinals $\gi<\gi'\in x$ such that
$\cof^K(d_{\gi,b_1})=\cof^K(d_{\gi',b_1})$.   An easy
induction, using lemma~\ref{t.dcases}, shows that
$\cof^K(d_{\gi,b'})=\cof^K(d_{\gi',b'})$ for all $b'\gsless b_1$, and  in 
particular $\cof^K(d_{\gi,b_0})=\cof^K(d_{\gi',b_0})$.  This is
impossible since it implies that 
\[
\cof(d_{\gi',b_0})=\cof(d_{\gi,b_0})\le
d_{\gi,b_0}<d_{\gi',b_0},
\]
 contradicting the fact that  $b_0$ is in $I$
and hence $d_{\gi',b_0}$ is regular.	
\end{proof}

\paranumber{22h}
\newcommand\propp{\operatorname{P}}
\theoremstyle{plain}
\newtheorem*{propPk}{Property $ \propp(k)$}
Since the set $y\subset\vec b$ is unbounded in $\vec b$, we have
$\tcf(\prod y)=\tcf(\prod\vec b)=\gl$.  Since $y$ is closed under the
operation $b\mapsto\gs_b$, the conclusion of lemma~\ref{cnt-seq} is
still true  with  $\vec b$ replaced by $y$.   If we  let $k$ be such that
$\gl=\gk^{+(k+1)}$ in $K$ then $y$ witnesses the truth of $\propp(k)$:
\begin{propPk}
There is a countable sequence $\vec b$ of regular cardinals of $K$, cofinal in
$\gk$, along with nondecreasing functions $f_b\in K$ and
ordinals~$\gg_b$, $\ga_b$, and $\gs_b$ for $b\in \vec b$ 
such that $\gg_b<b\le\gg_b^{+k}$, and  $\gs_b\in\vec b$ for
all but boundedly many $b\in\vec b$, and  
 and so that
the set $\cl$ defined below has true cofinality $\gk^{+(k+1)}$:

\paranumber{22i}
The sequence $\vec c\in\prod b$ is in $\cl$ if for some precovering
set $Y$, and all sufficiently large $b\in\vec b$,  
\[c_{b}=\indiscbelong(f_b(c_{\gs_b}), \gg_b,\ga_b).\]

\paranumber{22j}
Furthermore, any strictly increasing, non-cofinal subsequence of $\cl$
of cofinality greater than $\gk^+$ has a least upper bound in
$\prod\vec b$.
\end{propPk}
We will prove by induction that $\propp(m)$ is false for all $m>1$, 
contradicting the observation that  $\propp(k)$ is true and hence
finishing the proof of theorem~\ref{app-main-thm}.  
Notice that since lemma~\ref{t.dcases} and the claims above follow from
lemma~\ref{cnt-seq}, they are true for any witness $\vec b$ to
property $\propp(k)$ for any $k\ge 1$.  This makes it  easy to see
that 
$\propp(1)$ is false, since in 
that case we always have $\cof^K(d_{\gi,b})\le d_{\gi,b}<b<\gg_b^{+}$,
so that case~(2c) can never hold, contrary to the first claim above.  

\paranumber{22k}
Now we complete the proof by showing  that $\propp(m)$ implies
$\propp(m-1)$ for all $m>1$. 
Suppose that $\vec b$, $\vec \ga$, $\vec \gg$, and
$\vec f$ witness the truth of $\propp(m)$.  Then $\tcf(\prod\vec
b)=\gk^{+(m+1)}$, so we can 
let $\vec c$ be the least
upper bound of the first $\gk^{+m}$ sequences from $\cl$.  
Then
$c_b< b\le\gg_b^{+m}$ for all $b\in\vec b$, so
$\xi_b=\cof^K(c_b)<\gg_b^{+m}$.  Let
$b'_{b}=\gg_b+\xi_b\le\gg_{b}^{+(m-1)}$, and let $\tau_b\in K$ be a
continuous, unbounded and increasing map from $b'_b$ into $c_b$.  Set
$\tilde\tau_b=i_{\gg_b,\ga_b}(\tau_b)$ and set
$f'_b=\tilde\tau_{b}^{-1}\circ f_b\circ \tau_{\gs_b}$.  Then $\vec
b'$, $\vec \gg$, $\vec\ga$ and $\vec f'$ witness the truth of
$\propp(m-1)$, as required.

\paranumber{22l}
It follows by induction that $\propp(m)$ is false for all $m\ge1$,
contradicting 
$\propp(k)$.  This contradiction 
completes the proof of theorem~\ref{app-main-thm}(3), assuming
lemma~\ref{t.dcases}. 
\end{proof}

\paranumber{23}
In the rest of this subsection we finish the proof of
theorem~\ref{app-main-thm} by proving 
lemma~\ref{t.dcases}.  First we need a preliminary lemma:

\begin{lemma}\label{*cond}
If $\vec d\in\cd$ then for any \pcs\ $Y$ with $\vec d\in Y$, the
following equation holds for all but boundedly many $b\in \vec b$:
\begin{equation}
d_b=\cases
a^Y(f_b(d_{\gs_b}),\gg_b)\le\lis^Y(f_b(d_{\gs_b}),\gg_b)&\text{if $b\in I$}\\
\indiscbelong^Y(f_b(d_{\gs_b}),\gg_b,\ga_b)&\text{if $b\notin I$}.
\endcases\tag{$*$}
\end{equation}
\end{lemma}
\begin{proof}
Let $\vec d$ be any member of $\cd$, and let $Y$ be a \pcs\ with $\vec
d\in Y$.

First we will prove the inequality in the case $b\in I$.
Suppose to the contrary that there are unboundedly many $b\in I$ such that
$\xi_b=\lis^Y(f_b(d_{\gs_b}),\gg_b)<d_b$.  Then from the definition of
$\cd$ there is a sequence
$\vec c\in\cl$,
with $\vec c\lessb\vec d$, such  that $\xi_b<\vec c_{b}$ for all but
boundedly many of those $b\in I$ such that $\xi_b<d_b$.
But this is impossible, since then
\[
\xi_b<c_{b}=
\lis^Y(f_{b}(c_{\gs_b}),\gg_b)<\lis^Y(f_{b}(d_{\gs_b}),\gg_b)=\xi_b
\]
for each such $b$. This contradiction shows that
$d_b\le\lis^Y(f_b(d_{\gs_b}),\gg_b)$ for all but boundedly many
$b\in I$.

\smallskip
Now we prove the identity in both cases. 
Define
the sequence $\vec \xi$ by $\xi_b\eqdef a^Y(f_b(d_{\gs_b}),\gg_b)$ if $b\in
I$ and $\xi_b\eqdef\indiscbelong^Y(f_b(d_{\gs_b}),\gg_b,\ga_b)$ if
$b\notin I$. 
We need to show that $d_b=\xi_b$ for all but boundedly many $b\in\vec b$.  
We will show first that $d_b\le\xi_b$.

\paranumber{24}
If, to the contrary, $d_b>\xi_b$ for unboundedly many $b\in\vec b$
then there is a sequence 
$\vec c\in\cl\cap\prod\vec d$ such that $\xi_b\le c_b$ for
unboundedly many $b\in\vec b$.  
Then $\vec c\in Y$, since $\cl$ and $\vec b$ are in $Y$.  
Now $c_{\gs_b}<d_{\gs_b}$ implies
that $c_b=\lis^Y(f_b(c_{\gs_b}),\gg_b)<\lis^Y(f_b(d_{\gs_b}),\gg_b)=\xi_b$ for
all sufficiently large $b\in I$ and
$c_b=\indiscbelong^Y(f_b(c_{\gs_b}),\gg_b,\ga_b)<
\indiscbelong^Y(f_b(d_{\gs_b}),\gg_b,\ga_b)=\xi_b$ for every sufficiently
large $b\notin I$.  Since this contradicts the choice of $\vec c$ we
must have $\xi_b\ge d_b$ for almost every $b\in\vec b$.

\smallskip
\paranumber{25}
We now complete the proof of the lemma by showing that $d_b\ge\xi_b$
for all but boundedly many $b\in \vec b$.  Assume the contrary, that
$d_b<\xi_b$ for unboundedly many $b\in\vec b$.  We consider the cases
$b\in I$ and $b\notin I$ separately. 

Suppose first that
 $d_b<\xi_b=a^Y(f_b(d_{\gs_b}),\gg_b)$ for unboundedly many $b\in I$.
By the definition of an accumulation point it follows that for
unboundedly many $b\in I$  there is an 
ordinal $\zeta_b<f_b(d_{\gs_b})$ in $Y$ such that
$d_b\le\ell^Y(\zeta_b,\gg_b)$. 

We claim that  $f_b(\nu)<f_b(d_{\gs_b})$
for all sufficiently large $b\in I$ and all $\nu<d_{\gs_b}$.
Otherwise pick $\nu_b<d_{\gs_b}$ for unboundedly many $b\in\vec b$ so
that 
$f_b(\nu_b)\ge f_b(d_{\gs_b})$. Then
$S^{f,Y}_b(\nu_b)=\lis^Y(f_b(\nu_b),\gg_b)\ge d_b$, so that any member
of $\prod \vec d\cap\cl$ must be smaller than $\nu_b$ for all but boundedly
many $b$.  Since $\vec d=\lub(\prod d\cap\cl)$ it follows that
$\nu_b\ge d_b$.

Since  $f_b$ is continuous and $\xi_b<f_b(d_{\gs_b})$ 
it follows that there is
an ordinal  $\nu_{b}<d_{\gs_b}$
such that 
$\zeta_b<f_b(\nu_{b})<f_b(d_{\gs_b})$.  
Now pick $\vec c\in\cl\cap\prod d$
such that $c_{\gs_b}>\nu_{b}$ for all but boundedly many $b$ such that
$\nu_{b}$ is defined.  This is possible since $\cof(d_b)>\gw$ for
almost all $b\in\vec b$.  Then
\[
c_b=\ell^Y(f_b(c_{\gs_b}),\gg_b)\ge\ell^Y(f_b(\nu_{b}),\gg_b)
\ge\ell^Y(\zeta_b,\gg_b)\ge d_b
\] 
for all sufficiently large $b\in I$
such that $\nu_{b}$ is defined, contradicting the assumption that
$\vec c\in\prod\vec d$.  Thus $d_b=\xi_b$ for all but boundedly many
$b\in I$. 

\smallskip
\paranumber{26}
The argument for $b\notin I$ is similar.  If
$d_b<\xi_b=\indiscbelong^Y(f_b(d_{\gs_b}),\gg_b,\ga_b)$ for unboundedly 
many $b\notin I$ then for
unboundedly many $b\notin I$ there is an ordinal
$\nu_b<d_{\gs_b}$  such that 
$i^Y_{\gg_b,\ga_b}(d_b)<f_b(\nu_{b})<f(d_{\gs_b})$.  Choose  $\vec 
c\in\cl\cap\prod\vec d$ so that
$c_{\gs_b}>\nu_b$ whenever $\nu_b$ is defined.  Then
\[
c_{b}= \indiscbelong^Y(f_b(c_{\gs_b}),\gg_b,\ga_b)>
       \indiscbelong^Y(f_b(\nu_b),\gg_b,\ga_b) >d_b
\]
 for every sufficiently large $b$ such that $\nu_{b}$ is
defined, contradicting the assumption that $\vec c\in\prod \vec d$.

\paranumber{27}
This completes the proof that $\vec \xi=\vec d$, and hence of the
lemma. 
\end{proof}

The  next four lemmas correspond to the four cases in
lemma~\ref{t.dcases}.  The first is, by a wide margin, the most
difficult.

\begin{lemma}[lemma~\ref{t.dcases}, case~1]
Every $b\in I$ is regular in $K$, and
$\cof(d_{\gi,b})=\cof(d_{\gi,\gs_b})$ for almost all $(\gi,b)$ with
$b\in I$.
\end{lemma}

\paranumber{28}
\begin{proof}
Recall that every member of $I$ is a limit of principle indiscernibles
for $\ga_b$, and hence is either a principle indiscernible for $\ga_b$
or equal to the measurable cardinal $\ga_b$ of $K$.  In either case, $b$ is
regular in $K$.

Pick a $\gw_1$-closed \pcs\ $Y\prec H_{(2^\tau)^+}$ such that
everything relevant, including $H_\tau$ and  $\seq{\vec
d_\gi:\gi<\gw_1}$, is in $Y$.  Next pick, inside $Y$, a \pcs\ $Y_b$ for
each $b\in\vec b$ so that $Y_b\prec H_{\tau}$, $\gs_b\subset Y_b$, and
$Y_b$ contains all of the sequences which have been defined.  There
exists  \pcs{s} $Y'$ with $\gs_b\subset Y'$ since $\gk$ is a strong
limit 
cardinal and hence $\gg_b^{\gd}<\gk$, and we can find the sequence
$\seq{Y_b:b\in I}$ inside $Y$
since we have  strengthened the
usual requirement of $Y\prec H_\tau$ to $Y\prec
H_{(2^\tau)^+}$.

By lemma~\ref{*cond} there is, for each
$\gi<\gw_1$, an ordinal $\nu_\gi<\gk$ such that 
$d_{\gi,b}=a^Y(f_b(d_{\gi,\gs_b}),\gg_b)$ if $b\in I\setminus \nu_\gi$ and
$d_{\gi,b}=\indiscbelong^Y (f_b(d_{\gi,\gs_b}),\gg_b,\ga_b)$ if $b\notin
(I\cup \nu_\gi)$.  Since $\cof(\gk)=\gw$ there is a fixed $\nu$ such
that 
we can take $\nu_\gi=\nu$ for uncountably many $\gi<\gw_1$.
 By restricting ourselves to this
uncountable subset and removing $\vec b\cap\nu$ from $\vec b$ 
we can assume \iwlog\ that $d_{\gi,b}=a^Y(f_b(d_{\gi,\gs_b}),\gg_b)$ or
$d_{\gi,b}=\indiscbelong^Y(f_b(d_{\gi,\gs_b}),\gg_b,\ga_b)$  whenever
$\gi<\gw_1$ and $\gs_b\in\vec b$.

\paranumber{29}
Define $d_b\eqdef\sup_\gi (d_{\gi,b})<b$ for each $b\in I$, so that
$\vec d=\lub\set{d_{\gi}:\gi<\gw_1}$.
Since 
$\vec d_\gi$ satisfies condition~$(*)$ of lemma~\ref{*cond}
for each $\gi<\gw_1$, the sequence $\vec d$ must also satisfy
condition~($*$).
  Since $Y$ is $\gw_1$-closed it follows by
lemma~\ref{semi-cont} that 
$d_b=a^Y(f_b(d_{\gs_b}),\gg_b)=\lis^Y(f_b(d_{\gs_b}),\gg_b)$. In
particular, if $b\in I$ then $\seq{d_{\gi,b}:\gi<\gw_1}$ is a \pris\
for the constant sequence $d_b$

We will find functions 
$g_b\colon d_{\gs_b}\to O(\ga_b)$ in $K$ so that
$d_{\gi,b}=a^Y(g_b(d_{\gi,\gs_b}),\gg_b)$ for $\gi<\gw_1$.  In addition
the functions $g_b$ will be continuous, nondecreasing, and
$\range(g_b)$ will be cofinal in $g_b(d_{\gi,\gs_b})$. 

\paranumber{30}
To see that this implies the lemma, notice that the properties above
imply that 
$\set{\lis^{Y_b}(g_b(\nu),\gg_b):\nu<d_{\gi,\gs_b}}$ 
is cofinal in
$d_{\gi,b}\cap Y_b$ for all but countably many $\gi<\gw_1$.  Thus it
will be sufficient to show that $Y_b$ is cofinal in $d_{\gi,b}$, for
all but boundedly many $b\in\vec b$. Suppose to the contrary that $Y_b$ is
bounded in $d_{\gi,b}$ for unboundedly many $b\in\vec b$.  Since everything
under consideration, including $\seq{Y_b:b\in\vec b}$, is in $Y$, the
upper bounds $\xi_b=\sup (Y_b\cap d_{\gi,b})$ are in $Y$.
Since $d_{\gi,b}= a^Y(g_b(d_{\gi,\gs_b}),\gg_b)$, it
follows that there is $\nu_b<d_{\gi,\gs_b}$ such that
$\lis^Y(g_b(\nu_b),\gg_{b})\ge\xi_b$.  But $\nu_b\in\gs_b\subset Y_b$,
and hence $\lis^{Y_b}(g_b(\nu_b),\gg_b)\in Y_b$.
Thus
$\xi_b<\lis^{Y_b}(g_b(\nu_b),\gg_b)\le\sup(Y_b\cap d_{\gi,\gs_b})$,
contrary to the choice of $\xi_b$.  This contradiction shows that
$Y_b$ is cofinal in $d_{\gi,b}$ and hence
$\cof(d_{\gi,b})=\cof(d_{\gi,\gs_b})$. 

\paranumber{31}
The functions $f_b$ are continuous and increasing, and since
$\cl^f(\vec b)$ is cofinal in $\vec d$ the range of $f_b$
is cofinal in $f_{b}(d_{\gi,\gs_b})$  for almost all $b$, for each
$\gi<\gw_1$.  
As our first approximation to $g_b$, 
define $g_b^*\colon \gg_b\to \SO(d_{b})$ by letting
$g_b^*(\nu)$ be the least ordinal $\nu'$ such that
$i_{d_b,\ga_b}(\nu')\ge f_b(\nu)$.  Then $g^*_b$ is continuous and
nondecreasing since $f_b$ is, and the range of $g^*_b$ is cofinal in
$g^*_b(d_{\gi,\gs_b})\cap Y$ since
$d_{\gi,b}=a^Y(f_b(d_{\gi,\gs_b}),\gg_1)$ implies that the range of 
$i_{d_b,\ga_b}$ is cofinal in $f_b(d_{\gi,b})\cap Y$.  
\smallskip
There are two problems with $g^*_b$: first, it is not in  either $K$
or in $Y$, and second, it is cofinal in $g^*_b(d_{\gi,\gs_b})\cap Y$, not in
$g^*_b(d_{\gi,\gs_b})$.   
We will attack the second problem by going back to the proof of the
covering lemma, working in the preimage of the collapse map $\pi$.  

\paranumber{32}
Define $\bar g^*_b$  by letting $\bar g_b^*(\nu)$ be the
least  ordinal $\nu'$ such that $\bari_{\bar
d_b,\bar\ga_b}(\nu')\ge\pi^{-1}(f_b)(\nu)$.  Then $g^*_b$, as a set of
ordered pairs, is equal to $\pi\image\bar g^*_{b}$.  Now $\bar g^*_b$ is
defined from the iterated ultrapower $b_{\bar
d_b,\bar\ga_b}\colon\barmse m_{\bar d}\to \barmse m_{\bar\ga_b}$, but
it only requires a finite part of the iterated ultrapower: the initial
ultrapower by $\bare_{\bar d}$
together with the support of $\pi^{-1}(f_b)$.  Thus 
$\bar g^*_{b}$ can be defined inside $\barmse m_{\bar d_b}$.
Now define  
\[\pi_{d_b}\colon \barmse m_{\bar d_b}\longrightarrow\mse
m_{d_b}\eqdef\ult\bigl(\barmse m_{\bar d_b},\pi,\sup(\pi\image\len(\bare_{\bar
d_b})\bigr).
\]
Then $\pi_{d_b}(\bar g^*_b)$ is the desired extension of
$g^*_b$.
Unfortunately there is no reason to believe $\mse m_{d}$ is in $K$, so
we don't know that $\pi_{d_b}(\bar g_b^*)$ is in either $Y$ or $K$.
However $\pi_{d_b}(\bar g_b^*)$ is cofinal in $g^*(d_{\gi,\gs_b})$, so 
$\cof(g^*(d_{\gi,\gs_b}))<\gg_b$ in
$V$ and hence by elementarity in $Y$.

\paranumber{33}
We now proceed as in the proof of lemma~\ref{ftoi-lem}(5).
Since $\cof^Y(g^*_b(d_{\gi,\gs_b}))<\gg_b$ the covering lemma,
 applied in $Y$, implies that
there is a function $k\colon d_b\to \SO(d_b)$ in $K\cap Y$ such that
$\range (k)$ 
is closed and  is cofinal in each of the ordinals
$g^*(d_{\gi,\gs_b})$. for $\gi<\gw_1$.
Then $\bar k=\pi^{-1}(k)$ is in
$\bark$ and hence is in  $\barmse m_{\bar d}$.
Thus we can define a function $\bar s$ in $\barmse m_{\bar d}$
by letting $\bar s(\nu)$ be the least ordinal $\nu'$ such that $\bar
k(\nu')\ge \bar g^*(\nu)$. Then $\bar s\in \bark$,
since $\bark$ and $\barmse m_{\bar d}$ contain the same subsets of
$\bar d$, so the function $g_b=k\circ \pi(\bar s)$ is in $K$.  This
function $g_b$ has the required properties, and this
 completes the proof of case~1 of lemma~\ref{t.dcases}.
\end{proof}

\begin{lemma}
[lemma~\ref{t.dcases}, case~2a]
\label{less-gg}
 For almost all $(\gi,b)$ such that $b\notin I$ and
$\cof^K(d_{\gi,b})<\gg_b$ we have
$\cof^K(d_{\gi,b})=\cof^K(d_{\gi,\gs_b})$.
\end{lemma}

\paranumber{34}
\begin{proof}
This lemma, as well as the next two, depend of the following
calculation.  Each of the identities holds for almost all pairs
$(\gi,b)$ which satisfy the hypothesis of this lemma.
\begin{align*}
\pi^{-1}(\cof^K(d_{\gi,b}))
 &=\cof^{\bark}(\bar d_{\gi,b})\\
 &=\cof^{\barmse m_{\bar\gg_b}}(\bar d_{\gi,b})
    &&{\parbox[t]{2in}{since $\cof^K(d_{\gi,b})\le\gg_b$  and\\
      $\ps^{\barmse
           m_{\bar\gg_b}}(\bar\gg_b)=\ps^{\bark}(\bar\gg_b)$}}\tag{i}\\
 &=\bari_{\bar\gg_b,\bar\ga_b}(\cof^{\barmse m_{\bar\gg_b}}(\bar d_{\gi,b}))
    &&\text{since $\cof^K(d_{\gi,b})<\gg_b$}\tag{ii}\\
 &= \cof^{\barmse m_{\bar\ga_b}}(\bari_{\bar\gg_b,\bar \ga_b}(\bar
          d_{\gi,b}))\\ 
 &=\cof^{\barmse m_{\bar\ga_b}}(\pi^{-1}(f_b(d_{\gi,\gs_b})))
    &&\text{since $d_{\gi,b}=\gb^Y(f_b(d_{\gi,\gs_b}),\gg_b,\ga_b)$}\\
 &= \cof^{\bark}(\pi^{-1}(f_b(d_{\gi,\gs_b})))
    &&\text{since $\cof(\bar d_{i,\gs_b})<\bar\ga_b$}\\
 & =\pi^{-1}(\cof^K(f_b(d_{\gi,\gs_b})))\\
\intertext{so}
\cof^{K}(d_{\gi,b})&=\cof^K(f_b(d_{\gi,\gs_b}))\\
  &=\cof^K(d_{\gi,\gs_b})
   &&\text{since $\range(f_b)$ is cofinal in $f_b(d_{\gi,\gs_b})$.}
\end{align*}
\end{proof}

\paranumber{35}
\begin{lemma}
[lemma~\ref{t.dcases}, case~2b]
 $\cof^K(d_{\gi,b})\not=\gg_b$ for almost all $(\gi,b)$ such that
$b\notin I$. 
\end{lemma}
\begin{proof}
All of the first sequence of equalities in the proof of
lemma~\ref{less-gg} still hold in this case except for line~(ii).  In
this case we get 
\(
\bari_{\bar\gg_b,\bar\ga_b}(\cof^{\barmse m_{\bar\gg_b}}(\bar d_{\gi,b}))=
\bari_{\bar\gg_b,\bar\ga_b}(\bar\gg_b)=\bar\ga_b
\).  
The rest of the equalities in this sequence still hold, so 
\[
\cof^K(f_b(d_{\gi,\gs_b}))=\bari_{\bar\gg_b,\bar\ga_b}(\cof^{\barmse
m_{\bar\gg_b}}(\bar d_{\gi,b}))=
\ga_b,
\] 
but this is impossible since 
$\cof^K(f_b(d_{\gi,\gs_b}))=\cof^K(d_{\gi,\gs_b})<\ga_b$.
\end{proof}

\begin{lemma}
[lemma~\ref{t.dcases}, case~2c]
There is an $\gi_0<\gw_1$ such that for all $\gi,\gi'>\gi_0$, for all
but boundedly many $b\in\vec b$, if
$\cof^K(d_{\gi,b})=\cof^{K}(d_{\gi',b})$ 
then $\cof^K(d_{\gi,\gs_b})=\cof^K(d_{\gi',\gs_b})$.
\end{lemma}
\begin{proof}
Again, consider the sequence of equalities from the proof of 
lemma~\ref{less-gg}.  In this case, lines~(i) and~(ii) both fail.
Since $\cof^{\bark}(\bar d_{\gi,b})=\cof^{\bark}(\bar
d_{\gi',b})$ and $\barmse m_{\bar\gg_b}$ is larger than $\bark$ the
argument for line~(i) shows that
\(
\cof^{\barmse m_{\gg_b}}(\bar d_{\gi,b})=\cof^{\barmse
m_{\gg_b}}(\bar d_{\gi',b})
\).  
Then the argument for line~(ii)  gives
\(
\bari_{\bar\gg_b,\bar\ga_b}(\cof^{\barmse m_{\bar\gg_b}}(\bar
d_{\gi,b}))
=\bari_{\bar\gg_b,\bar\ga_b}(\cof^{\barmse m_{\bar\gg_b}}(\bar
d_{\gi',b}))
\).  
The rest of the identities remain valid, so that 
$\cof^K(f_b(d_{\gi,\gs_b}))=\cof^K(f_b(d_{\gi',\gs_b}))$ and hence 
$\cof^K(d_{\gi,\gs_b})=\cof^K(d_{\gi',\gs_b})$.
\end{proof}

This completes the proof of theorem \ref{app-main-thm}(3).\qed

\subsection{Further results}
In this subsection we extend the results of the two previous
subsections.
The
first result, theorem~\ref{amt1}, completes the proof of
theorem~\ref{app-main-thm} by  strengthening the conclusion from $o(\gk)=2^{\gk}$ to
$o(\gk)=2^{\gk}+\cof(\gk)$ in the case $\cof(\gk)>\gw_1$.   The
second, theorem~\ref{gch-main-cor}
shows that if $\cof(\gk)=\gw$ 
then we can strengthen the conclusion from $o(\gk)=2^{\gk}$ to
$o(\gk)=2^{\gk}+1$ if either $\gk<\aleph_{\gk}$ or the GCH holds below $\gk$. 
\begin{lemma}
\label{eta_b=gl}
Suppose that $\gk$ is a strong limit cardinal with
$\cof(\gk)=\gd<\gk$, and that $2^{\gk}=\gl>\gk^+$ where $\gl$ is
regular and if $\gl$ is a
successor cardinal then the predecessor of $\gl$ has cofinality
greater then $\gk$.
If $\cof(\gk)=\gw$ then also assume that there is $m<\gw$ such that
$\set{\ga<\gk:o(\ga)=\ga^{+m}}$ is bounded in $\gk$.

Then either $\gk$ is a limit of accumulation points for $\gl$, or
$\gk$ is a limit of indiscernibles for extenders $\ce_\gg$ on $\gk$
with $\gg\ge\gl$.
\end{lemma}

\paranumber{36}
\begin{proof}
We claim that $\eta_b\ge\gl$ for all but boundedly many $b\in\vec b$.
Suppose the contrary. Since $\cof(\gl)>\gk$ it follows that there is
$\eta<\gl$ such that 
$\eta_b<\eta$ for all but boundedly many $b\in\vec b$.  Now the
function $f_b$ used to define $\cl^f$ had range contained in $\eta_b$,
so we can restrict ourselves to functions $f$ with
$\range(f)\subset\eta$.  There are $(\eta^{\gk})^K<\gl$ many such
functions. 

\paranumber{37}
The only use of the hypothesis $(o(\gk)^{\gk})^K<\gl$ in the proof of
theorem~\ref{app-main-thm} came in the proof of
corollary~\ref{tcfb=gl}, where this hypothesis was used to show that
there is a single function $f$ such that $\cl^{f}$ is cofinal in
$\prod\vec b$.  The reason was that there were only
$(o(\gk)^{\gk})^K<\gl$ 
relevant functions $f$, while $\tcf(\prod(\vec b))=\gl$ is greater
than  $(o(\gk)^{\gk})^K$.  Thus the conclusion of 
corollary~\ref{tcfb=gl} is true under
our assumption that $\eta_b<\gl$ for cofinally many $b\in\vec b$.  In
the rest of the proof of theorem~\ref{app-main-thm} we showed that
lemma~\ref{Ufclf-cofinal} leads to a contradiction.  Hence the falsity
of our current claim would lead to the same contradiction, and the
claim must be true.

We now consider two cases.  We have $\ga_b=\gk$ for
cofinally many $b\in\vec b$.  
If $b\in\vec b$ has $\ga_b=\gk$ and is a limit of principle
indiscernibles then $b$ is an accumulation point for $\eta_b$.  
If $b\in\vec b$ has $\ga_b=\gk$ and is not a limit of principle
indiscernible then 
$\gg_b$ is a
principle indiscernible for some $\eta'$ with
$\len(\ce_{\eta'})\ge\eta_b$ so that $\eta'>\eta_b\ge\gl$.
One of these cases must hold for cofinally many $b\in\vec b$, and the
lemma follows.
\end{proof}

\begin{theorem}
[theorem~\ref{app-main-thm}(1)]
\label{amt1}
Suppose that $\gk$ is a strong limit cardinal with 
$\gw_1<\gd=\cof(\gk)<\gk$, and that $2^\gk\ge\gl>\gk^+$, where if $\gl$ is a
successor cardinal then the predecessor of $\gl$ has cofinality greater
than $\gk$.   Then $o(\gk)\ge\gl+\gd$.

If there is an $n<\gw$ such that $\set{\ga<\gk:o(\ga)\ge\ga^{+n}}$ is
bounded in $\gk$ then the result is also true for
$\gd=\cof(\gk)=\gw_1$. 
\end{theorem}

\paranumber{38}
\begin{proof}
Let $\vec d$ be given by lemma~\ref{eta_b=gl}, so that every member of
$\vec d$ is either an accumulation point for $\gl$ or a principle
indiscernible 
for some $\eta\ge\gl$.  Then every uncountable limit point of $\vec d$
of uncountable cofinality is an accumulation point for $\gl$ and hence,
by lemma~\ref{semi-cont}, is a principle indiscernible for some
$\eta\ge\gl$.  Continuing by induction, any ordinal which is a limit
of
$\gw_1^{\ga+1}$ members of $\vec d$  is a principle indiscernible for some
$\eta\ge\gl+\ga$.  Thus $o(\gk)\ge\gl+\gd$.
\end{proof}

\paranumber{39}
In view of Silver's
fundamental result in
\cite{Silver.SCH} the next
observation is only of interest when $\cof(\gk)=\gw$.  As usual, all
successors are calculated in $K$ unless indicated otherwise.
\begin{theorem}
\label{gch-main-cor}
Suppose that $\gk$ is a strong limit cardinal of cofinality $\gw$ and 
there is a $k<\gw$ so that the set of $\nu<\gk$ such that
$o(\nu)>\nu^{+k}$ is bounded in $\gk$.
Suppose further that $o(\gk)=2^{\gk}>(\gk^{++})^V$.  Then
(i)~$\gk=\aleph_{\gk}$ and~(ii)
if $2^\gk =({\gk^{+m}})^V$ then 
$2^{\nu}\ge({\nu^{+(m-1)}})^V$ for cofinally many $\nu<\gk$.

\end{theorem}
\begin{proof}
Set $\gl=2^{\gk}$ and let $n\ge m$ where
$\gl=(2^{+m})^V=(2^{+n})^K$.   
Since by hypothesis $o(\gk)=\gl=2^{\gk}$, lemma~\ref{eta_b=gl} implies
that there is a cofinal sequence $\vec b=\seq{b_i:i\in\gw}$ of
accumulation points for $\gl$.  We can pick $\vec b$ so that for each
$i<\gw$ 
there is $\gg_i<b_i$ so that $b_i=a^Y(\gl,\gg_i)$ for any \pcs\ $Y$. 

\paranumber{40}
If $\gb<\gl$ then define $\vec d_{\gb}$ by $d_i=s^Y(\gb,\gg_i)$ for
any \pcs\ $Y$ with $\gb\in Y$.  Then $\gb\le\gb'$ implies $\vec
d_\gb\leq\vec d_{\gb'}$, so $\set{\vec d_{\gb}:\gb<\gl}$ witnesses
that $\tcf(\prod\vec\gb)=\gl$.

Define ordinals $c_{k,i}$ for each $k,i\in\gw$ by
recursion on $k$, setting
$c_{0,i}=\gg_i$ and $c_{k+1,i}=\lis^X(\gk^{+(n-1)},c_{k,i})$.  Now
define 
the sequence $\vec c=\sup_i\vec c_i$, that is, $c_i=\sup_{k<\gw}c_{k,i}$.
Then $\vec
c$ is an \is\ for $\gk$ and for each $i<\gw$ the sequence
$\vec c^*_i=\seq{c_{k,i}:k\in\gw}$ is an
\is\ for $c_i$.
 Lemma~\ref{semi-cont} implies that $c_i=s^X(\gb_i,a_i)$ for some
$\gb_i$ with $\gk^{+(n-1)}<\gb_i<\gl$.  In particular
$i_{c_i,\gk}(o(c_i))\ge\gk^{+(n-1)}$ and hence, using
lemma~\ref{semi-cont},
$o(c_i)\ge c_i^{+(n-1)}$
for all sufficiently large $i<\gw$.

\paranumber{41}
For each $i<\gw$ and $\gb<c_i^{+(n-1)}$ define $\vec d_{i,\gb}$ by
$d_{i,\gb,k}=\lis^{Y}(\gb,c_{i,k})$ for an appropriate \pcs\ $Y$.   Then
$\vec d_{i,\gb}\in\prod\vec c_i$, and for each $\gb<c_i^{+(n-1)}$
there is $\gb'$ such that $\gb<\gb'<c_i^{+(n-1)}$ so that $\vec
d_{i,\gb}\lessb\vec d_{i,\gb'}$.  It follows that there are
$c_i^{+(n-1)}$  distinct sequences  $d_{i,\gb}$, and hence $2^{c_i}\ge
c_i^{+(n-1)}$. 

\paranumber{42}
As usual, the cardinal $c_i^{+(n-1)}$ is computed in $K$.  
\begin{claim}
For almost all  $i<\gw$ 
\begin{equation}
\left|(c_i^{+(n-1)})^K\right|^{V}\ge\left(c_{i}^{+(m-2)}\right)^V.
\tag{$*$}
\label{e.gchthm}
\end{equation}
Furthermore, if equality holds then $\gw<\cof^V(c_i^{+})=\card{c_i}^{V}<c_i$.
\end{claim}
\begin{proof}
First note that if $0<s< n$ and $\gk^{+s}$ is a cardinal in $V$
then
$\prod_{i<\gw}c_i^{+s}$ has true cofinality $\gk^{+s}$.  Suppose that 
$0<s<s'<n$ and
$\gk^{+s}$ and $\gk^{+s'}$ are both cardinals in $V$.
 Then 
\[\tcf(\prod_{\gi}c_i^{+s})=\gk^{+s}\not=\gk^{+s'}=\tcf(\prod_{i}c_i^{+s'})
\]
and it follows that  $\cof(c_i^{+s})\not=\cof(c_i^{+s})$
for all but finitely many integers $i$.
But if $\cof(c_i^{+s})\le c_i$  then since $c_i$ is singular the
covering lemma, 
lemma~\ref{cl},
implies that $\cof(c_i^{+s})=\card{c_i}<c_i$.  If
$\card{(c_i^{+(n-1)})^K}^{V}<(c_{i}^{+(m-2)})^V$ then there are at
most $m-2$ distinct cofinalities available, out of the minimum $m-1$
needed,
for $\set{c_i^{+s}:0<s<n}$.
This contradiction proves the inequality of the claim.
Furthermore it shows that if the equality holds then $\card{c_i}<c_i$.
Since $\gk$ is a limit cardinal, $c_i>\gw_1$ for all sufficiently
large $i<\gw$ and it follows that
$\gw<\cof^V(c_i^{+})=\card{c_i}^{V}<c_i$, as claimed.
\end{proof}

To prove clause~(i) of the theorem, 
suppose to the contrary that $\gk<\aleph_{\gk}$, and let $\tau<\gk$
so that $\gk=\aleph_{\tau}$. Then there are only $\tau^{\gw}<\gk$ many
countable sequences of cardinals below $\gk$.  Since $\tcf\prod\vec
b=\gl>\gk$ it follows 
that there is $\vec \gg\lessb\vec b$ so that $b_i\le\gg_i^{+}$ for
cofinally many $i<\gw$.  We can modify the definition of the sequence 
$\vec c_i$, if necessary, so that $\vec \gg\lessb \vec c$.  
Since $\gl>(\gk^{++})^V$, the claim implies that there is some $s<\gw$ 
such that $c_i^{+s}$ is a cardinal in $V$ for infinitely many $i<\gw$,
and 
hence $c_i^{+s}\ge b_i$.  This is impossible, since there is a
sequence $\vec d$ of principle indiscernibles such that $\vec c<\vec
d<\vec b$, and every principle indiscernible is a limit cardinal of $K$.

To prove clause~(ii) of the conclusion,  notice first that if the
inequality~\eqref{e.gchthm} is strict then
$2^{c_i}>(c_{i}^{+(m-2)})^V$, so that the conclusion is true for
$\nu=c_i$ for almost all $i<\gw$.  If, on the other hand, equality
holds in~\eqref{e.gchthm} then
set $\xi=\card{c_i}<c_i$.  Then 
$\xi^{\gw}=c_i^{\gw}=(c_i^{+(n-1)})^K\ge({c_i^{+(m-2)}})^V$.  Since 
$\xi=\card{c_i}=\cof(c_i^+)$ is regular, there is $\nu<\xi$ such that
$\nu^{\gw}=\xi^\gw\ge(\xi^{+(m-2)})^V\ge({\nu^{+(m-1)}})^V$. 
\end{proof}

\paranumber{43}
We used the strong version of the weak covering lemma,~\ref{cl},
 which uses \pcs{}s which are not $\gw$-closed, to get that
$\cof(c_i^{+s})=\card{c_i}$ whenever 
$\card{c_i}>\gw_1$.  At the cost of some extra calculation it is
possible to use the weaker version of lemma~\ref{cl} which is refered
to in the remark following the statement of the lemma.  This version
implies that $(\cof(c_i^{+s}))^{\gw}\ge\card{c_i}$.

\paranumber{45}
\newcommand\awo{\aleph_{\gw_1}}
\newcommand\aw{\aleph_{\gw}}
\newcommand\pcf{\operatorname{pcf}} 
The next theorem is somewhat different but uses some of the ideas of
theorem~\ref{app-main-thm}. 
\begin{theorem} \label{aleph-gw1}
If 
$2^{\gw}<\aleph_{\gw}$  and  $2^{\aleph_\gw}>\aleph_{\gw_1}$
then there is a sharp for a model with a strong cardinal.
\end{theorem}
\begin{proof}
The proof depends on the following results of Shelah.  The definitions
may be found in \cite{shelah.card-arith}.
\begin{theorem}{\rm (Shelah, \cite{shelah.card-arith})}\label{shelah2}
\begin{enumerate}
\item
$\pcf\seq{\gw_n:n<\gw}=\set{\gk\le (\aleph_\gw)^{\gw}:\gk\text{ is
regular}}$. 
\item
Assume that $\vec a$ is a set of regular cardinals such that $2^{\card{
\vec a}}<\min(\vec a)$.  Then for every $\vec d\subset\pcf(\vec a)$ and
every $\mu\in \vec d$ there is a set $\vec d'\subset \vec d$ such that
$\card{\vec d'}\le \card{\vec a}$ and $\mu\in\pcf(\vec d')$.
\end{enumerate}
\end{theorem}
Let $A$ be the set of cardinals $\gd^{+}$ of $K$ below
$\aleph_{\gw_1}$ such that either $o(\ga)<\gd$ for all $\ga\le\gd$
or else $\gd$ is larger than every measurable cardinal of $K$ smaller
than $\aleph_{\gw_1}$.  The set $A$ is unbounded in $\awo$ since there
are no overlapping extenders in $K$.

\paranumber{46}
We claim that if $B\subset A$ with $\card B<\inf B$ then $\pcf(\prod
B)\le(\sup B)^{+}$.  To see this, let $\gk=\sup B$ and define, in $K$,
functions $a_f\in\prod A$ for each function $f\colon \gk\to\gk$ in $K$ by
setting, for $\nu={\gd^+}^K$ in $A$, $a_f(\nu)=\sup
(f\image\gd)\cap\nu$.
We will show that $\set{a_f\restrict B:f\in K}$ is
cofinal in $\prod B$.  If there is a largest measurable cardinal in
$K$ below
$\awo^V$ then this 
follows from lemma~\ref{cl}, the weak covering lemma.  
Otherwise if $\vec b\in \prod B$ then use the covering lemma,
together with the fact that 
proposition~\ref{inxE-bnd} implies  that $\nu$ cannot be
an indiscernible
since $o(\ga)^+<\nu$ for $\nu\in A$, to show
that there is a function $f\in K$ such that $b_\nu\in
f\image\gd$ whenever $\nu={\gd^{+}}^{(K)}$ is in $B$. Thus $\vec b
\lessb a_f$.

\paranumber{47}
Now let $A'=\set{\card\nu:\nu\in A}$.  Then $A'$ is unbounded in
$\awo$ and it follows by theorem~\ref{shelah2} that there is a
countable subset $B'$ of $A'$ such that $\awo\in\pcf(B')$.  Let
$B\subset A$ so that $B'=\set{\card\nu:\nu\in B}$.  Then for each
$\nu\in B$ the weak
covering lemma implies that $\cof(\nu)=\card\nu$, so that $\pcf(\prod
B')=\pcf(\prod B)$ and hence $(\sup B)^+<\awo\in\pcf(\prod B)$.  The
contradiction completes the proof of the theorem. 
\end{proof}

%

%

\paranumber{1}

\section{Open Problems}
There are a number of open problems which are related to results in
this paper.  The
most obvious questions concern the situation when $\gk$ has cofinality
$\gw$.  The most 
general question is whether the definability and uniqueness of
indiscernible sequences break down at $\gk^{\gw}$ for cardinals $\gk$
of cofinality $\gw$.  Since the first version of this paper, Gitik
\cite{gitik.hidden} has given a negative answer to this first question:
\begin{quest}
Is it still true if $\cof(\gk)=\gw$ that the notion of being an \is\ in $X$
for the constant sequence $\gk$ belonging to a sequence $\vec\gb$ is
independent of the \pcs\ $X$?
\end{quest}

The application concerning the singular cardinal
hypothesis may still be true, however.  
Since $o(\gk)=\gk^{+\gw}$ is enough to
give $2^{\gk}=\gk^{+(\gw+1)}$ the simplest unknown cases are the following:
\begin{quest}
If  $\gk$ is a strong limit cardinal with $2^\gk\ge\gk^{+(\gw+2)}$ then
must there be an inner model of $\exists\gk\,o(\gk)\ge\gk^{+(\gw+2)}$?
If $2^{\gk}=\gk^{+(\gw+1)}$ then must there be an inner model of 
$o(\gk)=\gk^{+\gw}$?
\end{quest}

\paranumber{2}
\begin{quest}
What is the exact consistency strength of $\cof(\gk)=\gw_1$ and
$2^{\gk}=\gl$ for regular $\gl>\gk^+$?
\end{quest}
By  theorem~\ref{app-main-thm} together with results of Woodin
\cite{cummings.GCH-successors} the answer lies
between $o(\gk)=\gl$ and $o(\gk)=\gl+\gw_1$.

\medskip
\paranumber{3}
A second problem concerns our use of $\gd$-closed \pcs s $X$.  In
Dodd and Jensen's work this assumption was weakened to $\gw_1\subset X$.
In \cite{mitchell-schimmerling.unclosecl} these methods have been
extended to the core models used in this paper,  but we do not
see how to avoid the use of $\gd$-closed \pcs s for 
the Gitik games in the proof of
lemma~\ref{aux-main-lem}.  The following can be
regarded as a test question.
\begin{quest}
Suppose that $\gk$ is singular, $2^{\gk}=\gl>\gk^{++}$ and
$2^{\ga}\le\gk^+$ for $\ga<\gk$.  Does it follow that there is an
inner model with $o(\gk)\ge\gk^{++}$?
\end{quest}

\paranumber{4}
The final question concerns what happens when there exist
overlapping extenders.  We give two possible test questions.

\begin{quest}
Suppose that $2^\gw<\aleph_\gw$ and $2^{\aleph_\gw}>\aleph_{\gw_1}$.
Does it follow that there is an inner model with a Woodin cardinal?
\end{quest}
\begin{quest}
Suppose that there is no model with a Woodin cardinal and that the Steel
core model \cite{steel.core-model} exists.  If $\gk$ is a singular
strong limit cardinal of uncountable cofinality such that
$2^{\gk}=\gl$ does it follow that $o(\gk)^\gk\ge\gl$ in $K$?
\end{quest}

%


\bibliographystyle{plain} 
\bibliography{new,logic}
\bigskip
\noindent
Department of Mathematics, Tel Aviv University, Tel-Aviv, Israel,\newline
\emph{gitik@math.tau.ac.il}

\medskip\noindent
Department of Mathematics, 
The University of Florida, Gainesville Florida,\newline
\emph{mitchell@math.ufl.edu}

\end{document}
\ifnotate
\closeout\notatefile
\newpage
\centerline{Index of notation}

\newcommand\notation#1#2#3#4{\line{#1\dotfill#3.#4\qquad#2}}
\renewcommand\usage{\notation}
\input \jobname.notation

\input x-notes.tex	
\fi
\end{document}